\documentclass[a4paper,fleqn]{cas-sc}
\usepackage{subcaption}
\usepackage[percent]{overpic}
\usepackage{bm}
\usepackage{mathtools}
\usepackage{tikz}
\usepackage{algorithmic}
\usepackage{algorithm}
\usepackage[utf8]{inputenc}
\usepackage{hyperref}
\usepackage{lineno}

\usepackage[numbers]{natbib}

\graphicspath{{figure/}} 

\renewcommand{\Re}{\text{Re}}

\begin{document}
\let\WriteBookmarks\relax
\def\floatpagepagefraction{1}
\def\textpagefraction{.001}

\title[mode = title]{Modeling and optimization for arrays of water turbine OWC devices}

\shorttitle{Modeling and optimization for arrays of water turbine OWC devices}

\shortauthors{M. Gambarini, G. Agate, G. Ciaramella, E. Miglio, S. Maran}

\affiliation[1]{organization={MOX, Dipartimento di Matematica, Politecnico di Milano},
addressline={Piazza Leonardo da Vinci 32},
postcode={20133},
city={Milano},
country={Italy}
}

\author[1]{Marco Gambarini}[orcid=0000-0001-9562-1456]
\ead{marco.gambarini@polimi.it}
\cormark[1]
\credit{Methodology, Software, Writing - Original Draft}

\author[2]{Giordano Agate}
\credit{Validation, Data Curation, Writing - Review \& Editing}

\author[1]{Gabriele Ciaramella}
\credit{Methodology, Writing - Review \& Editing}

\author[1]{Edie Miglio}
\credit{Methodology, Writing - Review \& Editing}

\author[2]{Stefano Maran}
\credit{Conceptualization, Writing - Review \& Editing, Supervision}

\affiliation[2]{organization={Ricerca sul Sistema Energetico - RSE S.p.A.},
addressline={Via R. Rubattino 54},
postcode={20134},
city={Milano},
country={Italy}
}

\cortext[cor1]{Corresponding author}


\begin{abstract}
Wave energy conversion is emerging as a promising technology for generating energy from renewable sources. Large-scale implementation of this technology requires the installation of parks of devices.
We study the problem of optimizing the park layout and control for wave energy converters of the oscillating water column type. As a test case, we consider a device with a semi-submerged chamber and a Wells turbine working in the liquid phase. First, a novel model based on a nonlinear ordinary differential equation is derived to describe the behavior of the water column and used to estimate the power matrix of an isolated device. Then, its linearization is derived in order to enable the fast simulation of large parks with a high number of devices. 
The choice of the hydrodynamic model allows obtaining the gradient of the power with respect to the positions through an adjoint approach, making it especially convenient for optimization. We consider in particular the case of interaction with the piles of a floating wind energy plant.
The results from the developed computational framework allow us to draw interesting conclusions that are useful when designing the layout of a park. In particular, we observe that interaction effects can be significant even in parks made up of devices of small size, which would exhibit negligible diffraction and radiation properties in isolated conditions, if the number of devices is large enough. Moreover, results show that wave reflection from the piles of an offshore platform can have positive effects on energy production.
\end{abstract}


\maketitle

\section{Introduction}
Despite having been studied and developed over many decades, wave energy conversion is having a very limited use. In recent years, however, interest has been surging, and installations are expected to increase. The EU Offshore Renewable Energy Strategy targets 100 MW of wave and tidal installations by 2025 and at least 1 GW by 2030. Some wave energy converter (WEC) device concepts have reached a mature stage of development, and institutions have set specific goals on the reduction of the levelized cost of energy from waves \cite{Tapoglou2022}.
Wave energy has been indicated as complementary to wind energy, in that combined installations can provide a more reliable power output, with a reduction in the annual hours of zero production with respect to using a single source \cite{Kalogeri2017}. Furthermore, many ports are implementing strategies to reduce their carbon footprint by switching to renewable power sources, and wave and tidal energy are being considered among the possible choices \cite{Bonamano2023, Parhamfar2023}.

Economically feasible implementation of wave energy requires building arrays of devices. This poses the problem of designing park layouts and devising control strategies that maximize energy production.
As reported in the recent reviews \cite{Yang2022, Golbaz2022}, most works on layout optimization for WEC parks have considered numbers of bodies in the order of the tens, up to about 100, and have used mainly meta-heuristic methods.
An exception is \cite{Goeteman2015}, in which arrays in the order of the hundreds of bodies are simulated, and parameter exploration is performed, by adopting a single-body approximation: mutual hydrodynamic interactions between objects due to diffraction are neglected; interactions due to radiation are instead taken into account. 
Regarding complementarity with wind energy, aero-hydro-servo studies of small parks of WEC combined with wind turbines have recently been performed \cite{Cao2023, Zhu2024}. Their main findings are that the presence of WECs has a small effect on the power output of the wind turbine, a beneficial effect in terms of alleviation of dynamic loads and pitch motion amplitude, and a possible unfavorable effect in terms of heave motion amplitude.

The aim of this work is to present a modeling and optimization strategy for parks of oscillating water column devices (OWC).
As a test case, we consider a water turbine oscillating water column converter called WaveSAX and developed by RSE \cite{Bonamano2023, Peviani2022}. The device was initially conceived for integration in coastal structures (e.g., harbours and ports), and then a second version was developed for its offshore use, possibly in connection with an offshore wind park. It consists of a vertical pipe in which water moves upward and downward, forced by the wave motion. Inside the pipe, a Wells turbine is positioned in the liquid phase and its shaft is connected to an electrical generator. The Wells turbine is of bi-directional type, i.e., the rotor rotates in the same direction during both the ascending and the descending phase of water motion. The main advantages of the device are its low cost and its modularity, as it can be installed individually or in batteries of several elements. The device has been registered at the European Patent Office\footnote{Patent Document n. 2 848 802 B1, European Patent Bulletin 2016/23.}.

We introduce a lumped-parameter mathematical model for the description of a single device, that is simple enough to allow fast computation, but still able to capture the most important features of the apparatus. It is constructed by considering the integral equations of mass and momentum for the dynamics of the water column, linear potential wave theory for the problem of hydrodynamic interaction with the surrounding environment \cite{Newman2018}, and a description of the Wells turbine through characteristic curves. The resulting nonlinear, time-domain model is used for the determination of the optimal turbine rotational speed for each sea state of interest. The corresponding linearized problem, that can be solved in the frequency domain, is used to build the multibody model adopted for the simulation and optimization of arrays. 
Layout optimization is then tackled, limiting our scope to cylindrical devices. For such geometry, a semi-analytical numerical model is available. It is presented in \cite{Yilmaz1998, Child2010} and based on the interaction theory introduced in \cite{Kagemoto1986}. It allows relatively fast computation of oscillation amplitudes and average power for very large arrays of bodies (up to 150 in our test cases), and it does not adopt the single-body approximation mentioned before, taking instead into account the complete linear hydrodynamic problem, including diffraction. Moreover, such formulation allows the computation of the gradient of the cost function through an adjoint approach. Standard gradient descent methods (see, e.g., \cite{Nocedal2006}) can then be used to seek a solution to the minimization problem.

We remark that the models developed in this work are not strictly linked to the specific device considered: the single-device model can be extended to any OWC device, and the park optimization method is suited for both OWC devices and point absorbers.

The main contributions of this work are the following:
\begin{itemize}
\item the definition of a mathematical model for a water-turbine OWC device; \item the study of the hydrodynamic interactions among the devices of a park, and of the park with the piles of an offshore structure; 
\item the implementation of an optimization method for park layouts and the critical discussion of the obtained results.
\end{itemize}

The article is structured as follows.
In Section~\ref{sec:isodev}, we consider a single device. A nonlinear time dependent model, derived in \ref{sec:intmodel}, describes the evolution of the water column elevation. It is influenced by the behavior of the Wells turbine, detailed in \ref{sec:wellsturb}, and by the external flow field, whose modeling is presented in \ref{sec:extmodel}. In Section~\ref{sec:contrpowmat}, the model is used to compute the power matrix of the device. A linearization of the model, derived in \ref{sec:linmod}, is used in \ref{sec:dimopt} to analyze the dependence of annual power on the dimensions of the device; the wave climates of two locations in the Mediterranean sea are considered. 

In Section \ref{sec:multibodywavesax}, we consider arrays of devices. A linear, frequency domain model for park analysis and optimization is introduced in \ref{sec:parkmodel}. Sections \ref{sec:1devpiles} and \ref{sec:2dev} present preliminar analyses of the interaction between a single device and the piles of an offshore platform, and of the interaction between two devices. Then, an optimization algorithm for the positions of the devices is introduced in \ref{sec:optmethod} and used in \ref{sec:optresults} to optimize parks with large numbers of converters. Finally, the results of optimization are verified using the nonlinear model in \ref{sec:linnonlin}.

\section{Isolated device}\label{sec:isodev}
\subsection{Mathematical model}
The operation of the device is described by two coupled mathematical models: an ordinary differential equation (ODE) for the water column and the linear potential model for the external wave problem. 
The equation for the evolution in time of the water column level $\zeta$ with incident waves of frequency $\omega$, to be derived in Sec.~\ref{sec:intmodel}, is
\begin{equation}
\left[\rho C(\zeta) S(\zeta) + A(\omega) \right] \ddot{\zeta} + \rho C(\zeta) \frac{dS}{d\zeta}\dot{\zeta}^2 + \frac{1}{2} \rho \dot{\zeta}^2 \left(1 - \frac{S^2(\zeta)}{S^2(z_1)} \right) + \Delta p(v_t, \omega_t) +  B(\omega)\dot{\zeta} +  \rho g \zeta = p_e(t).
\end{equation}
Here, $\rho$ is the water density, $C(\zeta)$ and $S(\zeta)$ are geometrical parameters, and $g$ is the gravity field. Frequency-dependent quantities $A(\omega)$, $B(\omega)$ and function $p_e(t)$ are hydrodynamic properties computed by solving the external problem, detailed in Sec.~\ref{sec:extmodel}. $\Delta p$ is the pressure differential due to the turbine, which is obtained from a nondimensional characteristic curve $C_a$ as 
\begin{equation}
\Delta p (v_t, \omega_t) = \frac{C_a K_a}{S_t} (v_t^2 + \omega_t r_t^2),
\end{equation}
where $v_t$ is the axial velocity, $\omega_t$ is the rotational speed, and $K_a$, $S_t$ and $r_t$ are geometric parameters. The behavior of the turbine will be detailed in Sec.~\ref{sec:wellsturb}. Once the coupled problem has been solved, mechanical power is computed as $P = \mathcal{T} \omega_t$, where $\mathcal{T}$ is the torque, obtained from the characteristic curve $C_t$ as
\begin{equation}
\mathcal{T} (v_t, \omega_t) = C_t K_a r_t (v_t^2 + \omega_t^2 r_t^2).
\end{equation}

\begin{figure}[h!]
\centering
\subfloat{
\begin{overpic}[width=0.3\textwidth]{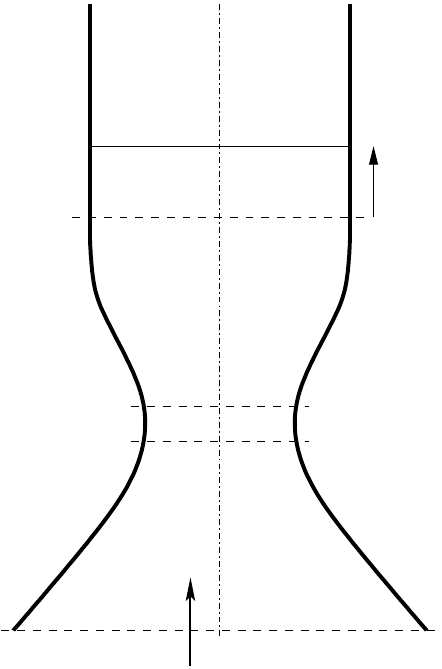}
\put(23,-2) {$Q$}
\put(-7, 5) {$S_1$}
\put(11,31) {$S_2$}
\put(11,38) {$S_3$}
\put(3,65) {$S_4$}
\put(58,78) {$\zeta$, $z$}
\end{overpic}}
\hspace{2cm}
\subfloat{
\vspace{1cm}
\input{fig-outer}
}
\caption{Domain of the internal (left; see Sec.~\ref{sec:intmodel}) and external (right; see Sec.~\ref{sec:extmodel}) problem}
\label{fig:scheme}
\end{figure}
The domains of the internal and external models are depicted in Fig.~\ref{fig:scheme}. The left panel shows the internal domain. 
Here, $S_1$ is the inflow section; sections $S_2$ and $S_3$ are located immediately below and immediately above the turbine section respectively, while $S_4$ is located at the equilibrium water level, corresponding to $z=0$. 
Finally, $Q$ is the volume flow rate through the inlet section.
The right panel shows a 2D slice of the domain $\Omega$ of the external wave problem. Its boundaries are the inflow section $\Gamma_f$, the external walls of the device $\Gamma_w$, the sea bottom $\Gamma_b$ and the mean free surface $\Gamma_s$.
The domain $\Omega$ is unbounded horizontally. 

To make park optimization less computationally demanding, the following linearized form of the equation for the evolution of the water column is used:
\begin{equation}
[\rho C(0) S(0) + A(\omega)] \ddot{\zeta} + [\Lambda S(0) + B(\omega)] \dot{\zeta} + \rho g \zeta = p_e(t).
\end{equation}
Its derivation is shown in Sec.~\ref{sec:linmod}.

\subsubsection{Internal model}\label{sec:intmodel}
For the water column, neglecting the effect of viscosity and considering incompressible flow, mass conservation requires
\begin{equation}
Q(t) = S(\zeta) \dot{\zeta},
\label{eq:Qdef}
\end{equation}
while the momentum balance equation is Euler's equation
\begin{equation}
\rho \dfrac{\partial \bm{u}}{\partial t}  +\rho (\bm{u} \cdot \nabla) \bm{u} = -\nabla p - \rho g \widehat{\bm{k}},
\label{eq:momentumbal}
\end{equation}
where $\bm{u}$ is the velocity and $\widehat{\bm{k}}$ is the upwards unit vector.
We introduce the quasi-1D approximation
$
\bm{u}(\bm{x}, t) \approx [0, 0, Q(t)/S(z)]$, $p(\bm{x}, t) \approx p(z, t),
$
which, replaced in \eqref{eq:momentumbal}, yields
\begin{equation}
\frac{\partial p}{\partial z} = -\frac{\rho \dot{Q}}{S(z)} + \frac{\rho Q^2(t)}{S^3(z)} \frac{\partial S(z)}{\partial z}  - \rho g.
\end{equation}
By integration, it is then possible to obtain the pressure difference between the section below the turbine $S_2$ and the inflow section:
\begin{equation}
p_2 - p_1 = -\rho \dot{Q} \int_{z_1}^{z_2} \dfrac{dz}{S(z)}  + \rho Q^2(t) \int_{z_1}^{z_2} \frac{\partial S(z)}{\partial z} \frac{{dz}}{S^3(z)} - \rho g (z_2 - z_1),
\label{eq:p2p1}
\end{equation}
and the same can be done for the pressure difference between the atmosphere and the section immediately above the turbine $S_3$:
\begin{equation}
p_{atm} - p_3 = -\rho \dot{Q} \int_{z_3}^{\zeta} \dfrac{dz}{S(z)} + \rho Q^2(t) \int_{z_3}^{\zeta} \frac{\partial S(z)}{\partial z} \frac{{dz}}{S^3(z)} - \rho g (\zeta - z_3).
\label{eq:patmp3}
\end{equation}
For the time being, we generically denote the pressure jump due to the turbine as $\Delta p$: then 
\begin{equation}
p_3 - p_2 = -\Delta p.
\label{eq:p3p2}
\end{equation} 
We also consider negligible the thickness of the turbine section: $z_2 \approx z_3$. 
The integrals appearing in the above expressions are specified through the definition
\begin{equation}
C(\zeta) = \int_{z_1}^\zeta \frac{dz}{S(z)}
\end{equation}
and the computation
\begin{equation}
\int_{z_1}^{\zeta} \frac{\partial S(z)}{\partial z} \frac{dz}{S^3(z)} = \int_{S(z_1)}^{S(\zeta)} \frac{dS}{S^3} = -\frac{1}{2} \left(  \frac{1}{S^2(\zeta)} - \frac{1}{S^2(z_1)} \right).
\end{equation}

The time derivative of the flowrate, computed from its expression \eqref{eq:Qdef}, involves the derivative of the section along the axial coordinate:
\begin{equation}
\dot{Q}(t) = S(\zeta)\ddot{\zeta} + \frac{dS}{d\zeta}\dot{\zeta}^2.
\label{eq:Qdot}
\end{equation}
The last contribution is zero if the section is uniform in the upper part of the device. From \eqref{eq:p2p1}, \eqref{eq:patmp3}, and \eqref{eq:p3p2}, we obtain the balance equation
\begin{equation}
\rho C(\zeta) S(\zeta) \ddot{\zeta}  + \rho C(\zeta) \frac{dS}{d\zeta}\dot{\zeta}^2 + \frac{1}{2} \rho \dot{\zeta}^2 \left(1 - \frac{S^2(\zeta)}{S^2(z_1)} \right) + \Delta p(v_t, \omega_t) + \rho g \zeta = p_1 - p_{atm} + \rho g z_1.
\label{eq:nonlindyn}
\end{equation}
We recognize that the forcing term is the pressure deviation from hydrostatic conditions at section 1:
\begin{equation}
p_1' = p_1 - (p_{atm} - \rho g z_1).
\label{eq:p1primedef}
\end{equation}

\subsubsection{Wells turbine model}\label{sec:wellsturb}
The Wells turbine consists of blades generated from the extrusion of a symmetric airfoil and mounted symmetrically on a hub. 
We shall call $\omega_t$ the angular velocity of the turbine, $r_t$ the distance between the axis of the turbine and the blade tips, $r_h$ the hub radius, $c$ the blade chord, $n$ the number of blades and $\mathcal{T}$ the torque exerted on the turbine by the flow. 
The area of the turbine's flow section, delimited by the hub and the duct walls, is $S_t = \pi (R_2^2 - r_h^2)$, where $R_2$ is the radius of the section at $z_2 = z_3$.
The hydraulic power lost by the flow due to the action of the turbine is $P_{hyd} = Q \Delta p$; the mechanical power available at the shaft is instead $P = \mathcal{T} \omega_t$. 

We now recall the definitions of a few common and useful nondimensional variables. The flow coefficient $\varphi$ is the ratio between the axial velocity $v_t = Q/S_t$ and the tip speed:
\begin{equation}
\varphi \equiv \frac{v_t}{\omega_t r_t}.
\label{eq:defflowcoeff}
\end{equation}
The Reynolds number relative to the chord is
\begin{equation}
\text{Re} \equiv \frac{\sqrt{v_t^2 + r_t^2\, \omega_t^2}\, c}{\nu} = \frac{r_t \omega_t c \sqrt{ \varphi^2 + 1}}{\nu},
\end{equation}
where $\nu$ is the kinematic viscosity.
We use the following definition\footnote{Another common definition in the literature is $ \Delta p^* \equiv \Delta p /(\rho \omega_t^2 r_t^2)$ (see \cite{Curran1997}).} of the nondimensional pressure jump \cite{Lekube2018}:
\begin{equation}
C_a \equiv \frac{\Delta p S_t}{K_a (v_t^2 + \omega_t^2 r_t^2)},
\label{eq:defdpstar}
\end{equation}
where $K_a = \rho c (r_t - r_h) n/2$. Likewise, the nondimensional torque is defined\footnote{We also report the alternative definition $\mathcal{T}^* \equiv \mathcal{T}/(\rho \omega_t^2 r_t^5)$.} as   
\begin{equation}
C_t \equiv \frac{\mathcal{T}}{K_a r_t (v_t^2 + \omega_t^2 r_t^2)}.
\label{eq:defTstar}
\end{equation}
In the following, we will consider a formulation based on $C_a$ and $C_t$, for which literature values are readily available \cite{Lekube2018}. For clarity, we report the characteristic equations
\begin{align}
\Delta p (v_t, \omega_t) &= \frac{C_a K_a}{S_t} (v_t^2 + \omega_t r_t^2), \label{eq:wellsdp} \\ 
\mathcal{T} (v_t, \omega_t) &= C_t K_a r_t (v_t^2 + \omega_t^2 r_t^2).
\label{eq:wellst}
\end{align}
Finally, the efficiency is 
\begin{equation}
\eta = \frac{P}{P_{hyd}}.
\end{equation}
We observe that factor $c(r_t - r_h)n$ in $K_a$ is the total area of the blades, if they are rectangular. Hence, by using the definition of solidity $\sigma$ as the ratio of blade area to the area of the turbine flow section $S_t$, we have $K_a = \rho \sigma S_t/2$ and thus the nondimensional pressure coefficient can be written as 
\begin{equation}
C_a = \frac{\Delta p}{\sigma q_\infty}, \quad q_\infty \equiv \frac{1}{2} \rho (v_t^2 + \omega_t^2 r_t^2),
\end{equation}
where $q_\infty$ is the dynamic pressure at the blade tips. Analogously,
\begin{equation}
C_t = \frac{\mathcal{T}}{\sigma S_t r_t q_\infty}.
\end{equation}
The pair of functions $C_a(\varphi)$, $C_t(\varphi)$, or alternatively $\Delta p^*(\varphi)$, $\mathcal{T}^*(\varphi)$, represent the characteristic curves of the turbine. 

In order to explain the working principle of the turbine and the shape of the characteristic curves, let us consider a single blade. The angle of attack of a blade section at a distance $r$ from the axis of rotation is $\alpha(r) = \arctan(\varphi r_t/r)$: it reaches its maximum value at the hub and its minimum at the blade tips. When the flow coefficient is zero, the blade is at zero incidence. 
In such conditions, since blades have symmetric airfoils, the lift force is zero, while the drag force is non-zero and produces a negative torque. 
If a constant angular velocity control strategy is chosen, then such negative torque must be compensated by a motor torque produced by feeding energy to the system. 
If, instead, the angular velocity is allowed to vary, then mechanical inertia may absorb this effect. 
As the flow coefficient increases, the lift force on blades increases. Both the nondimensional pressure jump and the nondimensional torque increase; correspondingly, the efficiency increases, until it reaches a plateau at values typically around 40\% to 60\%. 
As the angle of stall is reached, the lift force starts decreasing; this happens more or less abruptly depending on the shape of the airfoil (in particular, on its thickness) and on the turbulence intensity. 
Thick airfoils typically exhibit a soft stall behavior, produced by flow separating at the trailing edge; conversely, thin airfoils have a more abrupt stalling behavior, corresponding to leading edge flow separation \cite{Anderson1991}. Greater turbulence intensities can delay stall and mitigate its effects \cite{Swalwell2001}. 
Stall causes a reduction in nondimensional torque, while the nondimensional pressure difference keeps increasing. Hence, efficiency decreases until it becomes negligible for values of the flow coefficient around 0.3 to 0.4. 
We observe that $\varphi=0.4$ corresponds to an angle of attack of about $21.8^{\circ}$ at the blade tip, a very large value for symmetric airfoils.

An additional effect that impacts on the performance of the turbine is due to the value of the distance between the blade tips and the duct walls, called tip clearance or tip gap (TG). 
In \cite{Torresi2008}, the results of some simulations are reported, showing that the largest efficiencies can be reached for small values of tip clearance, of the order of 1\%. Intermediate values of leads to worse performance for small values of the flow coefficient, while tip clearances larger than about 10\% cause low efficiency in all operating conditions. 

Cavitation is another important, and potentially damaging, physical phenomenon for turbines immersed in water. 
It occurs when the pressure falls below the vapor pressure of the liquid \cite{Franc2005}. For water at $20^{\circ} \text{C}$, vapor pressure is $p_v = 2.34 \text{ kPa}$. 
The minimum pressure on the turbine can be approximately computed by considering the pressure coefficient
\begin{equation}
C_p = \frac{p - p_\infty}{1/2 \rho U_\infty^2}.
\end{equation}
This quantity, for large values of the Reynolds number and small values of the Mach number, does not depend on velocity $U_\infty$. 
Thus, for a turbine located at depth $z_t$ below mean sea level, the minimum absolute pressure can be computed as
\begin{equation}
p_\text{min} = p_\text{atm} + \rho g (\zeta - z_t) + \frac{1}{2} \rho (v_t^2 + \omega_t^2 r_t^2) C_{p,\text{min}},
\label{eq:pcavit}
\end{equation}
where $C_{p,\text{min}}$ can be obtained from literature data or from simulations either of the full blade (3D) or of its section (2D). 
The condition of incipient cavitation is $p_\text{min} = p_v$. 
It is reasonable to expect the minimum pressure to occur on the suction surface, close to the leading edge.\footnote{In this case, talking about upper surface and lower surface could be misleading, since their roles as suction surface and pressure surface are inverted at every half-cycle of oscillation of the water column.}
In our case, the turbine is not simply located at some depth below sea level: it is immersed in a water column of time-varying depth. Referring to Fig.~\ref{fig:scheme}, the freestream pressure is $p_2$ when the flow is ascending, and $p_3$ when the flow is descending. Using the results from Sec.~\ref{sec:intmodel}, we obtain
\begin{equation}
\begin{split}
p_\text{min} &= p_\text{atm} + \rho g (\zeta - z_t) + \frac{1}{2} \rho (v_t^2 + \omega_t^2 r_t^2) C_{p,\text{min}} 
\\& + \rho \left[\frac{dS}{d\zeta} \dot{\zeta}^2 + S(\zeta) \ddot{\zeta}\right][C(\zeta) - C(z_t)] + \frac{1}{2} \rho \left[ 1 - \frac{S^2(\zeta)}{S^2(z_t)}\right] + \begin{cases}
\Delta p(v_t, \omega_t) & \text{if $\dot{\zeta}>0$},
 \\ 0 & \text{otherwise}
\end{cases}
\end{split}
\end{equation}

To ensure that the turbine works correctly, it must also always remain immersed. 
The fact that the constraint of no cavitation is fulfilled does not imply that the constraint of immersed turbine is satisfied, nor vice versa.

In principle, to correctly model the behavior of the turbine, one would have to write a dynamic equation for its rotational motion, including, as forcing terms, the hydrodynamic torque and the electromagnetic torque. 
At this stage of the design of the device, full-scale data about turbine inertia and about the electrical generator are not known. 
For this reason, in our simulations we adopted the simplifying assumption of constant rotational speed, which in practice would need to be realized through an appropriate controller. 

\begin{figure}[h!]
\subfloat{
\includegraphics[width=0.45\textwidth]{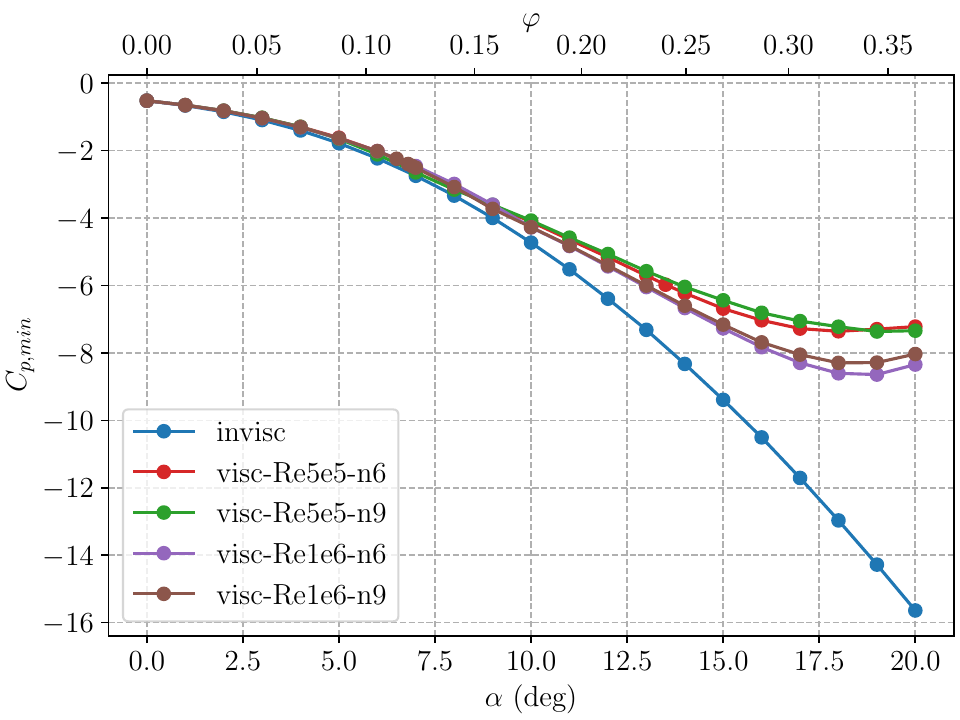}
}
\hfill
\subfloat{
\includegraphics[width=0.45\textwidth]{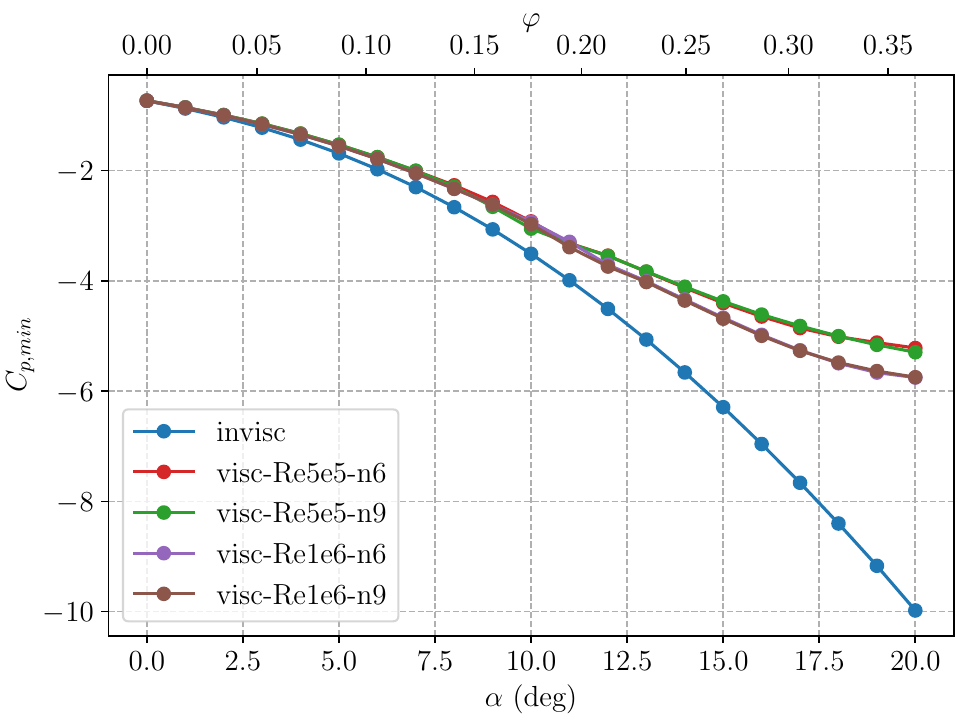}
}
\caption{Coefficient $C_{p,\text{min}}$ as a function of angle of attack and flux coefficient; NACA0015 (left), NACA0021 (right)}
\label{fig:cpminxfoil}
\end{figure}

Torque and pressure difference are computed from the data of \cite{Lekube2018}. 
Coefficient $C_{p,\text{min}}$, needed for cavitation computations, has been obtained from simulations of a single blade section performed using Xfoil \cite{Drela1989}. We repeated the computation for the NACA0021 airfoil, considered in \cite{Peviani2020}, and the NACA0015 airfoil, considered in \cite{Peviani2021}. 
Results are reported in Fig.~\ref{fig:cpminxfoil}. We first notice that the results of inviscid simulations differ greatly from the ones of viscous, turbulent simulations; we will refer to the latter. The minimum pressure coefficient is larger (in absolute value) for the NACA0015 airfoil, due to its smaller thickness, which induces a smaller radius of curvature of the streamlines and thus more intense depressions, being the pressure gradient linked to centripetal acceleration. 
The results indicate that the NACA0021 airfoil may be less prone to cavitation. We also notice that airfoils become more sensitive to cavitation as the Reynolds number increases. 
It is reasonable to expect that the data thus obtained yields cautionary estimates of the condition of incipient cavitation. Indeed, the turbine blades have a very small aspect ratio, which may lead to important tip effects. 
Such effects should reduce the intensity of local suction peaks thanks to the current generated between the pressure side and the suction side.

\subsubsection{External model}\label{sec:extmodel}
The pressure at the inflow section $p_1'$, appearing in the internal model \eqref{eq:nonlindyn}, can be computed by using the linear potential model \cite{Chakrabarti1987}. This is a very common modeling choice for point absorber devices and in general for floating bodies \cite{Newman2018}. 
The potential is written as the sum of the potential of ambient incident waves $\phi_i$, a diffraction potential $\phi_d$, produced by incident waves with the condition of zero volume flux through the inflow section $\Gamma_f$ (corresponding to $S_1)$, and radiation potential $\phi_r$, obtained, thanks to linearity, by multiplying the average velocity on section $S_1$ by the potential $\varphi_r$ corresponding to unit amplitude velocity fluctuations.
In the frequency domain, we have
\begin{equation}
\widehat{\phi} = \widehat{\phi}_i + \widehat{\phi}_d + \frac{\widehat{Q}}{S_1} \varphi_r.
\label{eq:phidecomp}
\end{equation} 
For any given frequency, \eqref{eq:phidecomp} allows writing the mean pressure on the inflow section as an integral:
\begin{equation}
\begin{split}
\widehat{p}'_1 =& i\omega \rho \dfrac{1}{S_1} \int_{\Gamma_f} \widehat{\phi} \, d\Gamma \\
=&\underbrace{i\omega\rho \dfrac{1}{S_1} \int_{\Gamma_f} (\widehat{\phi}_i + \widehat{\phi}_d) \, d\Gamma}_{\text{excitation pressure $\widehat{p}_e(\omega)$}} + \underbrace{(-i\omega \widehat{\zeta}) i\omega\rho \dfrac{S_4}{S_1} \dfrac{1}{S_1} \int_{\Gamma_f} \varphi_r \, d\Gamma.}_{\text{radiation pressure $\widehat{p}_r = -R(\omega)(-i\omega\widehat{\zeta})$}}
\end{split}
\label{eq:p1def}
\end{equation}
This quantity depends on the state of the water column through $\widehat{\zeta}(\omega)$, the Fourier transform of $\zeta(t)$. From the potential $\varphi_r$, the coefficients of added mass $A(\omega)$ and radiation damping $B(\omega)$ can be computed:
\begin{equation}
R(\omega) = -i\omega A(\omega) + B(\omega) \to \begin{dcases}
A(\omega) = \text{Im}\left[ i \rho \dfrac{S_4}{S_1} \dfrac{1}{S_1} \int_{\Gamma_f} \varphi_r \, d\Gamma \right] \\
B(\omega) = -\omega  \text{Re}\left[ i \rho \dfrac{S_4}{S_1} \dfrac{1}{S_1} \int_{\Gamma_f} \varphi_r \, d\Gamma \right].
\end{dcases}
\label{eq:radimpedance}
\end{equation}
Radiation damping describes the pressure contribution in phase with the flowrate. It corresponds to a net exchange of energy with the environment due to the generation of radiated waves by the device. 
Added mass describes a pressure contribution with a phase difference of $\pi/2$ with respect to the flowrate; thus, the time-averaged exchange of energy with the environment due to this term is zero. The diffraction potential is obtained from the solution of problem
\begin{equation}
\begin{cases}
\Delta \widehat{\phi}_d = 0 & \text{in $\Omega$}\\
\dfrac{\partial \widehat{\phi}_d}{\partial z} = 0  & \text{on $\Gamma_b$} \\[0.5em]
\dfrac{\partial \widehat{\phi}_d}{\partial z} - \dfrac{\omega^2}{g} \widehat{\phi}_d = 0 & \text{on $\Gamma_s$} \\[0.5em]
\dfrac{\partial \widehat{\phi}_d}{\partial n} = - \dfrac{\partial \widehat{\phi}_i}{\partial n} & \text{on $\Gamma_w \cup \Gamma_f$}, 
\end{cases}
\end{equation}
where
\begin{equation}
\widehat{\phi}_i = -i \dfrac{H}{2} \dfrac{g}{\omega} \dfrac{\cosh[k(z+h)]}{\cosh(kh)} \exp[ik(x\cos\theta + y\sin\theta)].
\label{eq:incwavepot}
\end{equation}
Here, $H$ is the wave height, $h$ is the water depth, $k$ is the wavenumber, $\theta$ is the wave direction and $(x,y)$ are the coordinates in the horizontal plane. We refer to Fig.~\ref{fig:scheme} for the nomenclature of the boundaries. The radiation potential satisfies
\begin{equation}
\begin{cases}
\Delta \varphi_r = 0 & \text{in $\Omega$}\\
\dfrac{\partial \varphi_r}{\partial z} = 0  & \text{on $\Gamma_b$} \\[0.5em]
\dfrac{\partial \varphi_r}{\partial z} - \dfrac{\omega^2}{g} \varphi_r = 0 & \text{on $\Gamma_s$} \\[0.5em]
\dfrac{\partial \varphi_r}{\partial n} = 0 & \text{on $\Gamma_w $} \\[0.5em]
\dfrac{\partial \varphi_r}{\partial n} = 1 & \text{on $\Gamma_f$}.
\end{cases}
\end{equation}
It is important to notice that in the radiation problem we have imposed uniform velocity directed upwards on the inflow section. 
In practice, on the inflow section the radial component of velocity may not be negligible, and the vertical component may not be uniform.
The linear potential model would allow to impose a non-uniform vertical velocity on the section. It would instead not be possible to impose a radial component of velocity, because only the component of velocity normal to the boundary (and thus, axial) appears in the boundary condition.

Since the domain is unbounded, a radiation condition needs to be enforced at infinity for both problems; at the numerical level, this is satisfied by the Green functions used. 
Hydrodynamic simulations have been carried out by employing the open-source code Capytaine \cite{Ancellin2019}, which solves the Laplace equation using the boundary element method.
\begin{figure}[h!]
\subfloat{
\includegraphics[width=0.45\textwidth]{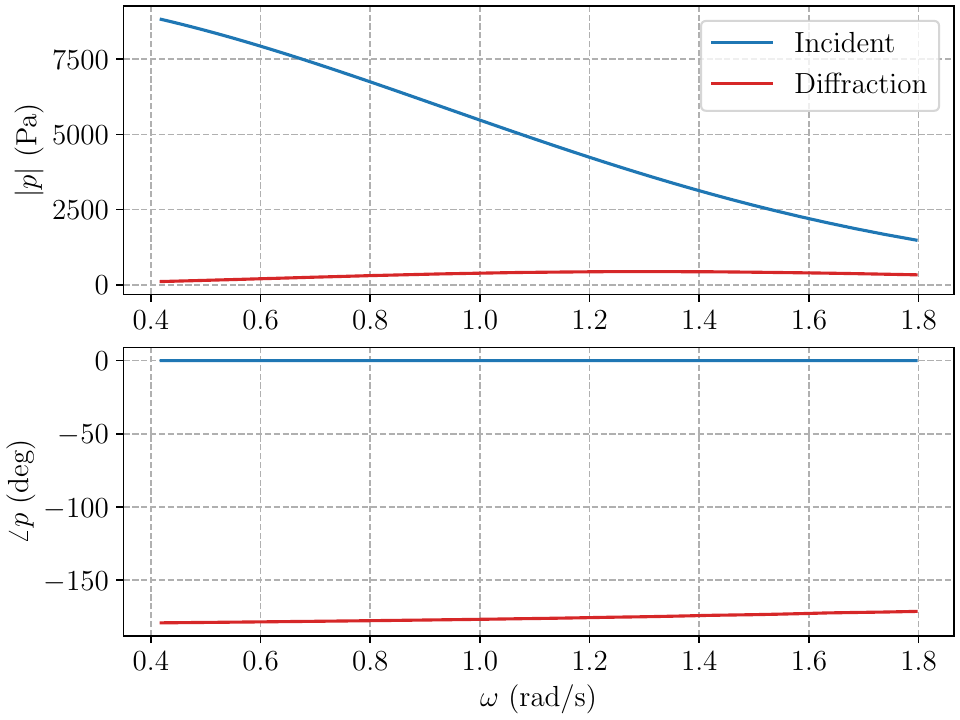}
}
\hfill
\subfloat{
\includegraphics[width=0.45\textwidth]{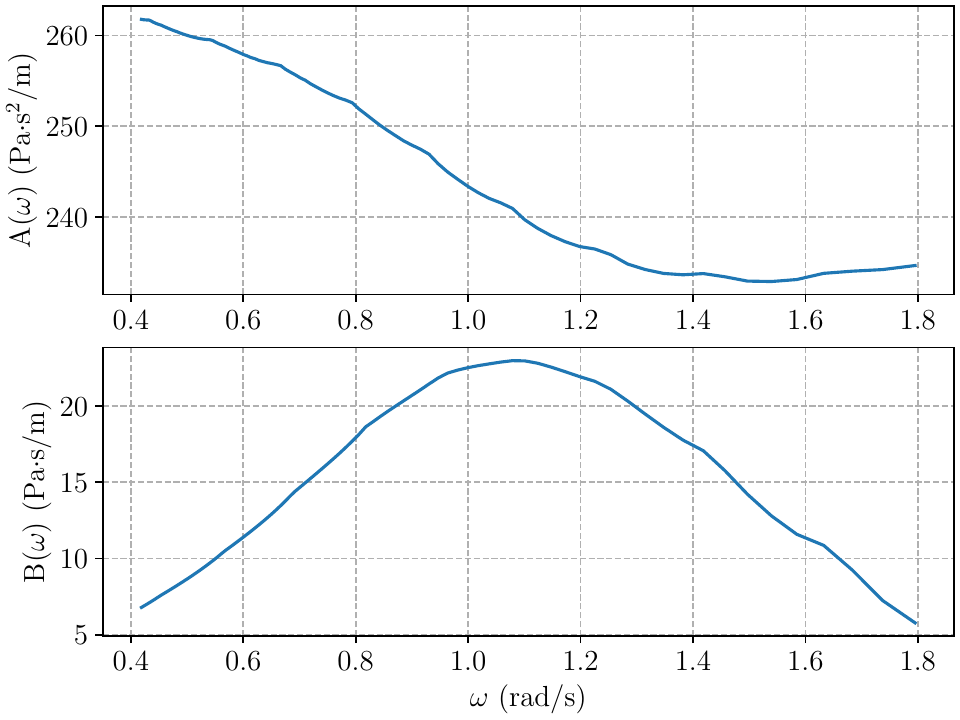}
}
\caption{Diffraction (left) and radiation data (right), from the solution of the external problem for an isolated device}
\label{fig:simdiffrad}
\end{figure}
The hydrodynamic properties of the considered device, obtained solving the external problem, are reported in Fig.~\ref{fig:simdiffrad}. We notice that, especially for low frequencies (long waves), the magnitude of the diffraction pressure is very small compared to the magnitude of the incident wave pressure. 
This is due to the small size of the device compared to the wavelength. The dominating contribution is thus due to the potential of incident waves \eqref{eq:incwavepot}. 
For $-h<z<0$, this is an increasing function of $z$. Hence, the mean pressure on the inflow section, which is the forcing term of the internal problem, decreases as the draft of the device increases.

Since the internal model is nonlinear, a forcing term due to a purely monochromatic wave would produce a response $\zeta(t)$ that is in general not monochromatic. To correctly represent this effect, one should write a convolution integral (in time) for the radiation term, corresponding to the product $R(\omega)\widehat{\zeta}(\omega)$ in the frequency domain. 
An accurate computation of such product would require solving the radiation problem for a potentially large number of frequencies, which would significantly increase the computational cost of simulation and optimization. 
We instead only use the values of added mass $A(\omega)$ and radiation damping $B(\omega)$ corresponding to the frequency $\omega$ of the incident monochromatic wave; the same is done for $\widehat{p}_e(\omega)$. Such approximation is reasonable due to the small magnitude of the radiation terms compared to the other terms of the equation. 
From \eqref{eq:p1def} and \eqref{eq:radimpedance}, and transforming back to the time domain while treating $A(\omega)$, $B(\omega)$ as constants as discussed above, we obtain $p_1'(t) = p_e(t) - A(\omega) \ddot{\zeta}(t) - B(\omega) \dot{\zeta}(t),$
which, substituted in the dynamic equation \eqref{eq:nonlindyn}, yields
\begin{equation}
\left[\rho C(\zeta) S(\zeta) + A(\omega) \right] \ddot{\zeta} + \rho C(\zeta) \frac{dS}{d\zeta}\dot{\zeta}^2 + \frac{1}{2} \rho \dot{\zeta}^2 \left(1 - \frac{S^2(\zeta)}{S^2(z_1)} \right) + \Delta p(v_t, \omega_t) +  B(\omega)\dot{\zeta} +  \rho g \zeta = p_e(t),
\label{eq:nonlindynrad}
\end{equation}
and where the forcing term is reconstructed from its complex coefficient as $
p_e(t) = \Re\left[ \widehat{p}_e(\omega) \exp(-i\omega t) \right].$
Equation \eqref{eq:nonlindynrad} is a nonlinear ordinary differential equation in the time domain, where the external forcing is a monochromatic wave of frequency $\omega$.
 
\subsubsection{Linearized model}\label{sec:linmod}
In this section, a linearization of the model introduced in Sec.~\ref{sec:intmodel} is derived and discussed. The advantages of adopting a linearized model are the possibility of obtaining analytical results that can be useful for an initial understanding of the behavior of the device, the lower computational cost due to the possibility of transforming the problem to the frequency domain and thus to an algebraic problem, and the possibility of coupling with the multi-body, semi-analytical hydrodynamic model presented in Sec.~\ref{sec:parkmodel}, which is especially suitable for optimization.

Linearization of the equation of motion \eqref{eq:nonlindyn} around the equilibrium condition $\zeta(t)=0$ yields
\begin{equation}
[\rho C(0) S(0) + A(\omega)] \ddot{\zeta} + [\Lambda S(0) + B(\omega)] \dot{\zeta} + \rho g \zeta = p_e(t),
\label{eq:lindyn}
\end{equation}
where $
\Lambda = \left. \frac{\partial \Delta p}{\partial Q} \right|_{Q=0}.$
To obtain this quantity from \eqref{eq:wellsdp}, we compute
\begin{equation}
\left. \dfrac{\partial \Delta p}{\partial v_t }\right|_{v_t=0} = \left[ \dfrac{\partial C_a}{\partial \varphi} \dfrac{\partial \varphi}{\partial v_t} \dfrac{K_a}{S_t}(v_t^2 + \omega_t^2 r_t^2) + \dfrac{2 C_a K_a v_t}{S_t} \right]_{v_t=0}  = \left. \dfrac{\partial C_a}{\partial \varphi} \right|_{\varphi=0} \dfrac{K_a \omega_t r_t}{S_t},
\end{equation}
from which, recalling that $Q = v_t S_t$,
\begin{equation}
\Lambda = \left.\dfrac{\partial \Delta p}{\partial Q} \right|_{Q=0} = \dfrac{1}{S_t}\left.\dfrac{\partial \Delta p}{\partial v_t} \right|_{v_t=0} = \left. \dfrac{\partial C_a}{\partial \varphi} \right|_{\varphi=0} \dfrac{K_a \omega_t r_t}{S_t^2}.
\label{eq:lambdadef}
\end{equation}
This results shows that the turbine has, up to first order, the role of a damping.

We now discuss the limits of validity  of the linearized model. In general, the linearization will lead to a good approximation if the oscillation amplitude is small. 
In addition, the terms multiplying $\dot{\zeta}^2$ in \eqref{eq:nonlindyn} will be small regardless of the magnitude of the oscillation amplitude if the variation of the area of the duct sections along $z$ is small close to $z=0$; under the same conditions, we have $C(\zeta)S(\zeta)\approx C(0)S(0)$. 
The linearization of the turbine's characteristic equation is acceptable if the flow coefficient is small, if, again, variations of section are negligible around $z=0$, and if $C_a(\varphi)$ is well approximated by a straight line. This last requirement is in general fulfilled, as shown in the used characteristic curves \cite{Lekube2018}. Expanding \eqref{eq:wellsdp} in $\dot{\zeta}$, one obtains
\begin{equation}
\Delta p = \frac{C_a'(0) K_a \omega_t r_t}{S_t^2} S(\zeta) \dot{\zeta} + \frac{C_a'(0) K_a}{\omega_t r_t S_t^5} S^3(\zeta) \dot{\zeta}^3 + o(\dot{\zeta}^3) 
= \Lambda S(\zeta) \dot{\zeta} (1 +  \varphi^2) + o(\dot{\zeta}^3),
\end{equation}
which shows that the ratio between the term of order 3 and the term of order 1 is $\varphi^2$. Typical values of $\varphi$ are less than or close to 0.3: the ratio is thus, at most, of order $1/10$ for conditions of interest.

Typical values of $\Lambda$ for the considered converter are of the order of $6 \cdot 10^3 \text{ Pa}\cdot\text{s/m}^3$. Since $S(0)$ is of the order of $1 \text{ m}^2$, and examining the values of radiation damping $B(\omega)$ from Fig.~\ref{fig:simdiffrad}, we observe that the latter is negligible compared to the contribution from the turbine. Likewise, in the mass term, term $\rho C(0) S(0)$ is about an order of magnitude larger than $A(\omega)$. By considering this approximation in the computation of the natural frequency of the linearized system, one obtains
\begin{equation}
\omega_0 = \sqrt{\frac{\rho C(0) S(0) + A}{\rho g}} \approx \sqrt{\frac{C(0)S(0)}{g}} \approx \sqrt{\frac{L}{g}},
\end{equation}
where the last step is an exact equality in the case of a tube of length $L$ and uniform section. The length of the duct is thus the dimension that influences the natural frequency the most.

In the frequency domain, the linearized dynamic equation reads
\begin{equation}
\left\{ -\omega^2 [\rho C(0) S(0) + A(\omega)] - i\omega[\Lambda S(0) + B(\omega)] + \rho g \right\} \widehat{\zeta} = \widehat{p}_e.
\label{eq:lindynfreqdom}
\end{equation}
The mean hydraulic power is 
\begin{equation}
P_\text{hyd} = \frac{1}{2}\Lambda \omega^2 S^2(0) |\widehat{\zeta}|^2.
\label{eq:hydpower}
\end{equation}
For the computation of the mean mechanical power, adopting a linearization of the nondimensional torque curve $C_t(\varphi)$ would be too inaccurate, since it would not be possible to describe the presence of the maximum point corresponding to stall. We instead consider an approximation given by an even polynomial
\begin{equation}
C_t(\varphi) = c_0 + c_2(\varphi) \varphi^2 + c_4(\varphi) \varphi^4 + \dots + c_{2d}(\varphi) \varphi^{2d},
\label{eq:torquepoly}
\end{equation}
whose degree $2d$ can be chosen in order to obtain a prescribed accuracy. This choice leads to the following expression for the instantaneous power:
\begin{equation}
P(t) = \mathcal{T}(t) \omega_t =  \left( c_0 + c_2 \frac{v_t^2(t)}{\omega_t^2 r_t^2} + c_4 \frac{v_t^2(t)}{\omega_t^2 r_t^2} + \dots   \right) K_a r_t \left( v_t^2(t) + \omega_t^2 r_t^2 \right) \omega_t,
\end{equation}
where the time dependence has been made explicit. Computing the time-average of this quantity involves integrals of the form
\begin{equation}
I_n = \frac{1}{2\pi} \int_{-\pi}^\pi \left[ \sin(x) \right]^{2n} {dx}.
\end{equation} 
Integration by parts yields the recursive expression 
\begin{equation}
I_n = \begin{cases}
0, & \text{if $n$ is odd,} \\
\frac{n-1}{n} I_{n-2}, & \text{if $n$ is even,}
\end{cases}
\end{equation}
from which we can obtain the explicit expression
\begin{equation}
I_n = \frac{(n-1)(n-3)\dots}{n(n-2)\dots} = \frac{(n-1)!!}{n!!}.
\end{equation}
The mean mechanical power can then be expressed as
\begin{equation}
P = \sum_{n=0}^d p_n \left(\omega |\widehat{\zeta}| \right)^{2n}, 
\end{equation}
with coefficients $p_n$ given by
\begin{align}
p_0 &= c_0 K_a \omega_t^3 r_t^3, \\
p_n &= I_n(c_{2(n-1)} + c_{2n}) K_a (r_t \omega_t)^{3 - 2n} \left(\frac{S(0)}{S_t} \right)^{2n}, \quad n=1,\dots,d-1, \\
p_d &= I_d c_{2(d-1)} K_a (r_t \omega_t)^{3-2d} \left(\frac{S(0)}{S(t)} \right)^{2d}.
\end{align}
By adopting this approach, passing to the time domain is not required for computing the mean power. The condition of low flux coefficient and thus low angle of attack, in which only the drag component of the hydrodynamic force is acting on the blades and power is negative, is represented by the condition $p_0<0$. Higher order terms describe the regimes of positive power and stall.

\subsection{Simulation and optimization}
\begin{figure}[h!]
\centering
\subfloat{
\includegraphics[width=0.48\linewidth]{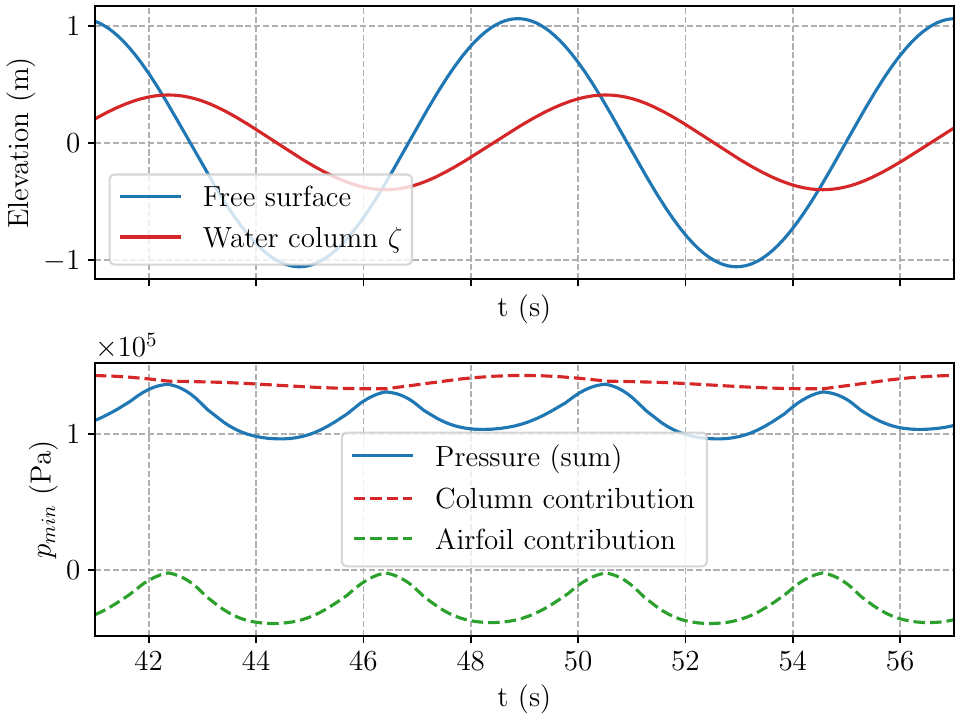}
}
\hfill
\subfloat{
\includegraphics[width=0.48\linewidth]{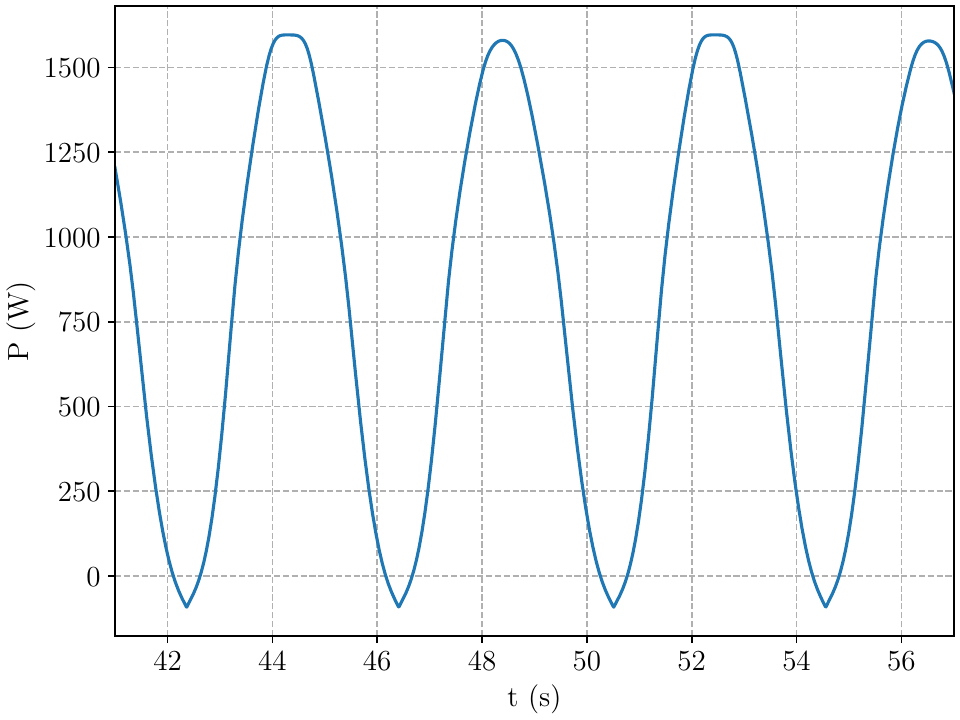}
}
\caption{
Water level, minimum pressure on blades (left) and extracted power (right) with 7-blade turbine, $H_s = 3$ m, $T = 8.15$ s.}
\label{fig:testsim}
\end{figure}

Figure \ref{fig:testsim} shows the time series of the main quantities of interest obtained from the solution of the dynamic equation \eqref{eq:nonlindynrad} for a reference wave condition. With reference to Fig.~\ref{fig:scheme}, the main dimensions of the device are $R_1 = 1.4\text{ m}$, $R_2 = 0.5\text{ m}$, $R_4 = 0.75 \text{ m}$, $z_2 = -3.650 \text{ m}$, $z_4 = -5.650 \text{ m}$.
A phase difference between the levels of the water column and of the free surface outside the device can be observed. The elevation head thus generated makes energy production possible. 
During each wave period, two maxima and two minima of power and pressure on blades occur. Every power maximum corresponds to a pressure minimum. 
The latter are observed when the water column is at the reference level ($\zeta = 0$), corresponding to a maximum value of flowrate. 
The pressure minima that occur when the internal free surface is ascending are less severe than the ones that occur when the internal free surface is descending. 
Intervals of negative power are observed around each maximum and each minimum point of the water column level. Power must then be provided to the system in order to keep the turbine in motion at constant speed.
Moreover, such maxima and minima occur approximately when the curves of internal and external water levels intersect. The intuitive explanation of this fact is that the tube can fill up only as long as the external level is higher than the internal level, and it can empty only as long as the internal level is higher than the external level. 
Positive power peaks appear clipped because the condition of maximum efficiency (optimal flow coefficient) is reached and slightly exceeded.


\subsubsection{Control optimization and power matrix estimation}\label{sec:contrpowmat}
\begin{figure}[h!]
\centering
\includegraphics[width=0.5\textwidth]{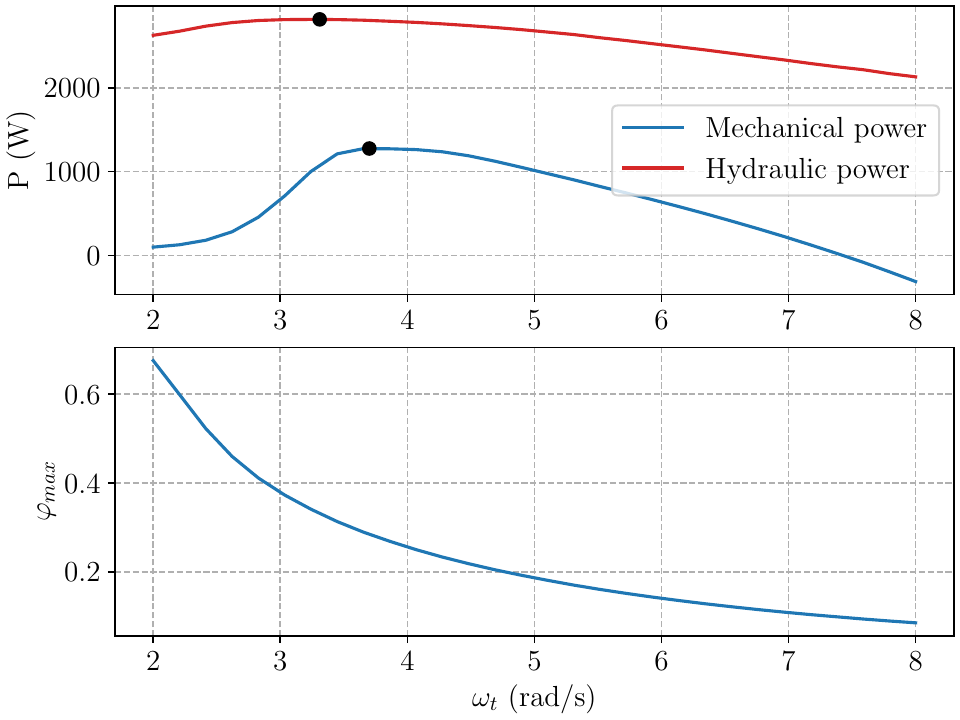}
\caption{
Mechanical power, hydraulic power and maximum flow coefficient as functions of the turbine's rotational speed}
\label{fig:mechvshydpower}
\end{figure}

From the linearized model \eqref{eq:lindynfreqdom} it is possible to compute analytically the value of $\Lambda$ that maximises the mean hydraulic power \eqref{eq:hydpower} when a monochromatic wave of angular frequency $\omega$ forces the system:
\begin{equation}
\Lambda_\text{max$P_{hyd}$}(\omega) = \frac{1}{S(0)} \sqrt{B^2(\omega) + \frac{\left[ -\omega^2 (\rho C(0) S(0) + A(\omega)) + \rho g \right]^2}{\omega^2}}.
\label{eq:lambdamaxphyd}
\end{equation}
The corresponding value of $\omega_t$ can be found by inverting \eqref{eq:lambdadef}. In practice, however, one does not seek to maximize the mean hydraulic power, but rather the mean mechanical power. For the latter it is not possible in general to determine the optimal value of $\omega_t$ analytically.
Nevertheless, the value of $\omega_t$ maximising hydraulic power can be used as the initial guess of an optimization process for mechanical power. The optimization problem was solved using function \texttt{minimize\_scalar} of \texttt{scipy}. A comparison of the two values is shown in Fig.~\ref{fig:mechvshydpower}. 
We observe that the maximum hydraulic power is reached for a turbine rotational speed that is lower than the one corresponding to maximum mechanical power. Indeed, the maximum flow coefficient reached when the rotational speed is set to the value of maximum hydraulic power corresponds to stall conditions and thus suboptimal efficiency. 
It can be further noted that while the rotational speed corresponding to maximum hydraulic power is independent on the height of the incident waves (which can be observed from \eqref{eq:lambdamaxphyd} and \eqref{eq:lambdadef}), the rotational speed that maximises the mechanical power must be chosen based on the flow coefficient, which depends on the oscillation amplitude and thus on the incident wave height.

The computation can be repeated for several sea states, obtaining the maximum power absorbable in each of them under the assumption of constant rotational speed. In this way, we obtain the power matrix of the device. 
In the literature, power matrices are generally defined in terms of peak period $T_p$, or energy period $T_e$, and significant wave height $H_s$ \cite{Babarit2012}. 
In the present work, each sea state has been approximated as a single monochromatic wave following the criteria described in \ref{sec:powmatparams}. 
The results are reported in \ref{sec:pmat}. We have reported, in addition to optimal power and turbine velocity, the maximum and minimum level of the water column and the minimum pressure on the turbine. The capture width ratio, a performance indicator, is also computed. It is defined as $
\text{CWR} = P/(2 R_1 J)$,
where $P$ is the mean mechanical power, $J$ is the energy flux of the ambient incident waves, and $R_1$ is the radius of the duct at the flow section, that is the section of maximum width. It can be observed that the oscillations of the internal free surface are small compared to the vertical scale of the device; in particular, their amplitude is lower for the 7-blade turbine than for the 5-blade turbine, due to the greater hydraulic resistance of the latter. The most energetic conditions lead to cavitation for the 5-blade turbine, while this effect is not present in the results for the 7-blade turbine.

As a comparison, a cylindrical device with constant section has been considered. Its radius is 75 cm, equal to the radius of the upper part of the original device, and it is equipped with a 7-blade turbine. Computations on this geometry yield a greater power than the original device. The capture width ratio, in particular, is significantly larger. For this reason, in the following we will consider only the device with constant section. In all simulated conditions, the minimum pressure is above the vapor pressure.

Finally, the computations for the constant section device have been repeated using the linearized model introduced in Sec.~\ref{sec:linmod}. The linear model slightly overestimates the power, especially for the most energetic wave conditions, while the results in terms of optimal rotational speed are comparable. 

\subsubsection{Influence of device dimensions on annual power}\label{sec:dimopt}
\begin{figure}[h!]
\centering
\subfloat{
\includegraphics[width=0.32\textwidth]{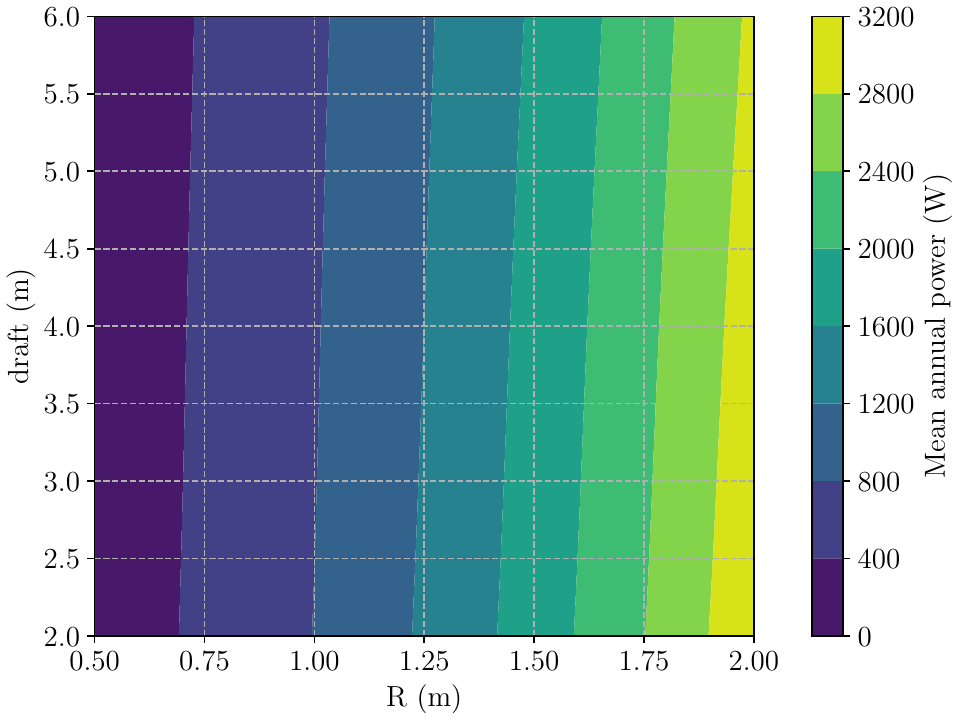}
}
\subfloat{
\includegraphics[width=0.32\textwidth]{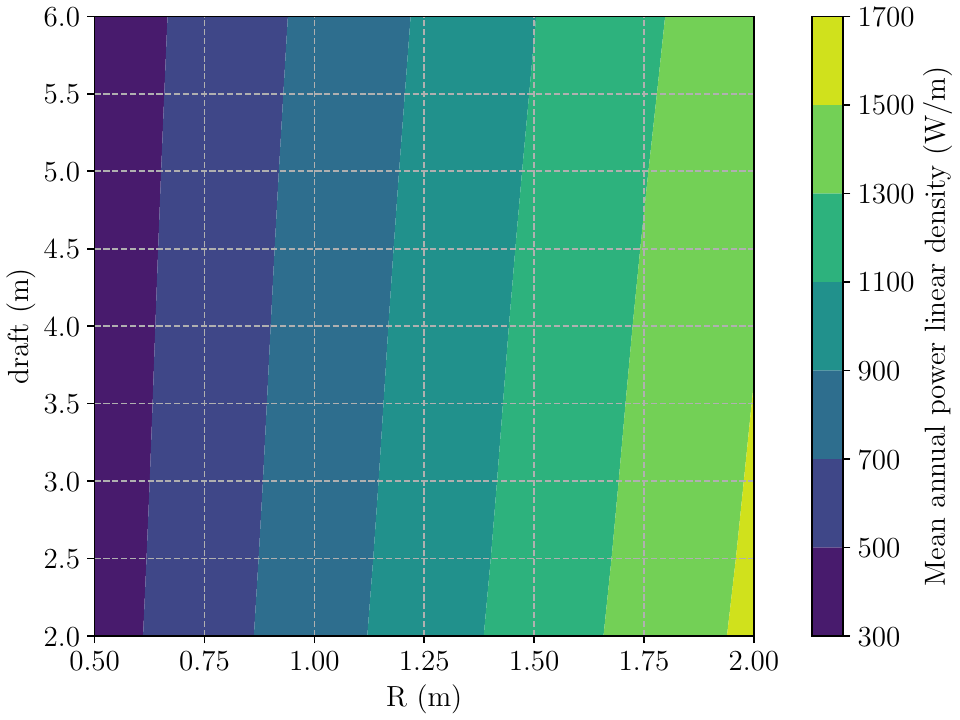}
}
\subfloat{
\includegraphics[width=0.32\textwidth]{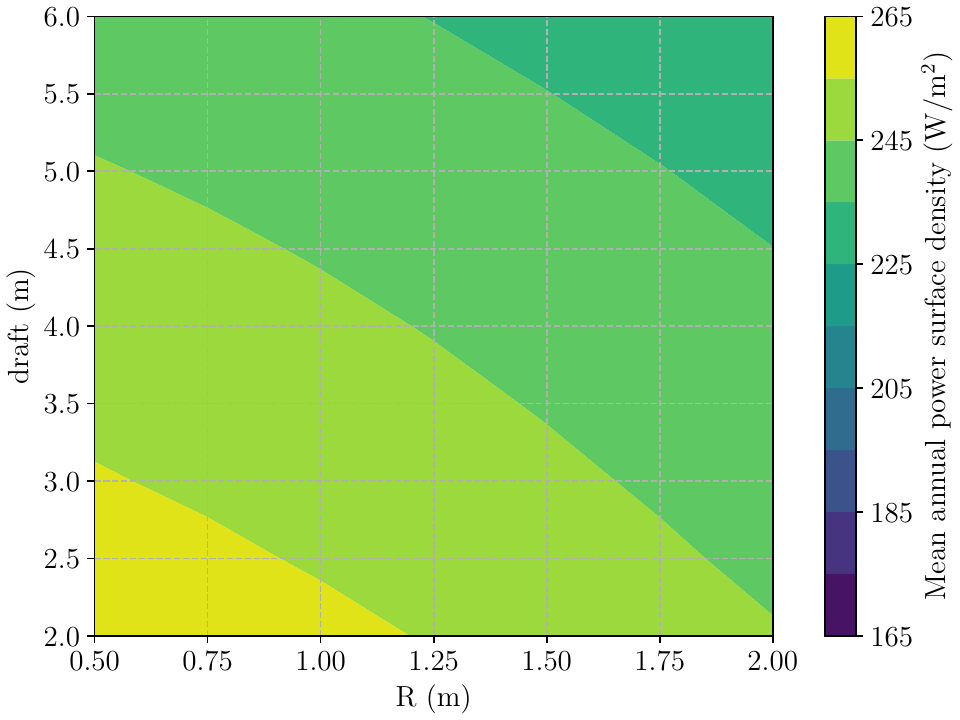}
}
\caption{Mean annual power for variable geometric parameters; wave climate of Alghero (unconstrained model)}
\label{fig:annpower-alghero}
\end{figure}

\begin{figure}[h!]
\centering
\subfloat{
\includegraphics[width=0.32\textwidth]{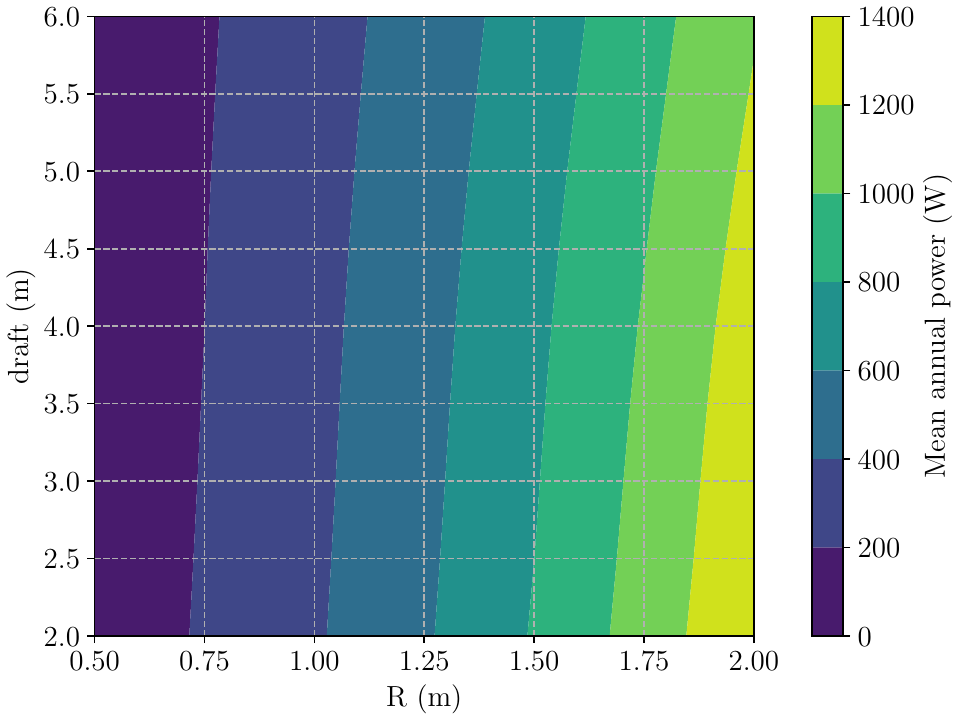}
}
\subfloat{
\includegraphics[width=0.32\textwidth]{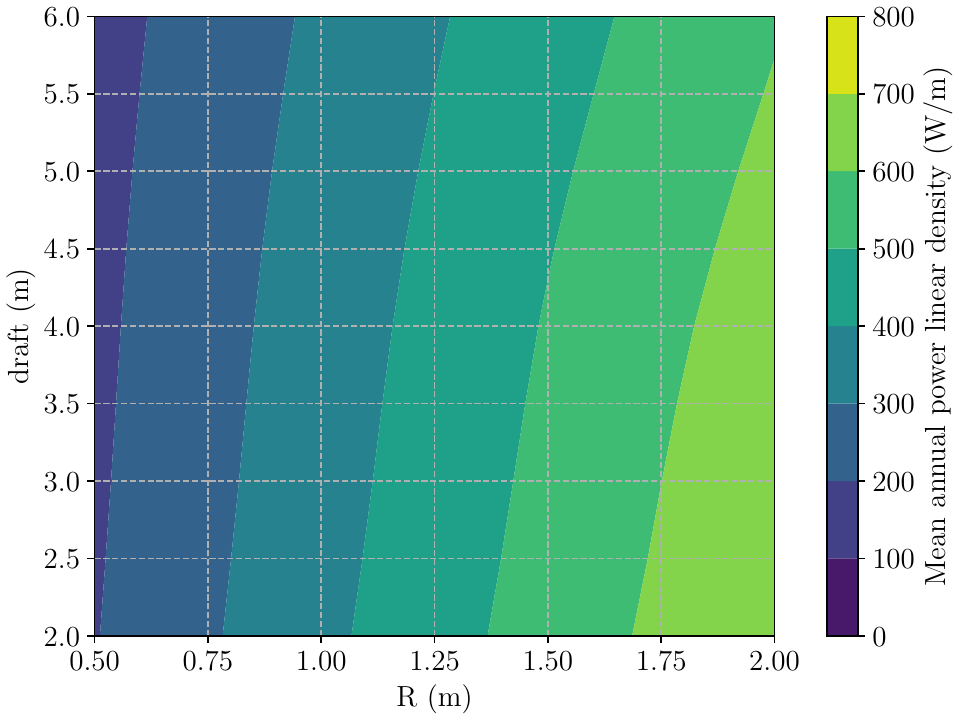}
}
\subfloat{
\includegraphics[width=0.32\textwidth]{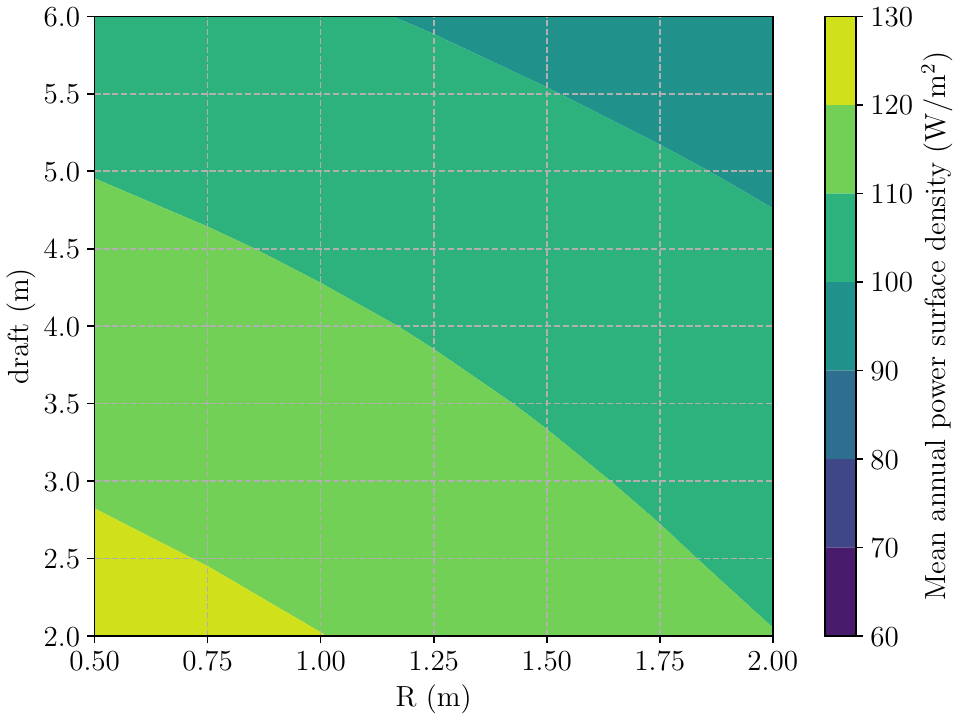}
}
\caption{Mean annual power for variable geometric parameters; wave climate of Civitavecchia (unconstrained model)}
\label{fig:annpower-civitavecchia}
\end{figure}

The mean annual power of a device installed at a specific location can be computed by multiplying term by term the power matrix by the matrix of relative occurrences of sea states at such location, and summing. 
The latter matrix is indicated in the literature as \emph{scatter diagram of wave statistics} \cite{Babarit2012} or \emph{resource characterization matrix} \cite{Bozzi2018}. Given the resource characterization matrix of a location, the mean annual power of a device will depend on the dimensions of the device and on the control strategy. In the following, we present a computationally cheap procedure for a parameter exploration on device dimensions. 
For simplicity, the linear model is used and the control strategy is set to an approximate optimum: we impose that, in each sea state, the maximum magnitude of the flow coefficient (that is the absolute value of its complex amplitude) corresponds to the condition of maximum nondimensional torque. 
This is equivalent to requiring that at each wave cycle the turbine blades reach the condition of incipient stall without exceeding it. We write such condition as $|\widehat{\varphi}| = \varphi_\text{opt}$. Mass conservation requires
\begin{equation}
S \omega |\widehat{\zeta}| = S_t |\widehat{v}_t|.
\end{equation}
The conditions above, together with the definition of the flow coefficient \eqref{eq:defflowcoeff}, yield
\begin{equation}
|\widehat{\zeta}| = \frac{|\widehat{v}_t| S_t}{\omega S} = \frac{\varphi_\text{opt} S_t r_t \omega_t}{\omega S}.
\end{equation} 
It is possible to determine the value of $\omega_t$ by solving the nonlinear equation
\begin{equation}
\left\{ [-\omega^2(A(\omega) + \rho L) + \rho g]^2 + \omega^2(B(\omega) + S\Lambda(\omega_t))^2 \right\}  \frac{\varphi_\text{opt} S_t r_t \omega_t}{\omega S} - |\widehat{p}_e|^2 = 0
\end{equation}
using Newton's method. The results obtained from the wave climates of two locations, Alghero and Civitavecchia, are reported in Fig.~\ref{fig:annpower-alghero}, \ref{fig:annpower-civitavecchia}.
The presented method has a very low computational cost, but it does not take into account the constraints of immersed turbine and absence of cavitation. For this reason, the obtained results may overestimate the power for some combinations of the parameters. We thus present an alternative method, in which constraints are enforced by penalization. 

A generic constraint will be expressed as $c_i \leq 0$. In the frequency domain, the constraint of immersed turbine can be expressed as
\begin{equation}
c_1 = |\widehat{\zeta}| - |z_t|.
\end{equation}
The vertical coordinate $z_t$ of the turbine is fixed at 60\% of the draft, indicated as $d$: $z_t = -0.6 d $.
For the constraint of no cavitation, we introduce a cautionary estimate of the minimum pressure on blades. Starting from \eqref{eq:pcavit} and defining the minimum over $\varphi$ of the minimum pressure coefficient
\begin{equation}
\tilde{C}_{p, \text{min}} = \min_\varphi C_{p, \text{min}}(\varphi),
\end{equation}
one obtains the lower bound $\tilde{p}_\text{min}$:
\begin{equation}
p_\text{min}(t) \geq \tilde{p}_\text{min} (t) \equiv p_\text{atm} + \rho g (\zeta(t) - z_t) + \dfrac{1}{2} \rho \left( v_t^2(t) + \omega_t^2 r_t^2 \right) \tilde{C}_{p,\text{min}}.
\end{equation}
We now express the time evolution of the level of the water column and of the axial velocity from their complex amplitudes (their absolute phases are irrelevant for our purposes):
\begin{equation}
\zeta(t) = |\widehat{\zeta}| \cos(\omega t),
\quad
v_t(t) = \dot{\zeta}(t) \dfrac{S}{S_t} = - \frac{\omega |\widehat{\zeta}| S}{S_t} \sin(\omega t).
\end{equation}
It is possible to separate a constant pressure term\footnote{This is not the mean pressure, since the term $v_t^2$, that is included in the fluctuations, does not have zero mean.} $\overline{p}$ from fluctuations $\tilde{p}'$:
\begin{equation}
\overline{p} = p_\text{atm} - \rho g z_t + \dfrac{1}{2} \rho \omega_t^2 r_t^2 \tilde{C}_{p, \text{min}},
\end{equation}
and compute the time derivative of the fluctuations:
\begin{equation}
\tilde{p}'_{min}(t) = \sin(\omega t) \left[ -\rho g |\widehat{\zeta}| \omega + \dfrac{\rho \omega^3 |\widehat{\zeta}|^2 S^2}{S_t^2} \tilde{C}_{p, \text{min}} \cos(\omega t) \right].
\label{eq:ptildemin}
\end{equation}
The first pair of stationary points occurs for $\sin(\omega t)=0$. It corresponds to  $\tilde{p}_\text{min}(t) = \overline{p} \pm \rho g |\widehat{\zeta}| $, i.e., a local maximum and a local minimum. We take the minimum and introduce the constraint function
\begin{equation}
c_2 = -\overline{p} + \rho g |\widehat{\zeta}| + p_v.
\end{equation}
An additional stationary point exists if the square bracket vanishes in \eqref{eq:ptildemin}. Since $\cos(\omega t) \geq -1$ and $\tilde{C}_{p, \text{min}}<0$, this can happen only if 
\begin{equation}
\tilde{C}_{p, \text{min}} \leq - \dfrac{g S_t^2}{S^2 \omega^2 |\widehat{\zeta}|}.
\label{eq:ptildemincond}
\end{equation}
We thus introduce the constraint
\begin{equation}
c_3 = \begin{cases}
0, & \text{if \eqref{eq:ptildemincond} holds,}\\
-\overline{p} - \dfrac{\rho g^2 S_t^2}{\omega^2 S^2 \tilde{C}_{p,min}} - \dfrac{1}{2} \rho \dfrac{|\widehat{\zeta}|}{S_t^2} \sqrt{(\omega^2 |\widehat{\zeta}|S^2 \tilde{C}_{p, min})^2 - (g S_t^2)^2} + p_v & \text{otherwise.}
\end{cases}
\end{equation}
Finally, for modeling reasons, we impose that the flow coefficient does not exceed the limits of validity of the polynomial approximation of the nondimensional torque curve \eqref{eq:torquepoly}. The corresponding constraint function is
\begin{equation}
c_4 = |\widehat{\varphi}| - \varphi_\text{max, model}.
\end{equation}
All constraints are enforced by a quadratic penalty function \cite{Nocedal2006}:
\begin{equation}
Q_\mu(\omega_t) = -P(\omega_t) + \mu \sum_{i=1}^4 [c_i(\omega_t)]_+^2.
\end{equation}
This is the cost function to be minimized; the single variable optimization problem is solved, again, using function \texttt{minimize\_scalar}. Convergence of the quadratic penalty method to the true constrained minimum is guaranteed in the limit $\mu \to \infty$; thus, optimization is repeated for increasing values of $\mu$ until a tolerance on the violation of constraints is satisfied.

\begin{figure}[h!]
\centering
\subfloat{
\includegraphics[width=0.32\textwidth]{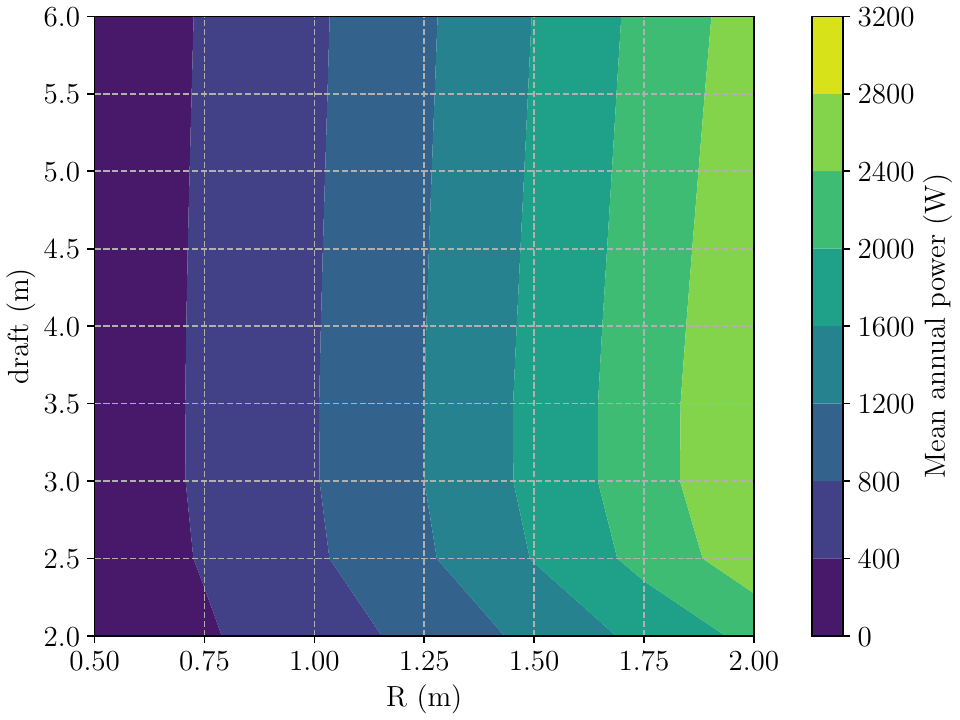}
}
\subfloat{
\includegraphics[width=0.32\textwidth]{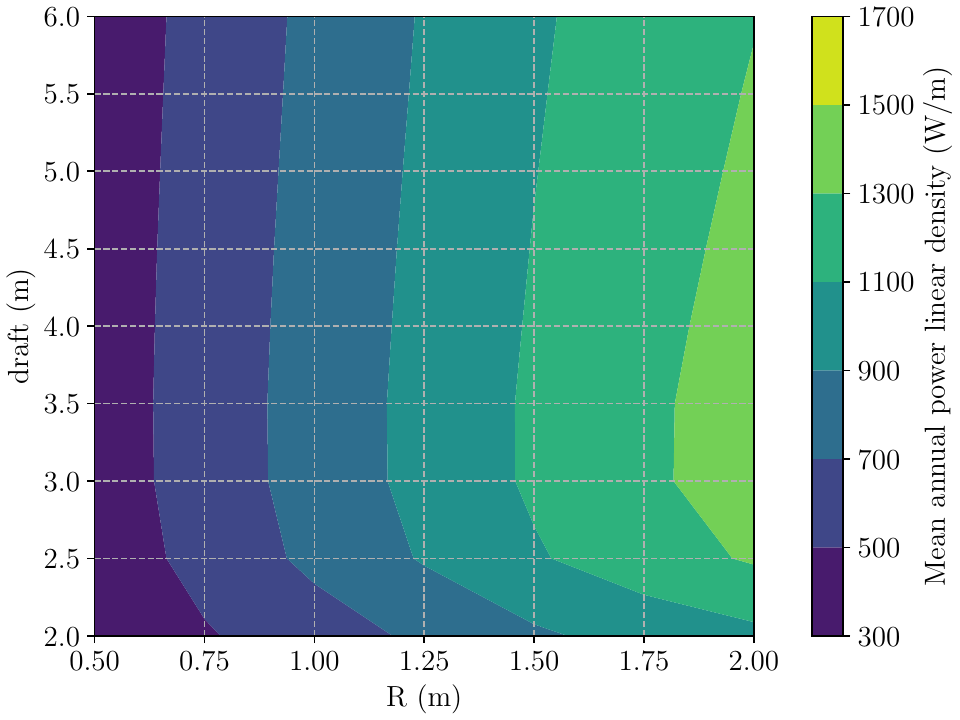}
}
\subfloat{
\includegraphics[width=0.32\textwidth]{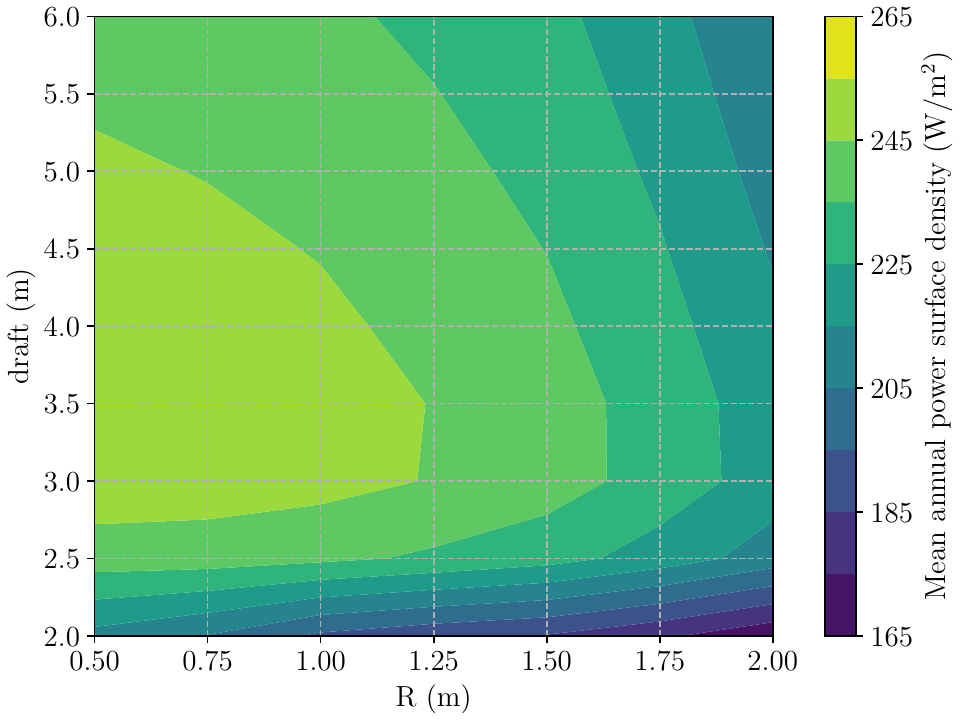}
}
\caption{Mean annual power for variable geometric parameters; wave climate of Alghero (constrained model)}
\label{fig:annpower-alghero-constr}
\end{figure}

\begin{figure}[h!]
\centering
\subfloat{
\includegraphics[width=0.32\textwidth]{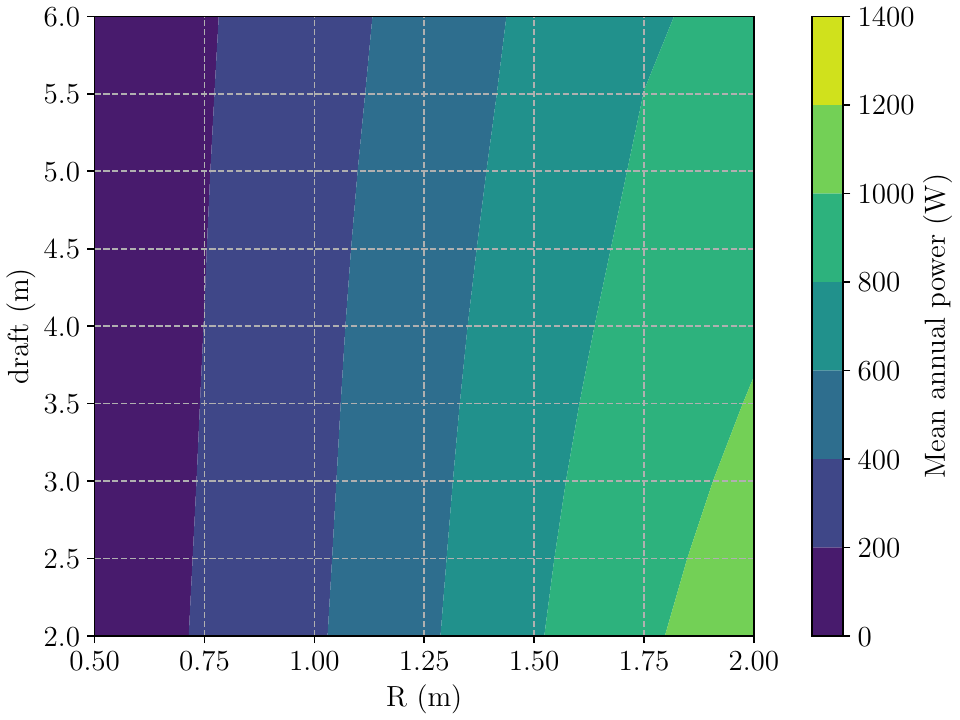}
}
\subfloat{
\includegraphics[width=0.32\textwidth]{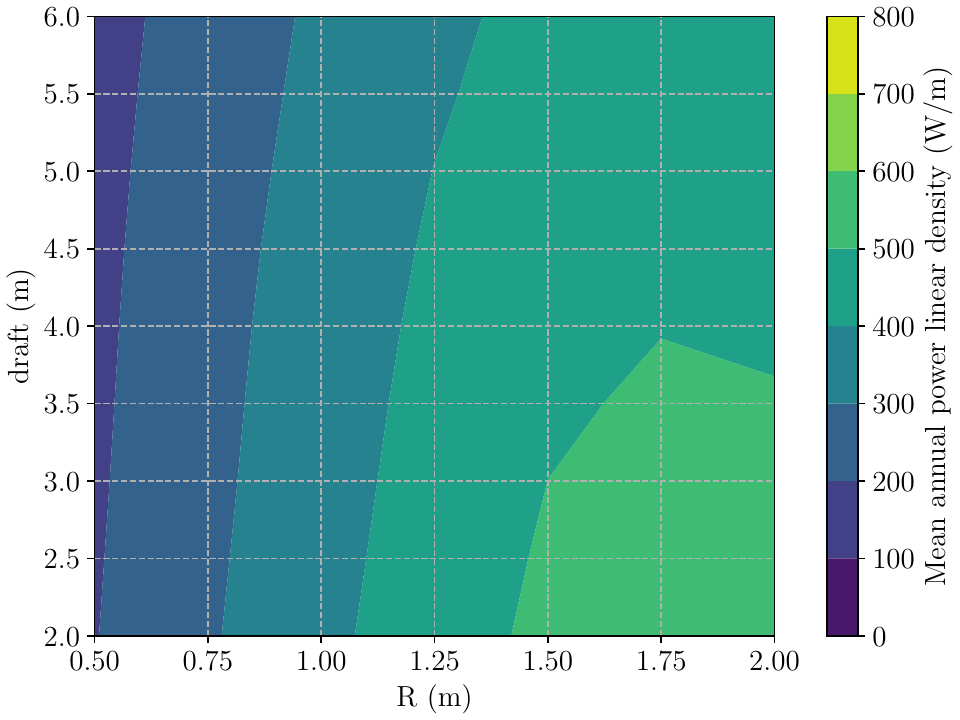}
}
\subfloat{
\includegraphics[width=0.32\textwidth]{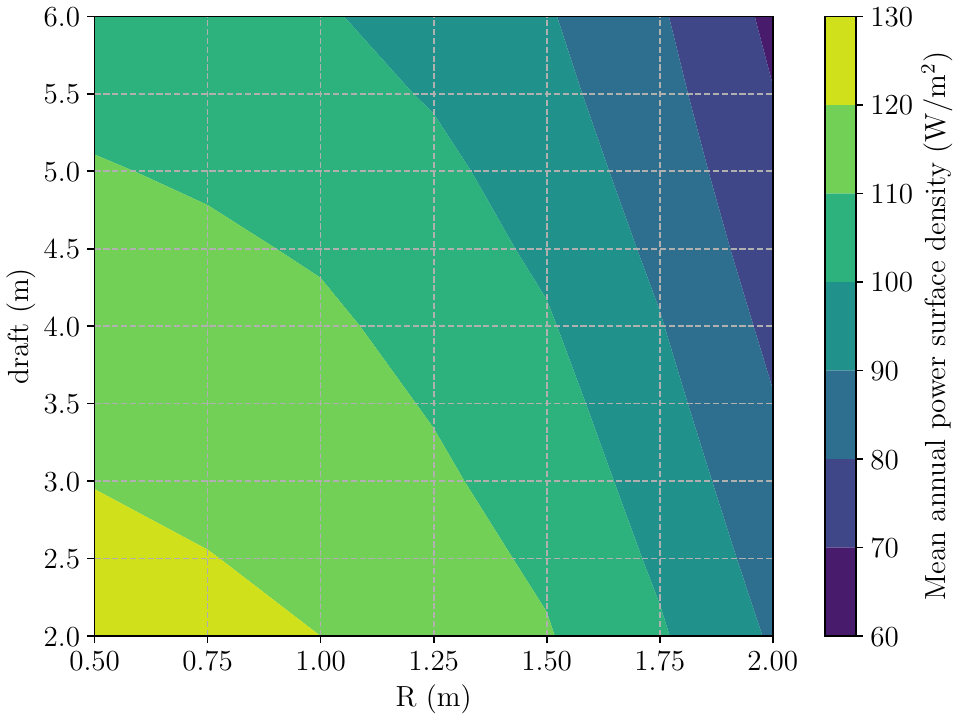}
}
\caption{Mean annual power for variable geometric parameters; wave climate of Civitavecchia (constrained model)}
\label{fig:annpower-civitavecchia-constr}
\end{figure}

Results are shown in Fig.~\ref{fig:annpower-alghero-constr}, \ref{fig:annpower-civitavecchia-constr}.
The results of the constrained model for the wave climate of Alghero, more energetic than the one of Civitavecchia, are especially different from the ones obtained with the unconstrained model. Indeed, in the unconstrained case, the oscillation amplitude of the water column obtained in the conditions of maximum turbine efficiency corresponds to a violation of the immersion constraint for the geometries with smallest drafts. Optimal conditions are thus obtained for intermediate values of the draft. The wave climate of Civitavecchia is instead not energetic enough for this effect to be noticeable.

The results are reported in terms of individual power of a device, linear power density and surface power density. They indicate that a larger device would be more efficient in case of installation in single units or in linear arrays, while smaller devices would yield better performance per unit area if, for instance, an installation with multiple rows is devised. The choice of a specific dimension should be based on the location of installation, on mechanical considerations on the feasibility of large turbines, and on economic considerations on the cost of a single device and its dependence on scale. Caution needs to be applied in interpreting the results of this section, since diffraction and radiation become important even for very small devices, if they are closely packed. 

The result that small devices may yield a larger power per unit area may not seem intuitively convincing. For this reason, a simple computation for devices that are small compared to the wavelength is reported here. In this approximation, for a single device, we can neglect diffraction and radiation. The excitation force is then only due to incident waves, whose potential is \eqref{eq:incwavepot}. If the device is a small floating cylinder of radius $r$ and draft $d$, then the potential can be considered uniform on its bottom surface, situated at $z = -d$. If the deep water limit is considered, then $\cosh[k(z+h)]/\cosh(kh) \approx \exp(kz)$. We thus have
\begin{equation}
m = \rho \pi r^2 d, \quad k = \rho g \pi r^2, \quad |\widehat{f}_e| \approx \rho g \frac{H}{2} \pi r^2 \exp(-kd).
\end{equation}
Since radiation is negligible, the dynamic equation in the frequency domain with damping control coefficient $c$ reduces to
\begin{equation}
(-\omega^2 m - i \omega c + k) \widehat{\zeta} = \widehat{f}_e.
\end{equation}
It can be verified that maximum average power is obtained for 
\begin{equation}
c_\text{opt} = \frac{|-\omega^2 m + k|}{\omega}, \quad |\widehat{\zeta}|^2 = \frac{|\widehat{f}_e|^2}{2 (-\omega^2 m + k)^2}.
\end{equation}
By substitution, we obtain the corresponding value of power:
\begin{equation}
P_\text{opt} = \frac{1}{2} \omega^2 c |\widehat{\zeta}|^2 = \frac{\rho g \omega H^2}{16} \frac{\pi r^2}{1 - kd} \exp(-2kd) \approx \frac{\rho g \omega H^2}{16} \pi r^2 (1 - kd).
\end{equation}
The last expression gives at least a partial explanation of the obtained results. It predicts a quadratic increase of power with respect to the radius, which corresponds to a linear increase of power per unit width. This is quite well reproduced by the numerical results. It also predicts that power per unit area is independent on the radius and weakly decreasing with the draft. The numerical results for the unconstrained case do, indeed, show a decrease with the draft; however, an increase with decreasing radius is also observed. In the constrained case, the dependence on draft changes, as already observed, because the constraint of immersed turbine is enforced. 

\section{Arrays of devices}\label{sec:multibodywavesax}
In this section, we analyze configurations with multiple devices,  considering in particular the integration with floating offshore wind energy plants. We consider the mutual hydrodynamic interactions between different devices, and the interactions with the piles of the wind energy platform. Geometry and dimensions are consistent with the ones of the Windfloat project\footnote{\url{https://www.boem.gov/sites/default/files/about-boem/BOEM-Regions/Pacific-Region/Renewable-Energy/11-Kevin-Bannister---BOEM-Workshop.pdf}.}.

\subsection{Mathematical model} 
\label{sec:parkmodel}
Park hydrodynamics are described using the semi-analytical model proposed in \cite{Yilmaz1998, Child2010}. The model is based on separation of variables for the Laplace equation in cylindrical coordinates, resulting in an infinite series which needs to be truncated for purposes of numerical computation. Only cylindrical objects can be simulated. The hydrodynamic behavior of each cylinder is described by a diffraction transfer matrix, which is independent of the position and can hence be computed once and for all. The interactions between cylinders are described by basis transformation matrices, which instead depend on the positions, and whose analytical expressions can be explicitly differentiated with respect to the coordinates of the center of mass of each device. Truncation of the series leads to the linear, block system
\begin{equation}
\begin{bmatrix}
M_{\gamma\gamma} & M_{\gamma\zeta} \\
M_{\zeta\gamma} & M_{\zeta\zeta} 
\end{bmatrix}
\begin{bmatrix}
\widehat{\bm{\gamma}} \\ \widehat{\bm{\zeta}}
\end{bmatrix}
=
\begin{bmatrix}
\bm{h}_1 \\ \bm{h}_2
\end{bmatrix},
\label{eq:stateprob}
\end{equation}
whose first block-row is the hydrodynamic problem and the second block-row describes the dynamics of the bodies. $\widehat{\bm{\gamma}}$ is the vector of coefficients of the series describing the hydrodynamic solution. All blocks of matrix $M$ and of right hand side vector $\bm{h}$ depend on the coordinates of the devices.

\subsection{Simulation and optimization}
\subsubsection{Interaction between a single device and a platform}\label{sec:1devpiles}
\begin{figure}[h!]
\centering
\subfloat{
\begin{tikzpicture}
\node (image) at (0, 0)
{\includegraphics[width=0.45\textwidth, trim= 0.5cm 0 0 0, clip]{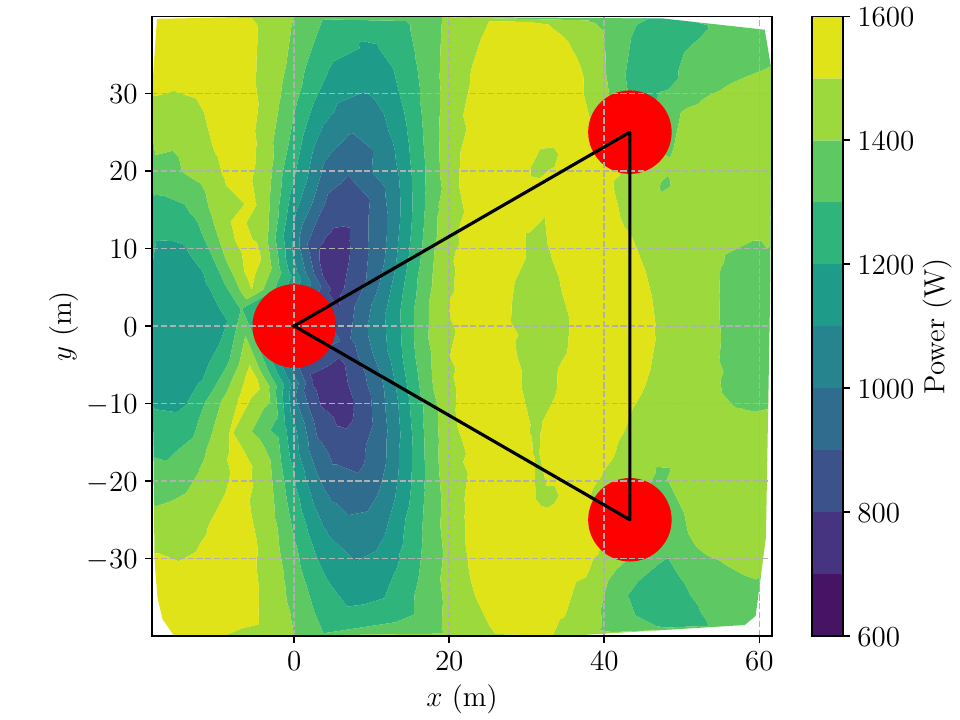}};
\draw[->, line width=1.5](-2.35,0.25)--(-1.85,0.25);
\end{tikzpicture}
}
\hfill
\subfloat{
\begin{tikzpicture}
\node (image) at (0, 0) {\includegraphics[width=0.45\textwidth, trim= 0.5cm 0 0 0, clip]{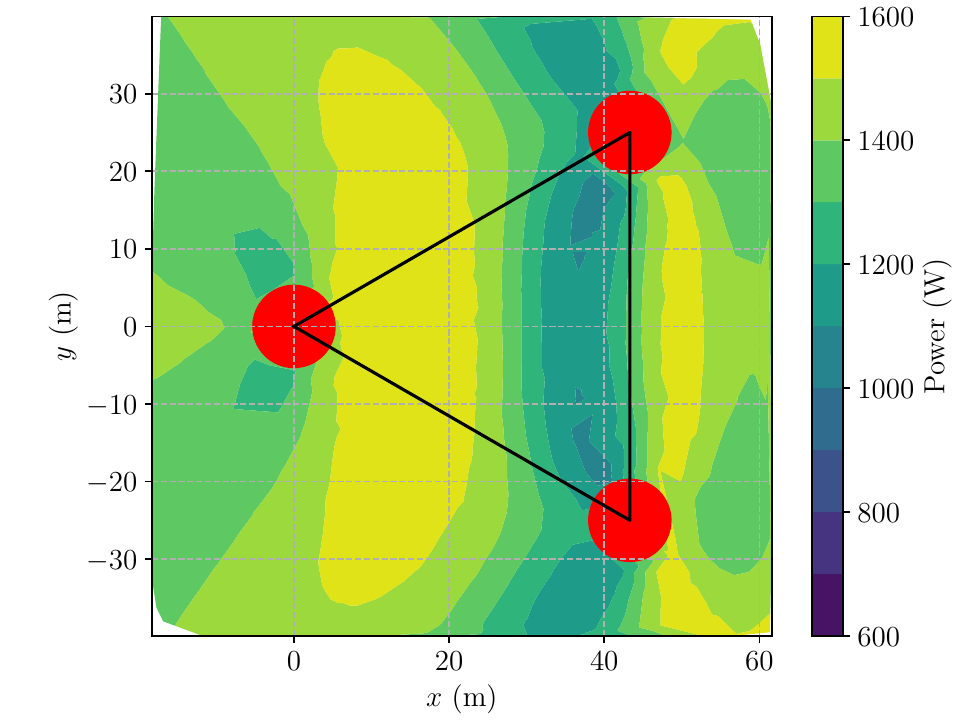}};
\draw[->, line width=1.5](1.8,0.25)--(1.3,0.25);
\end{tikzpicture}
}
\caption{
Map of powers for a single device interacting with the piles, for waves coming respectively from the left and from the right, $H_s$ = 3 m, $T$ = 8 s}
\label{fig:powermap1bodypiles}
\end{figure}

We first analyze the effect of the presence of the piles of the structure on the performance of a single device. The motion of the piles is assumed to be negligible. In order to realize this effect, the piles are assigned very large values of stiffness. This can be interpreted as the equivalent stiffness of the moorings.
A more accurate approach, which has not been explored in the present work, would require taking into account the forces acting on the wind turbine, in order to model the displacement of the floating platform on which it is mounted. 
Simulations are carried out with a wave condition typical of the Mediterranean as forcing term, and with the piles located at the vertices of an equilateral triangle. With reference to Fig.~\ref{fig:powermap1bodypiles}, we consider horizontal wave vectors, in both possible directions.

In the case of waves coming from the left, and considering only the region inside the triangle, an average power of 1470~W is obtained, while considering the whole square the average power is 1364~W. In the case of waves coming from the right, the average power in the internal region is 1312~W, while in the whole square it is 1396~W. The maximum power is 1540~W in both cases. Thus, if only the internal region is available, then the first configuration is more efficient; if, instead, also the external region in the proximity of the piles can be utilized, then the second configuration is better.
It has to be taken into account that the choice might be dictated by the necessities of the wind power plant: typically, the turbine is positioned on the downwind vertex of the platform, with respect to the dominant wind direction. In this case, in a location with waves mainly generated by the local wind and thus with wind vector approximately aligned with the wave vector, the second configuration would be chosen.

\subsubsection{Interaction between 2 devices}\label{sec:2dev}
\begin{figure}[h!]
\centering
\subfloat{
\begin{tikzpicture}
\node (image) at (0, 0) {
\includegraphics[width=0.45\textwidth, trim= 0 0 0 0, clip]{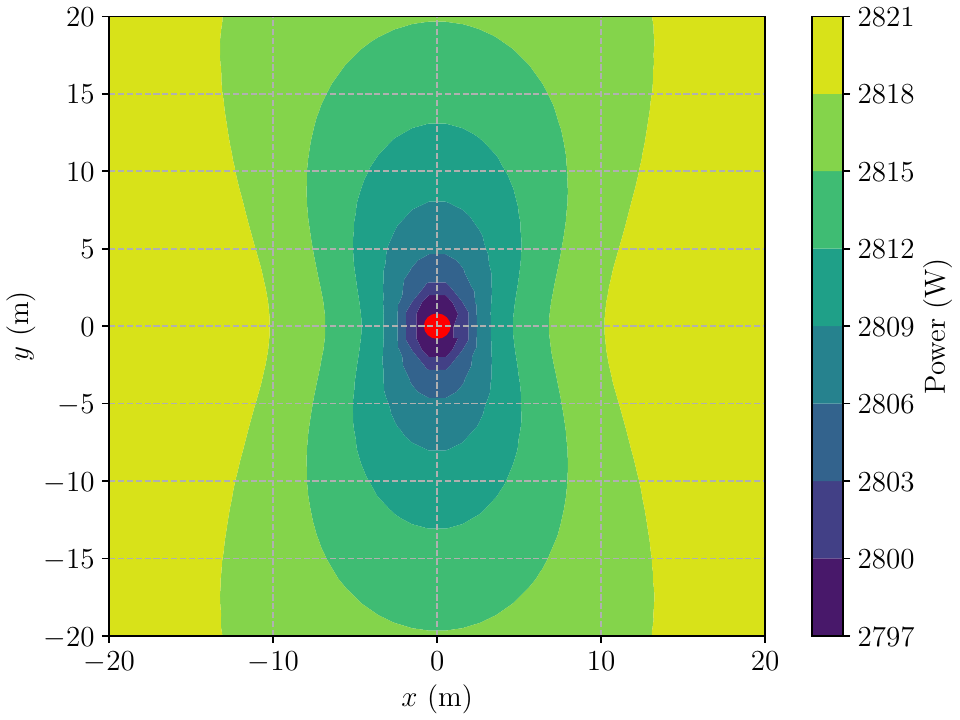}};
\draw[->, line width=1.5] (-2.,0.25)--(-1.5,0.25);
\end{tikzpicture}
}
\hfill
\subfloat{
\begin{tikzpicture}
\node (image) at (0, 0){
\includegraphics[width=0.45\textwidth, trim= 0 0 0 0, clip]{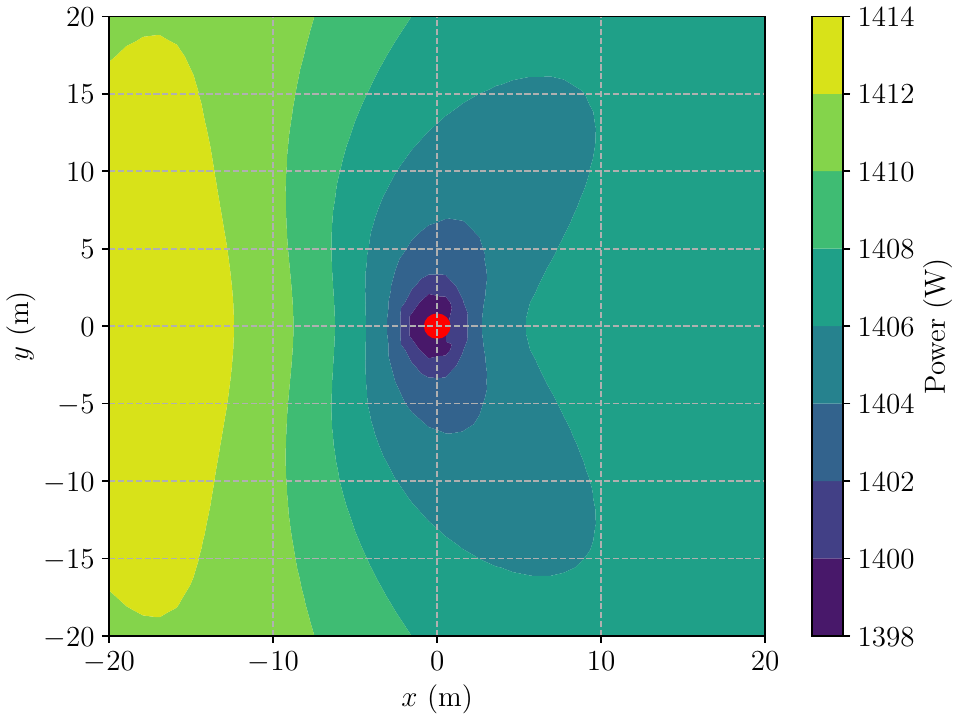}};
\draw[->, line width=1.5] (-2.,0.25)--(-1.5,0.25);
\end{tikzpicture}
}
\caption{
Array of 2 bodies, one in the origin and the other in $(x, y)$. Left: total power; right: power of the device in $(x, y)$}
\label{fig:powermap2bodies}
\end{figure}

As a complementary analysis with respect to the one of the previous section, we consider two interacting devices, in absence of piles. 
The aim is to quantify the effect of interaction while varying the distance. The results are shown in Fig.~\ref{fig:powermap2bodies}, with waves coming from the left. 
It can be observed that the upwave body gains an advantage thanks to the presence of the downwave body. The most unfavourable conditions are the ones with minimum distance between the two devices. The quantitative power variation may seem to indicate that the interaction effect is very small (less than 1\%). 
It has however been observed in the literature that the so-called park effect is relevant only when the number of devices is large enough (indicatively more than 10) \cite{Babarit2015}. We thus postpone this discussion to Sec.~\ref{sec:optresults}. Some further comparisons with 3D CFD simulations are reported in \cite{Agate2023}.

\subsubsection{Optimization method} \label{sec:optmethod}
The cost function is the power, with flipped sign:
\begin{equation}
J(\widehat{\bm{\zeta}}) = -P(\widehat{\bm{\zeta}}) = -\sum_{\ell=1}^{N_b} \sum_{n=0}^d p_n \left( \omega |\widehat{\zeta}_\ell| \right)^{2n}.
\label{eq:multibodycost}
\end{equation}
It must be minimized taking into account the constraints given by the state problem, the available domain and the requirement that devices do not overlap. 

\begin{figure}[h!]
\centering
\subfloat{
\begin{tikzpicture}
\definecolor{lightblue}{RGB}{173,216,230}
\path (0., 0.) coordinate (O);
\path (4.33, 2.5) coordinate (A);
\path (4.33, -2.5) coordinate (B);
\draw (O) -- (A) -- (B) -- (O);
\draw [fill=red, fill opacity=0.6] (O) circle (0.5);
\draw [fill=red, fill opacity=0.6] (A) circle (0.5);
\draw [fill=red, fill opacity=0.6] (B) circle (0.5);
\path (0.65, 0.202) coordinate (P1);
\path (3.83, 2.04) coordinate (P2);
\path (4.18, 1.84) coordinate (P3);
\path (4.18, -1.84) coordinate (P4);
\path (3.83, -2.04) coordinate (P5);
\path (0.65, -0.202) coordinate (P6);
\draw [fill=lightblue, dashed] (P1) -- (P2) -- (P3) -- (P4) -- (P5) -- (P6) -- (P1);
\end{tikzpicture}}
\hfill
\subfloat{
\begin{tikzpicture}
\definecolor{lightblue}{RGB}{173,216,230}
\draw [fill=red, fill opacity=0.6] (O) circle (2.125);
\path (0., 0.) coordinate (O);
\path (2.39, -1.38) coordinate(Q);
\path (4.42, 2.55) coordinate (A);
\path (4.42, -2.55) coordinate (B);
\path (2.76, 0.86) coordinate (P1);
\path (2.76, -0.86) coordinate (P6);
\path (4.42, 1.82) coordinate (P2);
\path (4.42, -1.82) coordinate (P5);
\draw [dashed] (Q) arc (-30:30:2.76);
\draw (O) -- (A);
\draw (O) -- (B);
\draw [fill=lightblue, dashed] (P5) -- (P6) -- (P1) -- (P2);
\end{tikzpicture}}
\caption{Complete convexified domain (left) and detail of the comparison with the original domain (right)}
\label{fig:convexdomain}
\end{figure}
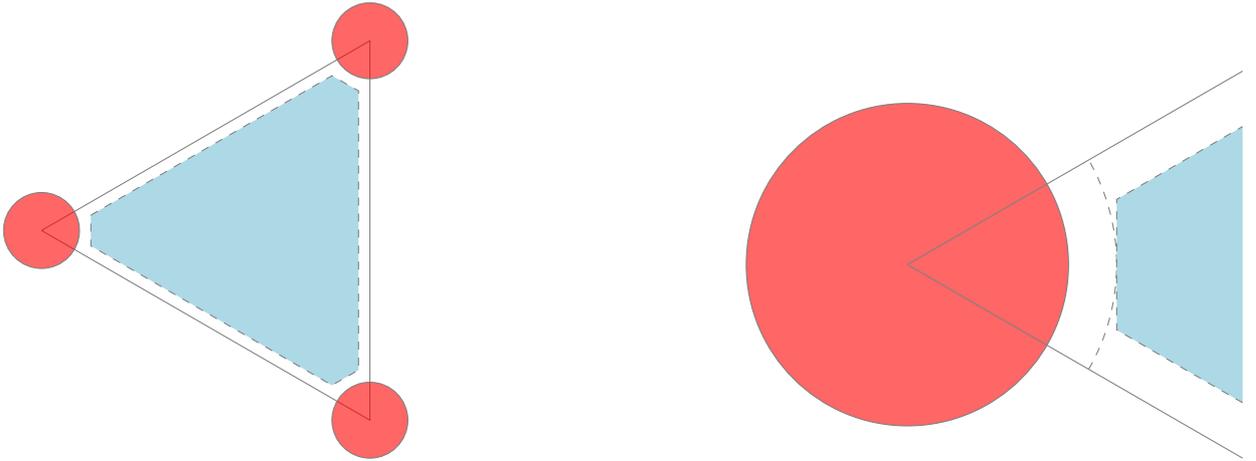

We consider as available area the interior of the triangle defined by the centers of the piles supporting the platform (see Fig.~\ref{fig:convexdomain}). 
The devices' centers can be placed anywhere in the light blue region bounded by the dashed line, so that their walls are at most tangent to the boundary defined by the solid line and the piles. Notice that, in addition, the arcs corresponding to the piles are substituted by their tangents in order to make the domain convex, thus allowing the use of a projected gradient algorithm for the enforcement of the constraint.

For the solutions to be physically meaningful, it is also necessary to impose that devices do not overlap. A minimum distance $d_\text{min}$ between the centers of the devices is imposed. This constraint corresponds to a non-convex admissible set. To see this, consider two devices positioned respectively in $(x_1, y_1)$ and $(x_2, y_2)$. The vector of optimization variables, which contains the coordinates of the centers of all devices, is $\bm{u} = (x_1, y_1, x_2, y_2) \in \mathbb{R}^4$. Two admissible configurations are $\bm{u}_1 = (0, 0, d_\text{min}, 0)$ and $\bm{u}_2 = (0, 0, 0, d_\text{min})$. Their convex combination with parameter 1/2 is $\bm{u}_3 = (0, 0, d_\text{min}/2, d_\text{min}/2)$, corresponding to a distance between the centers equal to $d_\text{min} \sqrt{2}$. Hence, $\bm{u}_3$ does not satisfy the constraint of no overlap and the admissible set is not convex.

The state equation is enforced by a Lagrangian approach to find a reduced gradient. Optimization is performed by employing a projected gradient algorithm \cite{Nocedal2006}, where projection is onto the convex admissible domain. A modified backtracking line search is executed at each iteration. It is based on the Armijo rule, with the additional requirement that steps are small enough to prevent the devices from colliding. 

To find the reduced gradient of the cost function, we introduce multiplier vectors $\bm{\lambda}$, $\bm{\mu}$ to enforce the state equation \eqref{eq:stateprob}, and define the Lagrangian function
\begin{equation}
\mathcal{L}(\bm{x}; \widehat{\bm{\gamma}},\widehat{\bm{\zeta}}; \bm{\lambda}, \bm{\mu}) = J(\widehat{\bm{\zeta}}) + \Re\left[\bm{\lambda}^H(M_{\gamma\gamma} \widehat{\bm{\gamma}} + M_{\gamma\zeta}\widehat{\bm{\zeta}} - \bm{h}_1) + \bm{\mu}^H(M_{\zeta\gamma} \widehat{\bm{\gamma}} + M_{\zeta\zeta} \widehat{\bm{\zeta}} - \bm{h}_2) \right].
\end{equation}
Differentiation with respect to the state variables $\widehat{\bm{\gamma}},\widehat{\bm{\zeta}}$
\begin{align}
D_\gamma \mathcal{L}(\delta \widehat{\bm{\gamma}}) &= \Re\left[ \delta \widehat{\bm{\gamma}}^H (M_{\gamma\gamma}^H \bm{\lambda} + M_{\zeta\gamma}^H \bm{\mu}) \right], \\
D_\zeta \mathcal{L}(\delta \widehat{\bm{\zeta}}) &= -\sum_{\ell=1}^{N_b} \sum_{n=0}^d 2 n p_n \Re[\delta \widehat{\zeta}_\ell^* (\widehat{\zeta}_\ell^*)^{n-1} \widehat{\zeta}_\ell^n] + \Re[\delta \widehat{\bm{\zeta}}^H (M_{\gamma\zeta}^H \bm{\lambda} + M_{\zeta\zeta}^H \bm{\mu})],
\end{align}
leads to the adjoint equations
\begin{equation}
\begin{bmatrix}
M_{\gamma\gamma}^H & M_{\zeta\gamma}^H \\[0.5em]
M_{\gamma\zeta}^H & M_{\zeta\zeta}^H
\end{bmatrix}
\begin{bmatrix}
\bm{\lambda} \\ \bm{\mu}
\end{bmatrix}
=
\begin{bmatrix}
0 \\ \widetilde{\bm{h}}_2
\end{bmatrix}, \quad \widetilde{h}_{2,\ell} = \sum_{n=0}^d 2 n p_n  (\widehat{\zeta}_\ell^*)^{n-1} \widehat{\zeta}_\ell^n.
\label{eq:adjprob}
\end{equation}
Differentiation with respect to the positions yields the components of the gradient $\nabla J$ with respect to the coordinates of devices:
\begin{equation}
\frac{\partial \mathcal{L}}{\partial x_\ell} = \Re \left[ \bm{\lambda}^H \left( \frac{\partial M_{\gamma\gamma}}{\partial x_\ell} \widehat{\bm{\gamma}} + \frac{\partial M_{\gamma\zeta}}{\partial x_\ell} \widehat{\bm{\zeta}} - \frac{\partial \bm{h}_1}{\partial x_\ell}
\right) + \bm{\mu}^H \left( \frac{\partial M_{\zeta\gamma}}{\partial x_\ell} \widehat{\bm{\gamma}} + \frac{\partial M_{\zeta\zeta}}{\partial x_\ell} \widehat{\bm{\zeta}} - \frac{\partial \bm{h}_2}{\partial x_\ell}
\right) \right],
\label{eq:gradcostpark}
\end{equation}
where derivatives of all blocks in the system can be computed analytically, as detailed in \cite{jacopo}.

\begin{figure}
\centering
\includegraphics[width=0.4\textwidth]{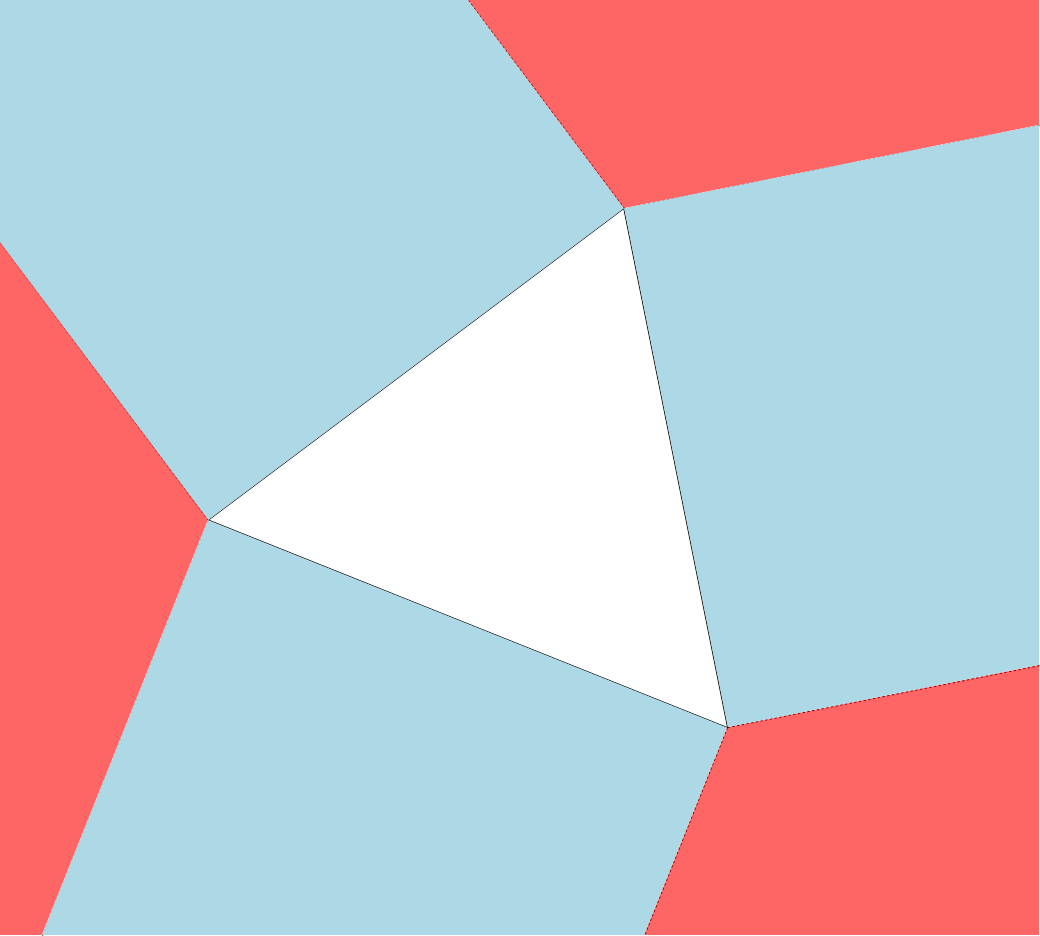}
\caption{Graphical representation of the projection algorithm}
\label{fig:projalg}
\end{figure}
The definition of the projector onto the admissible domain is based on the observation of Fig.~\ref{fig:projalg}. Consider a point $\bm{x}$ outside the white triangle (in general convex polygon) $\Omega$: the point of $\Omega$ at the minimum distance from $\bm{x}$ is either a vertex (if $\bm{x}$ is in one of the red regions) or it is on an edge (if $\bm{x}$ is in one of the light blue regions). To each edge of the polygon we can associate an inequality subconstraint, corresponding to belonging to an half-plane. The constraint of belonging to the polygon will then be given by the system of such subconstraints. 
If the orthogonal projection of $\bm{x}$ onto one of the lines corresponding to the violated subconstraints belongs to $\Omega$, then it is the projection of $\bm{x}$ onto $\Omega$ (light blue regions). Otherwise, the projection will be the vertex closest to $\bm{x}$ (red regions). We will denote with $\mathcal{P}(\bm{x})$ the projection of a point $\bm{x}$ onto $\Omega$, computed as described above.

Since the constraint of no overlap is not convex, it needs to be treated differently. It is taken into account during backtracking line search. If the distance between two devices is less than $d_\text{min}$, the step is reduced until an admissible condition is obtained. This corresponds to substituting the gradient $\nabla J$ with a vector $D \nabla J$, where $D$ is a diagonal matrix with nonnegative elements. A descent direction is obtained, since $D$ is semi-positive definite. An alternative approach could be penalization. In such case, however, it would be required to compute the cost in non-admissible configurations, in which the hydrodynamic problem cannot be solved.

Algorithm~\ref{alg:optwavesaxpark} summarizes the steps described above. Here, $s$ is the stepsize, $\alpha$ is the Armijo constant and $k_d$ is the backtracking factor.

\begin{algorithm}[h!]
\caption{Optimization algorithm}
\label{alg:optwavesaxpark}
\begin{algorithmic}[1]
\REQUIRE $s$, tol, maxit, $k_d$, $\alpha$
\STATE Set $k=0$
\WHILE {$k<\text{maxit}$ {\bf and} $\text{err}>\text{tol}$} 
	\STATE Compute $\widehat{\bm{\gamma}}$ and $\widehat{\bm{\zeta}}$ by solving the state problem \eqref{eq:stateprob} with device positions $\bm{x}^k$
	\STATE Compute the cost $J^k$ using \eqref{eq:multibodycost}
	\STATE Compute $\bm{\lambda}$ and $\bm{\mu}$ by solving the adjoint problem \eqref{eq:adjprob}
	\STATE Compute the gradient $\nabla J^k$ using \eqref{eq:gradcostpark}
	\STATE Set $D=I$
	\STATE Perform the projected gradient step $\widetilde{\bm{x}} = \mathcal{P}(\bm{x}^k - s D \nabla J^k)$ \label{line:startstep}
	\WHILE {the set of overlapping pairs of points is not empty}
		\STATE Update $d_{\ell\ell} \leftarrow k_d d_{\ell \ell}$ for each point $(\widetilde{x}_\ell, \widetilde{y}_\ell)$ in an overlapping pair
		\STATE Set $\widetilde{\bm{x}} = \mathcal{P}(\bm{x}^k - s D \nabla J^k)$
		\STATE Solve the state problem with positions $\widetilde{\bm{x}}$ and compute the associated cost $\widetilde{J}$
	\ENDWHILE \label{line:endstep}	
	\WHILE {$\widetilde{J}-J^k > -\alpha \|\widetilde{\bm{x}} - \bm{x}^k\|^2 /(max_\ell d_{\ell\ell} s)$}
		\STATE Update $D \leftarrow k_d D$
		\STATE Repeat lines \ref{line:startstep}-\ref{line:endstep}
	\ENDWHILE 
	\STATE Set $\bm{x}^{k+1} = \widetilde{\bm{x}}$
	\STATE Compute $\text{err} = |J^k - \widetilde{J}|$ and update $k \leftarrow k+1$
\ENDWHILE 
\end{algorithmic}
\end{algorithm}

\subsubsection{Optimization results}\label{sec:optresults}
\begin{figure}[h!]
\centering
\subfloat{
\includegraphics[width=0.45\textwidth]{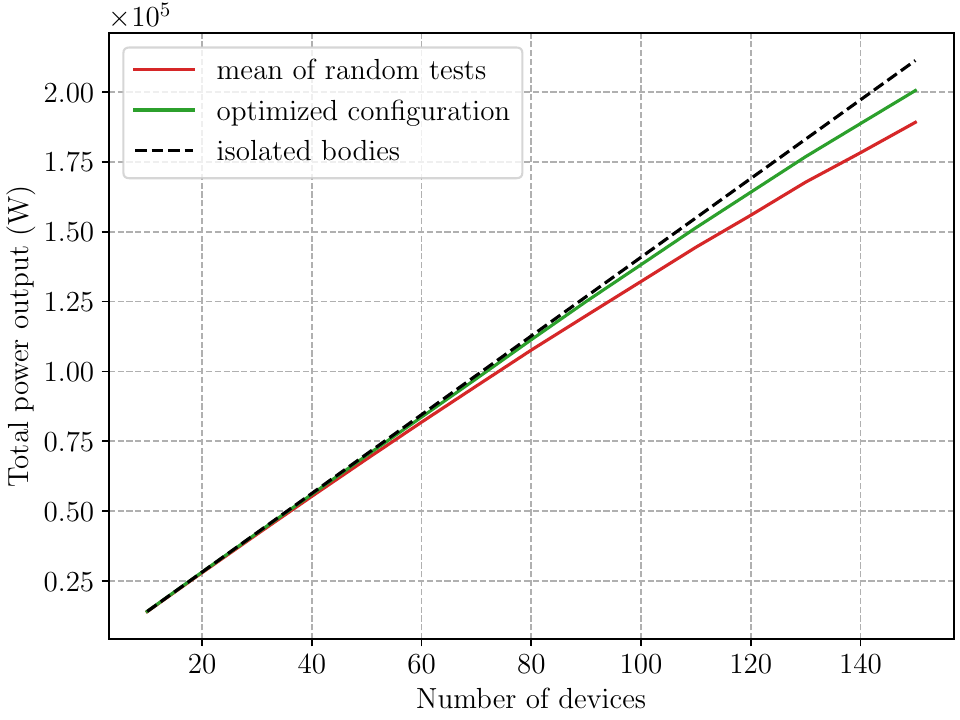}
}
\hfill
\subfloat{
\includegraphics[width=0.45\textwidth]{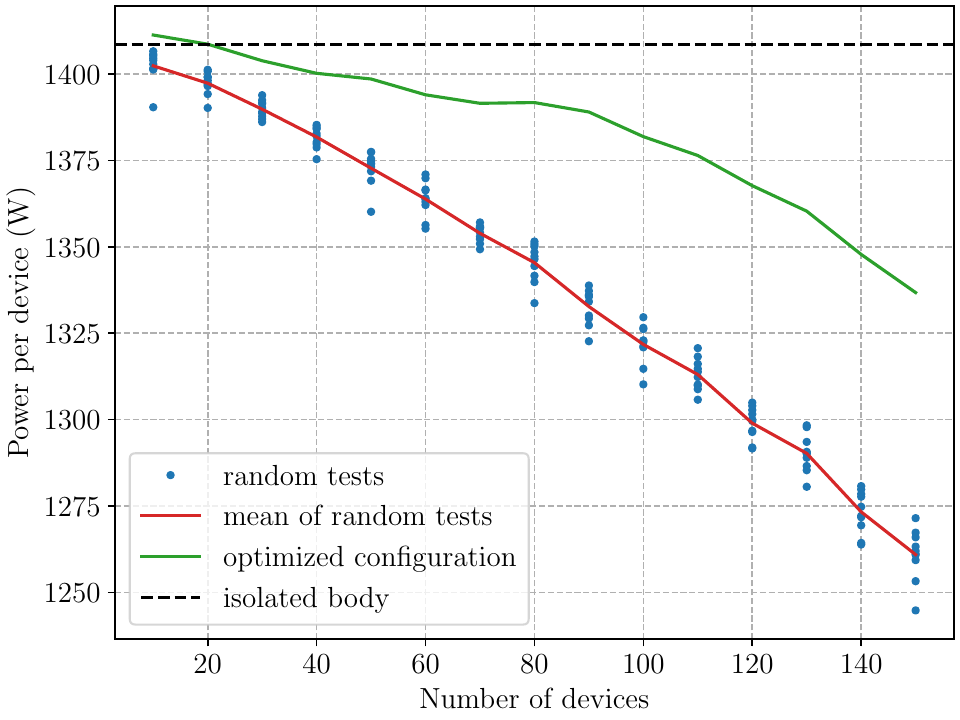}
}
\caption{Results for a triangular domain of edge 50 m, $H_s$ = 3 m, $T$ = 8 s}
\label{fig:poweredge50}
\end{figure}

The optimization strategy presented in Sec.~\ref{sec:optmethod} has been first tested on a triangular domain, without the support structure of the wind turbine, and imposing a minimum distance of 1.6~m between the centers of the devices. For a given number of devices to be installed in the domain, 10 random configurations are obtained by sampling from an uniform distribution with the constraint of minimum distance, and the total power is computed for each of them. Optimization is then carried out starting from the random configuration of maximum power. This process is repeated for increasing values of the number of bodies.
A monochromatic wave equivalent, in the sense explained in \ref{sec:powmatparams}, to the sea state with $H_s = 3\text{ m}$ and $T_e = 8\text{ s}$, and with a wave vector directed as the positive $x$ direction, is considered. The results in terms of total power of the array and average power per device are reported in Fig.~\ref{fig:poweredge50}. The total power is monotonically increasing, while the average power per device has a value greater than the one obtained in the isolated case for an array of 10 devices and less than the isolated value for more than 10 devices (interaction factor less than 1). 
The power gain obtained through optimization is significantly larger than the dispersion of power from random tests. This indicates the efficiency of the proposed optimization strategy\footnote{In other words, finding a configuration as efficient as the optimum randomly would probably require a very large number of trials.}.

\begin{figure}[h!]
\centering
\subfloat{
\includegraphics[height=0.32\textwidth, trim={1cm 0 2cm 0}, clip]{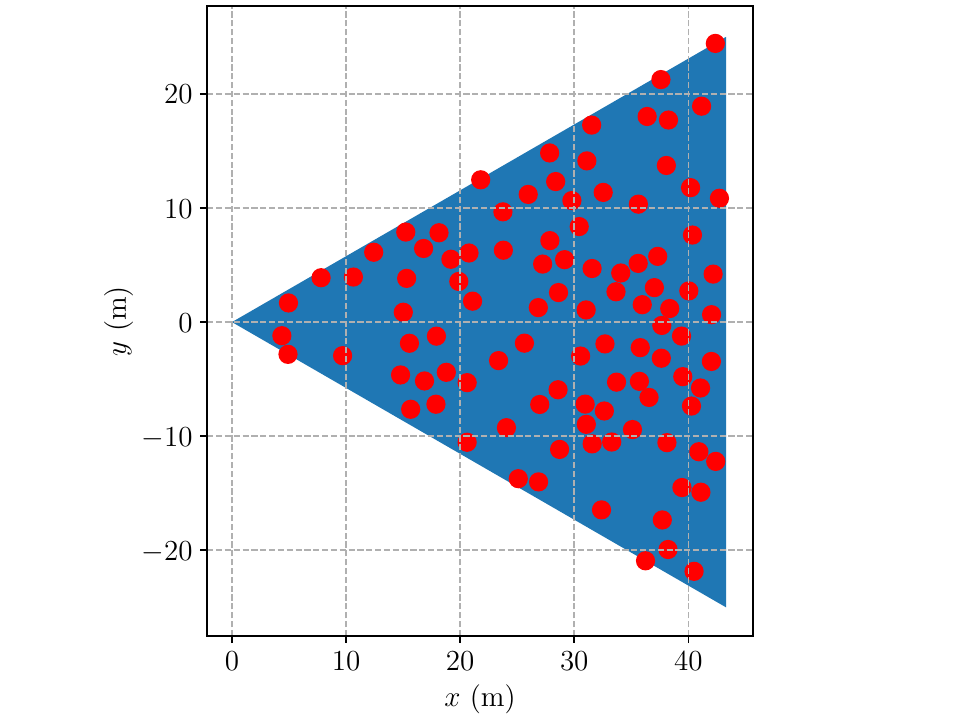}
}
\subfloat{
\includegraphics[height=0.32\textwidth, trim={2.5cm 0 2cm 0}, clip]{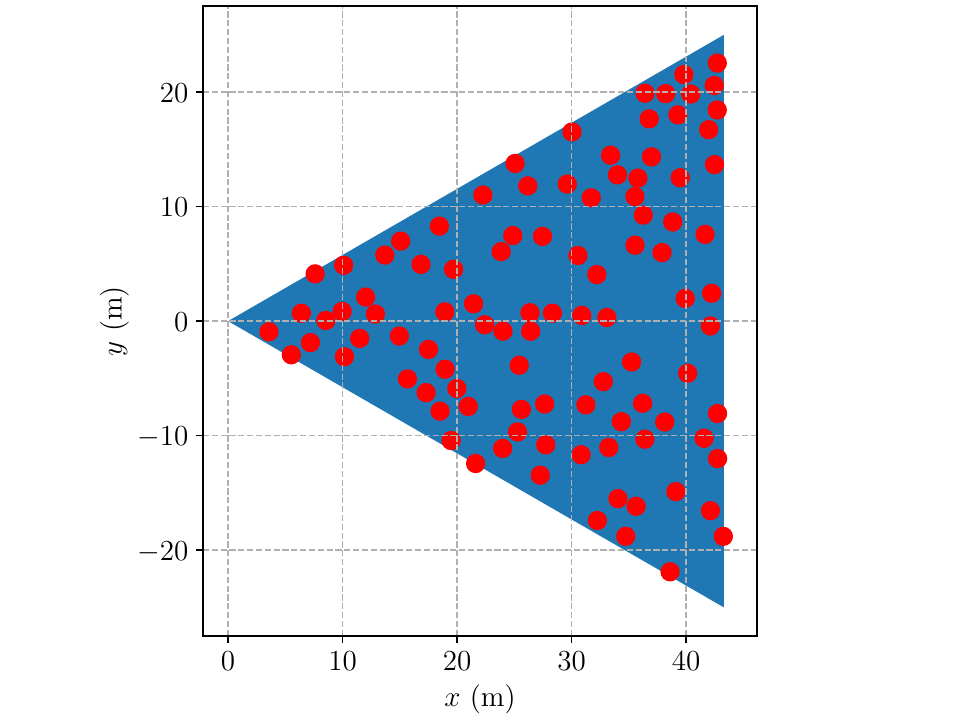}
}
\subfloat{
\includegraphics[height=0.32\textwidth, trim={2.5cm 0 2cm 0}, clip]{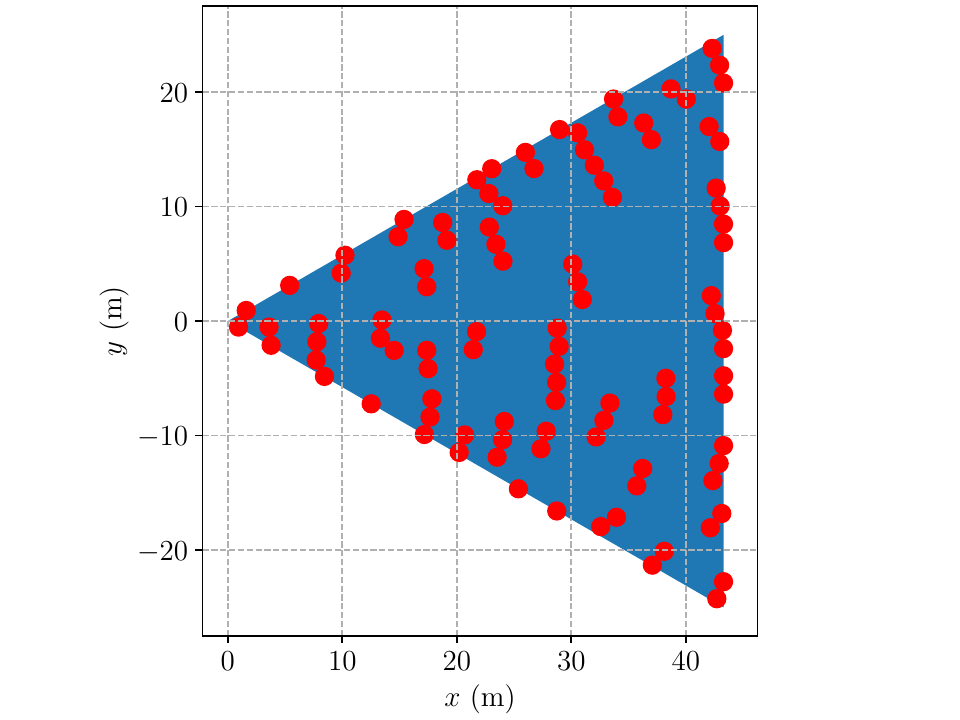}
}
\caption{
Test with 100 devices, $H_s$ = 3 m, $T$ = 8 s. From the left: worst (108.5 kW) and best (109.4 kW) configuration from random tests; optimized configuration (112.4 kW)}
\label{fig:n100edge50}
\end{figure}

\begin{figure}[h!]
\centering
\subfloat{
\includegraphics[width=0.31\textwidth]{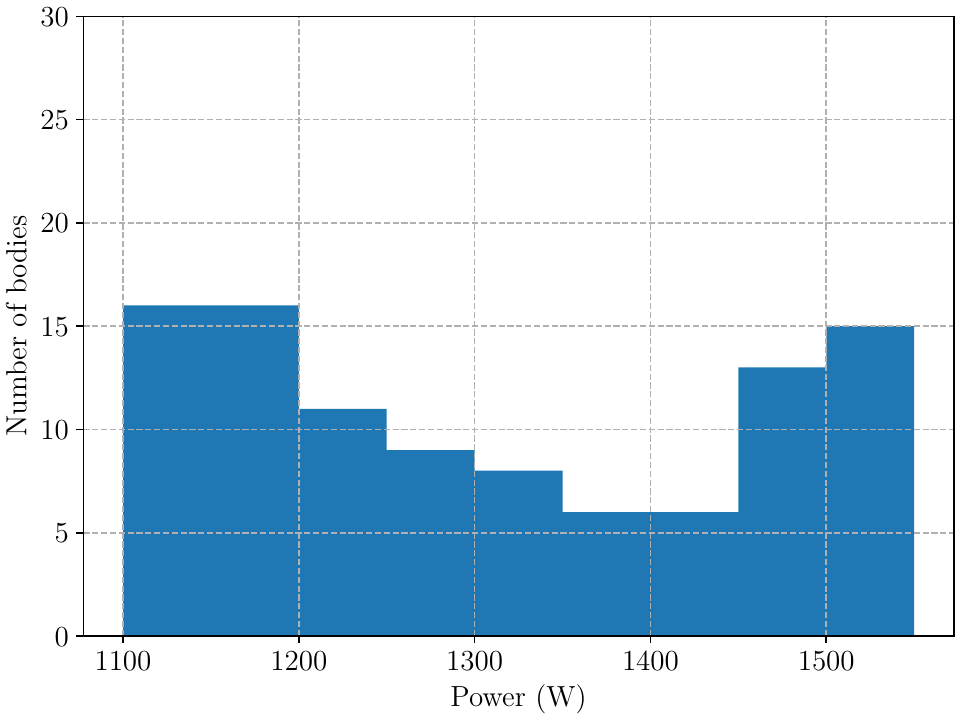}
}
\hfill
\subfloat{
\includegraphics[width=0.31\textwidth]{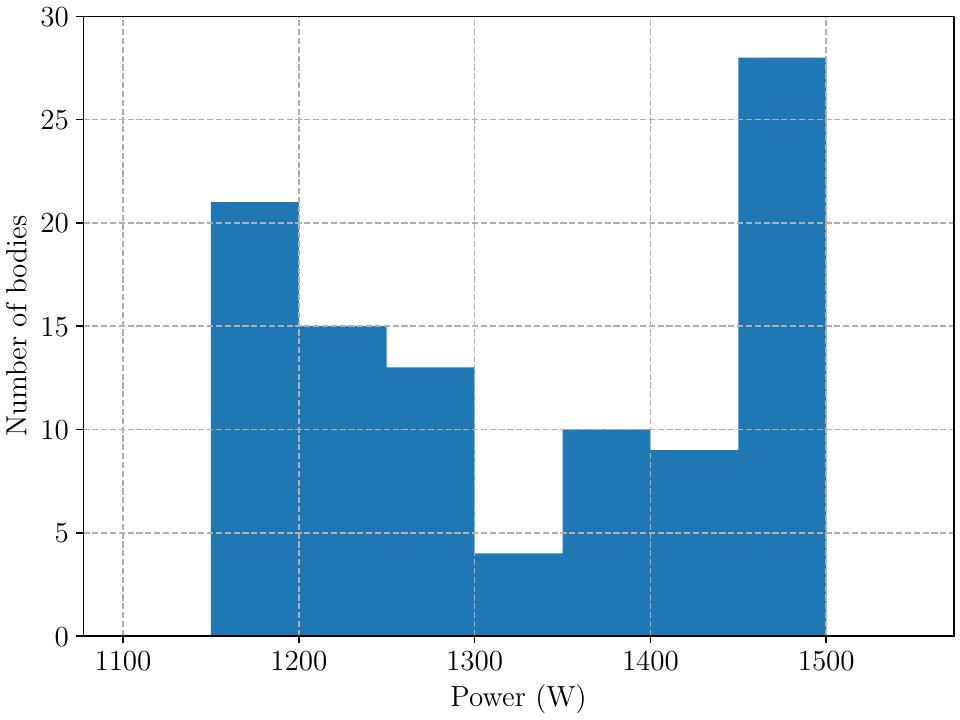}
}
\hfill
\subfloat{
\includegraphics[width=0.31\textwidth]{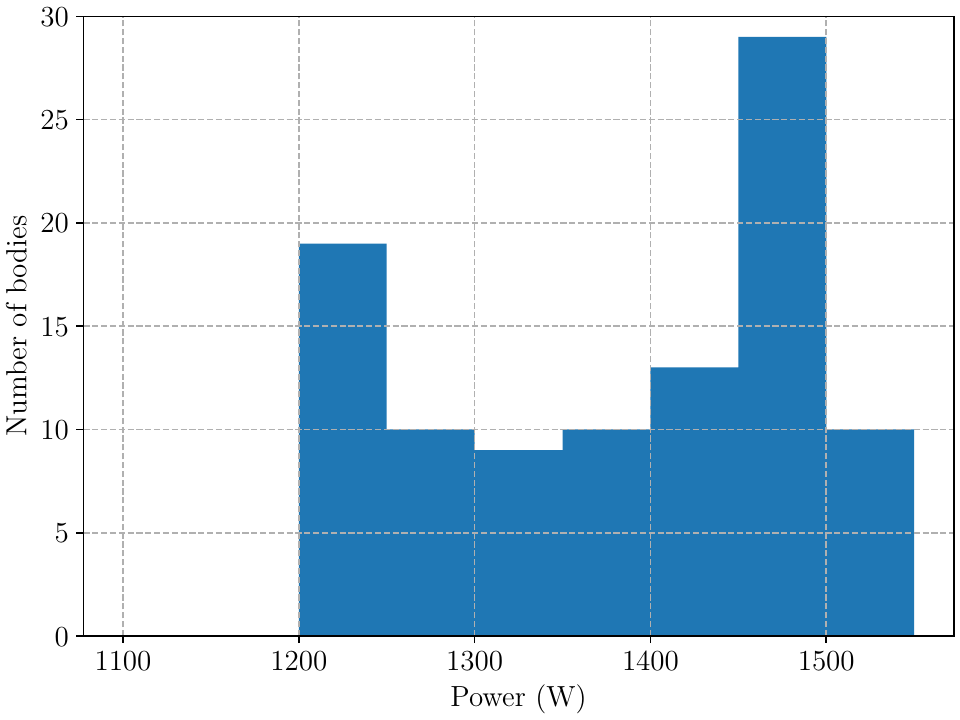}
}
\caption{
Histograms of the powers of single devices, triangular domain, 100 device test. From the left: worst and best configuration from random tests; optimal configuration}
\label{fig:powerhist}
\end{figure}

The best and worst configurations from random tests and the optimized configuration are shown in Fig.~\ref{fig:n100edge50} for the case of 100 devices.
A common observation from numerical tests is that the best configurations in terms of power have a greater density of devices in the upwave region with respect to the worst performing configurations. In the optimal configurations, we observe the formation of lines of closely packed devices.
Histograms of the powers of devices for the same test are reported in Fig.~\ref{fig:powerhist}. The most apparent feature of the histograms is the wide variability in the power values. Since no bodies other than the devices are present, such variability is produced by the mutual  hydrodynamic interactions. This shows that, as expected, the park effect becomes important when the number of devices is large, even if the hydrodynamic interactions in a single pair of devices appear almost negligible (see Sec. \ref{sec:2dev}) and even if, as in this case, the characteristic dimension of the devices is much smaller than the wavelength. 
It can also be observed that the left tail of the distribution in the least favourable condition is absent in the optimized condition. In general, the distribution is shifted to the right by the optimization process.

\begin{figure}[h!]
\subfloat{
\includegraphics[width=0.45\textwidth]{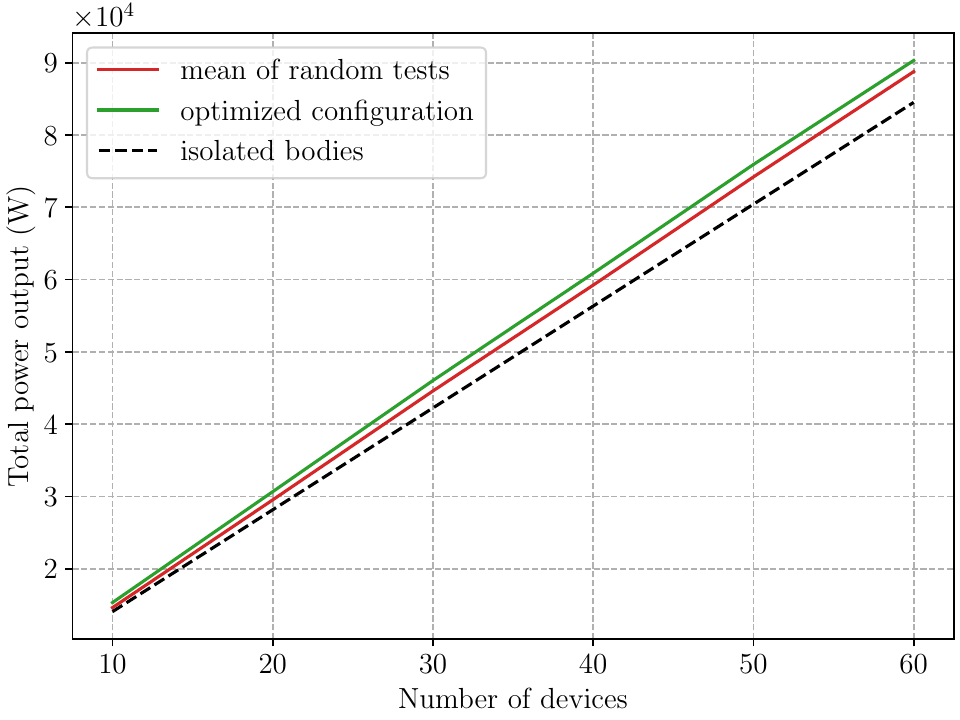}
}
\hfill
\subfloat{
\includegraphics[width=0.45\textwidth]{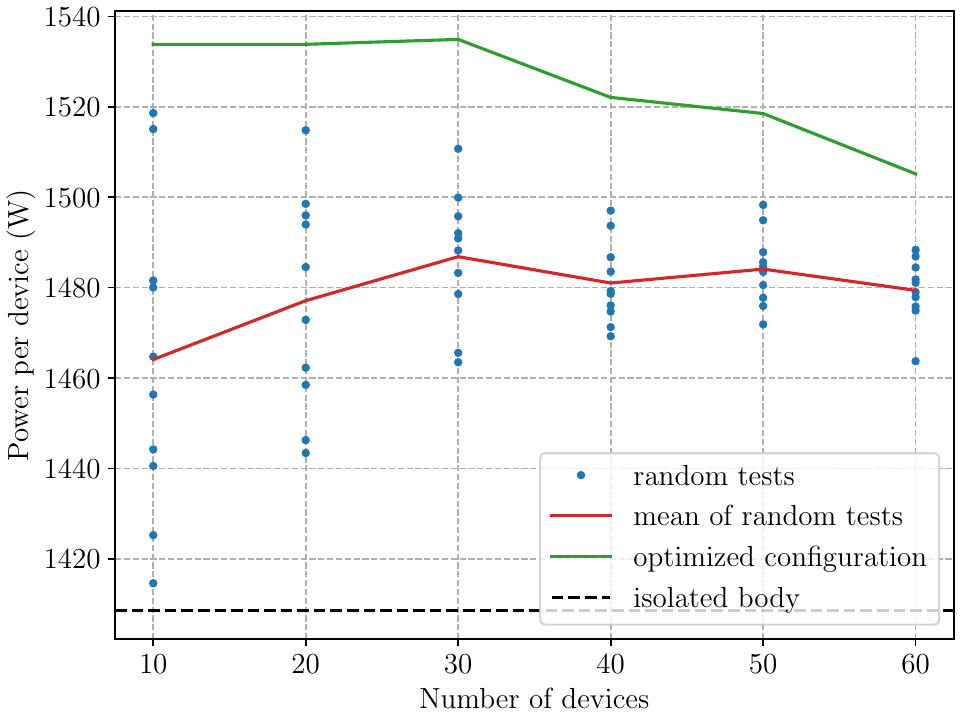}
}	
\caption{Results for the truncated triangular domain}
\label{fig:powertritrunc}
\end{figure}

\begin{figure}[h!]
\centering
\subfloat{
\includegraphics[trim={1cm 0 1cm 0.6cm},clip,width=0.45\textwidth]{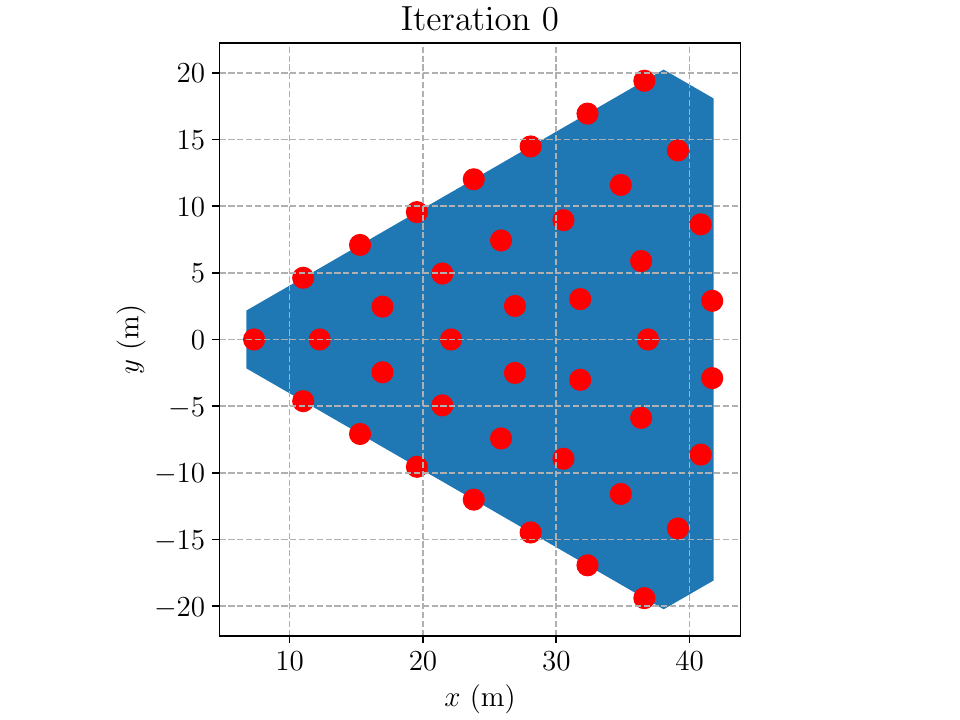}
}
\subfloat{
\includegraphics[trim={1cm 0 0.1cm 0},clip,width=0.46\textwidth]{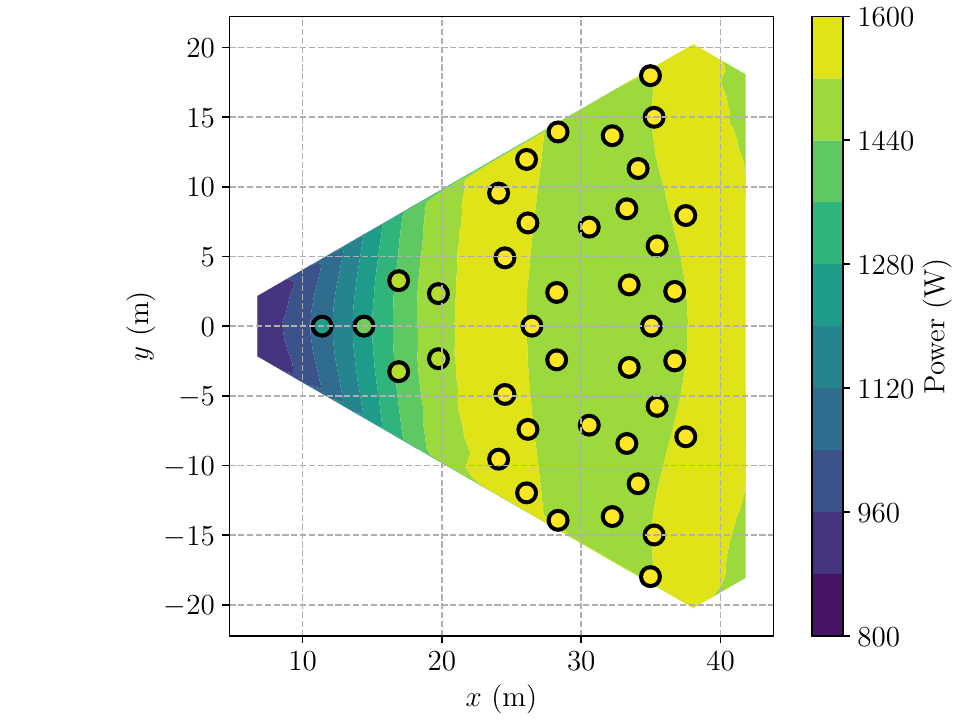}
}
\caption{Initial (left) and optimized (right) configuration; truncated triangular domain, 40 devices. In the right panel, the map of the power of a single device in the domain is reported in the background. The color inside each circle represents the actual power of the corresponding device, taking hydrodynamic interactions into account.}
\label{fig:configtritrunc}
\end{figure}

\begin{figure}[h!]
\centering
\subfloat{
\includegraphics[width=0.45\textwidth]{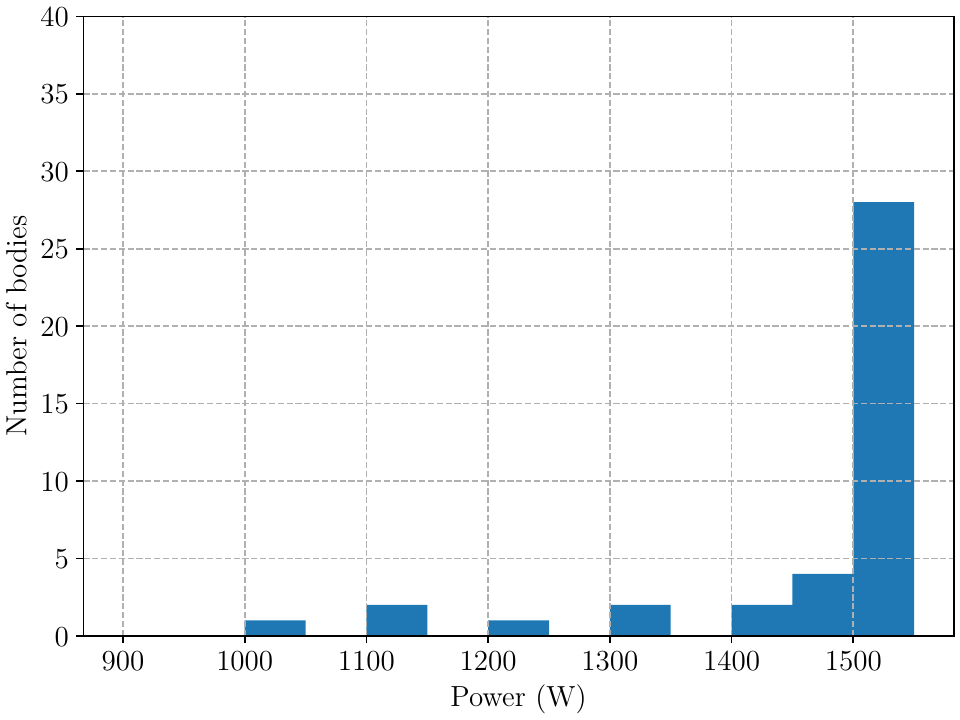}
}
\hfill
\subfloat{
\includegraphics[width=0.45\textwidth]{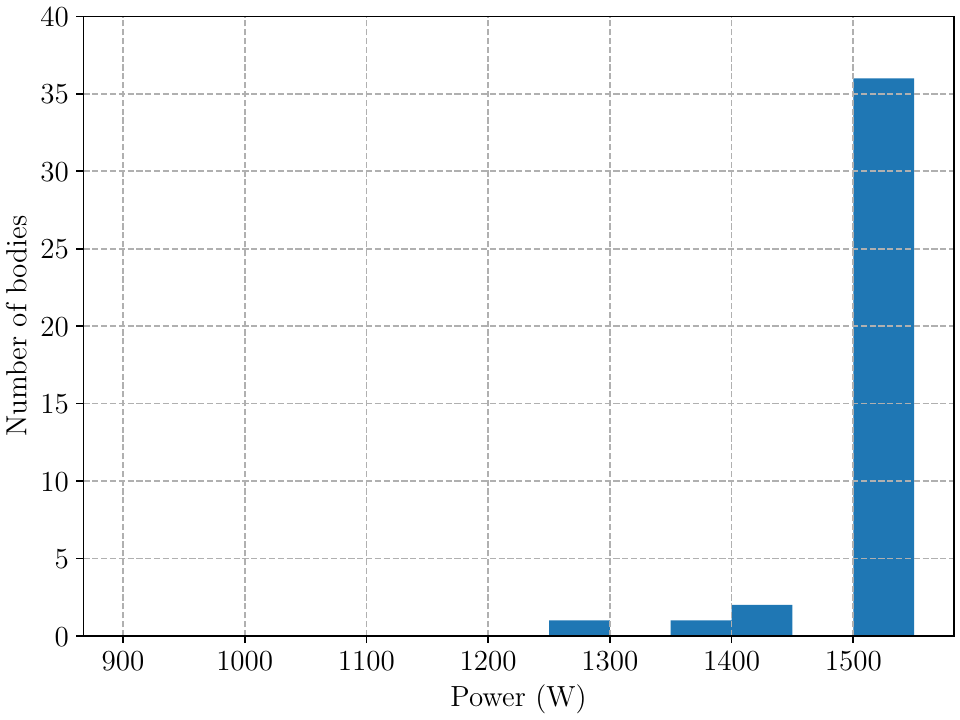}
}
\caption{Histograms of the powers of single devices, triangular domain, 100 device test. Initial (left) and optimized (right) configuration}
\label{fig:histtritrunc}
\end{figure}

In the next tests, the truncated triangular domain has been considered, taking into account the interaction with the piles of the structure on which the wind turbine is mounted. A constraint of minimum distance equal to 3~m between the centers of the devices has been imposed, and the wave vector is directed along the positive $x$ direction. 
The results in terms of power are shown in Fig.~\ref{fig:powertritrunc}. In this case, the mean power of the devices is always greater than the power of an isolated device, thanks to the favourable interaction with the piles. The constraint of minimum distance between the devices significantly reduces the admissible set as the number of devices increases. This, in turn, reduces the dispersion of the results of random tests and the power gain obtained through optimization. We observe that the mean power in the optimized configurations is very close to the maximum power of a single device inside the triangle, that is equal to 1540~W (see Sec.~\ref{sec:1devpiles}).
Optimization runs have been also performed starting from symmetric initial conditions, in order to evaluate the difference with respect to random initial conditions. No significant differences are observed in terms of power. An example of symmetric configuration with 40 devices is shown in Fig.~\ref{fig:configtritrunc}. In the right panel of the figure, corresponding to the optimal configuration, we have superimposed the optimal configuration over the map of the power of a single device in the domain. The color inside each circle represents, instead, the power absorbed by each device considering the hydrodynamic interactions, with the same color scale as the background. It can be observed that the power output of most of the downwave devices is close to the maximum obtainable by a single device. The most upwave converters achieve a production that is greater than the value expected in the single device case, presumably thanks to favourable interactions with the downwave objects.
Power histograms are reported in Fig.~\ref{fig:histtritrunc}. As in the case of the triangular domain without piles, the effect of optimization is to trim the left tail of the power distribution. As a consequence, the number of devices in the bin of greatest power is increased.

\subsubsection{Comparison with the nonlinear model}\label{sec:linnonlin}
The results of optimization, performed using the linear model described in the previous sections, are now verified by recomputing the powers of all explored configurations with a nonlinear model.
We have considered the multibody version of \eqref{eq:nonlindynrad} for a device of constant section $S$:
\begin{equation}
\rho C(\zeta_\ell) S \ddot{\zeta}_\ell + \sum_{m=1}^{N_b} A_{\ell m}(\omega) \ddot{\zeta}_m + \Delta p (v_{t,\ell} \omega_t) + \sum_{m=1}^{N_b} B_{\ell m} (\omega) \dot{\zeta}_m + \rho g \zeta_\ell = p_{e,\ell}(t), \quad \ell=1,\dots,N_b,
\label{eq:multibody-nonlindyn}
\end{equation}
where hydrodynamic interactions between different devices are described by the off-diagonal terms of matrices $\bm{A}(\omega)$ and $\bm{B}(\omega)$; moreover, $p_{e,\ell}$ contains diffraction effects. Hydrodynamic properties are obtained by solving the hydrodynamic equations introduced in Sec.~\ref{sec:extmodel} with the methods described therein. Equation \eqref{eq:multibody-nonlindyn} is solved in the time domain over 8 wave periods, the instantaneous power of the park is computed at each time step and the mean power is obtained by averaging over the last 4 periods.

\begin{figure}[h!]
\centering
\subfloat{
\includegraphics[width=0.45\textwidth]{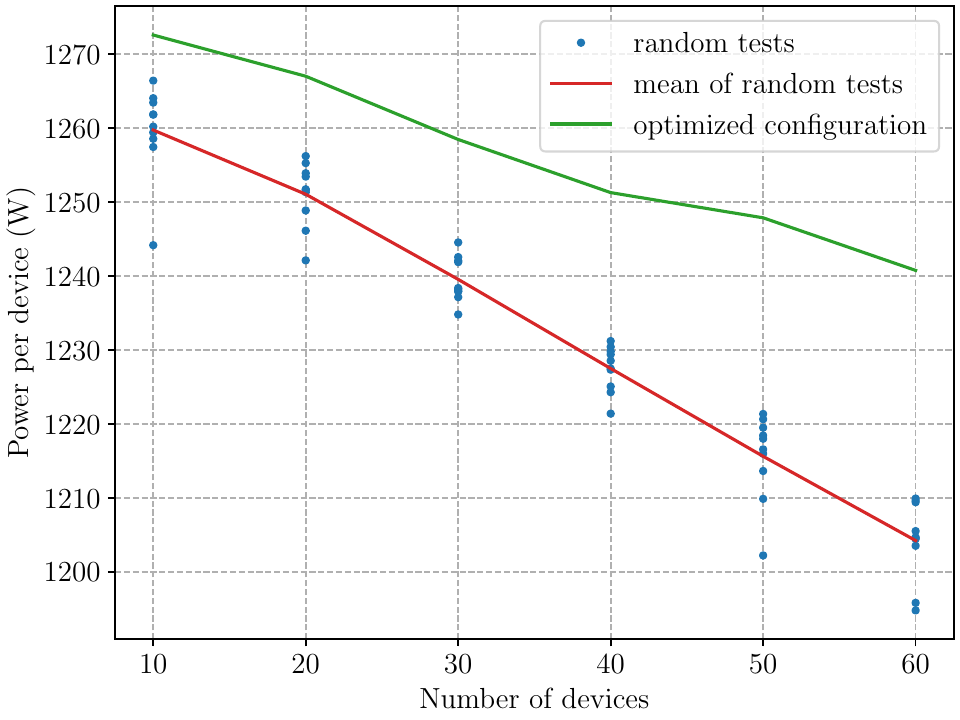}}
\hfill
\subfloat{
\includegraphics[width=0.45\textwidth]{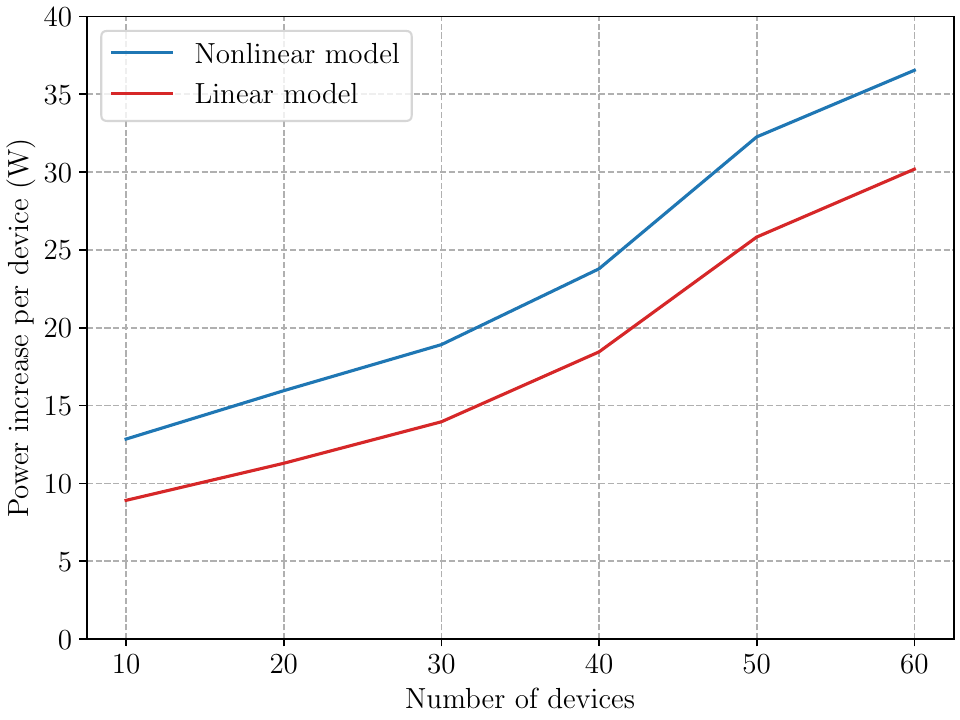}}
\caption{Recomputation of the power per device from the optimization over the triangular domain of edge 50 m using the nonlinear model (left; compare with Fig.~\ref{fig:poweredge50}) and power increase per device obtained through optimization and computed with the linear and nonlinear models (right)}
\label{fig:poweredge50-nonlin}
\end{figure}

The recomputed powers for the case of the triangular domain without piles (see Fig.~\ref{fig:poweredge50}) are reported in Fig.~\ref{fig:poweredge50-nonlin}. Simulations have been carried out up to 60 devices, since larger test cases would require a significant amount of memory for the hydrodynamic computations. While the absolute values of power differ by about 10\%, power in the optimized configuration is confirmed to be significantly greater than powers from random configurations. In particular, as shown in the right panel of Fig.~\ref{fig:poweredge50-nonlin}, the power gained through optimization is comparable between the linear and nonlinear results. Thus, the linear model may be an adequate surrogate of the nonlinear model for purposes of optimization.

\section{Conclusions}
A mathematical model for a water-turbine OWC device has been introduced, incorporating the behavior of the internal water column, of the Wells turbine, and the interaction with the outside environment. The model has enabled the optimization of the turbine's rotational speed for each sea state, and the estimation of the power matrix.

A linearization of the model has been then derived and used to explore the dependence of the annual power production on the dimensions of the device. Then, a generalization to arrays, achieved by limiting the treatment to cylindrical devices, has been introduced. Hydrodynamic simulations show that interactions with the piles of an offshore structure can be beneficial in terms of wave energy production, thanks to the reflections from the piles.  Interestingly, even though the diffraction properties of a single device are negligible with respect to mass and stiffness, a large power variability between the devices is observed, showing a significant mutual interaction effect.

Finally, an optimization method has been implemented starting from the linear WEC park model. The results show that optimization improves the power per device significantly by shifting the distribution of powers upwards. Moreover, an a-posteriori verification using the nonlinear model corroborates the choice of the linear model as a cheap surrogate, justifying its use in optimization. 

Finally, we remark that the optimization strategy presented here can be used for more general classes of WEC devices. Furthermore, it can be adapted to perform control co-design, i.e., simultaneous optimization of layout and control variables. This is the subject of a work in preparation \cite{noidopo}.

\section*{Acknowledgements}
M.G., E.M., and G.C. are members of the GNCS Indam group. The present research is part of the activities of ``Dipartimento di Eccellenza 2023-2027." This work has been financed for S.M. and G.A. by the Research Fund for the Italian Electrical System under the Three-Year Research Plan 2022-2024 (D.M. MITE n. 337, 15.09.2022), in compliance with the Decree of April 16th, 2018.

\printcredits

\clearpage
\appendix
\section{Power matrices and related results} \label{sec:pmat}
This section contains the main results of the simulation and control optimization of a single device, both in the original configuration shown in Fig.~\ref{fig:scheme} and in the constant section configuration. For the latter, a comparison with the linearized model is included. We refer to Sec.~\ref{sec:contrpowmat} for the discussion of the results.
\begin{figure}[hp!]
\subfloat{
\includegraphics[width=0.45\textwidth]{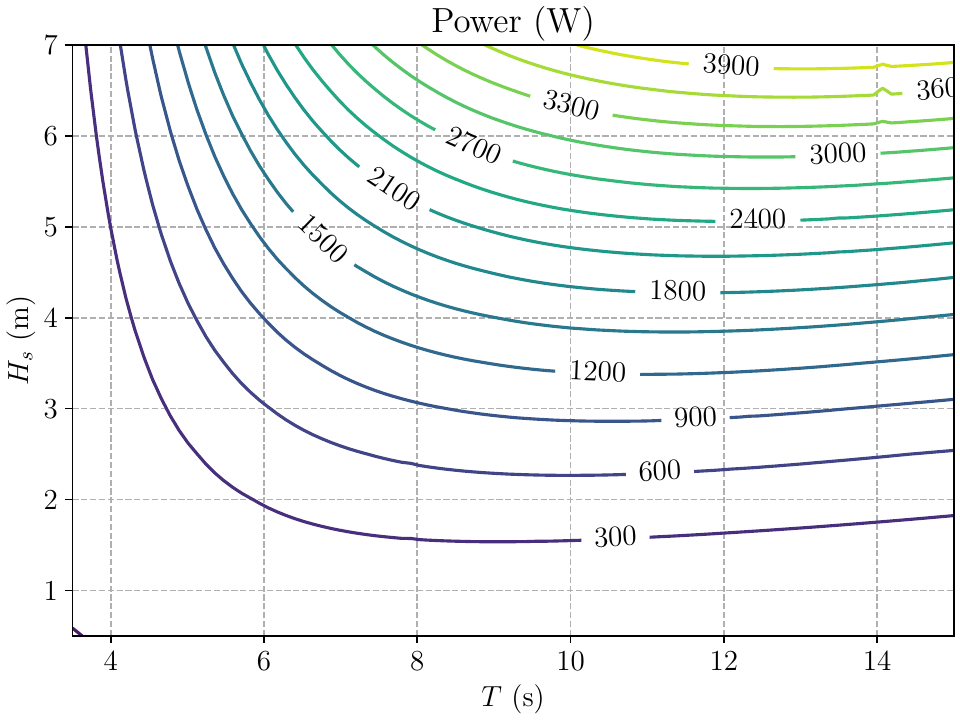}
}
\hfill
\subfloat{
\includegraphics[width=0.45\textwidth]{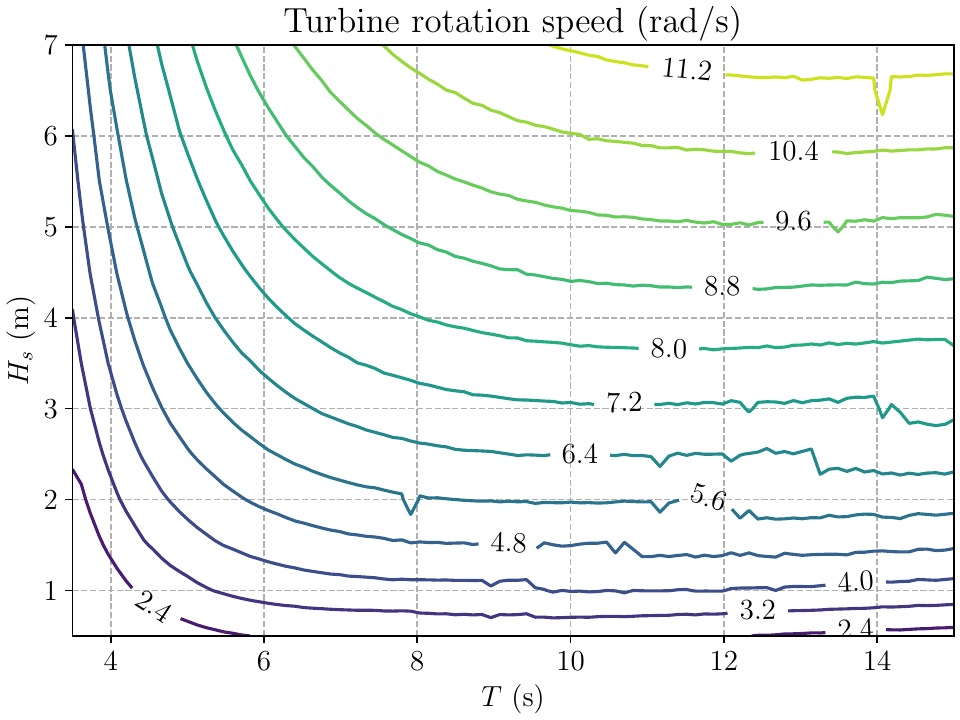}
}
\\
\subfloat{
\includegraphics[width=0.45\textwidth]{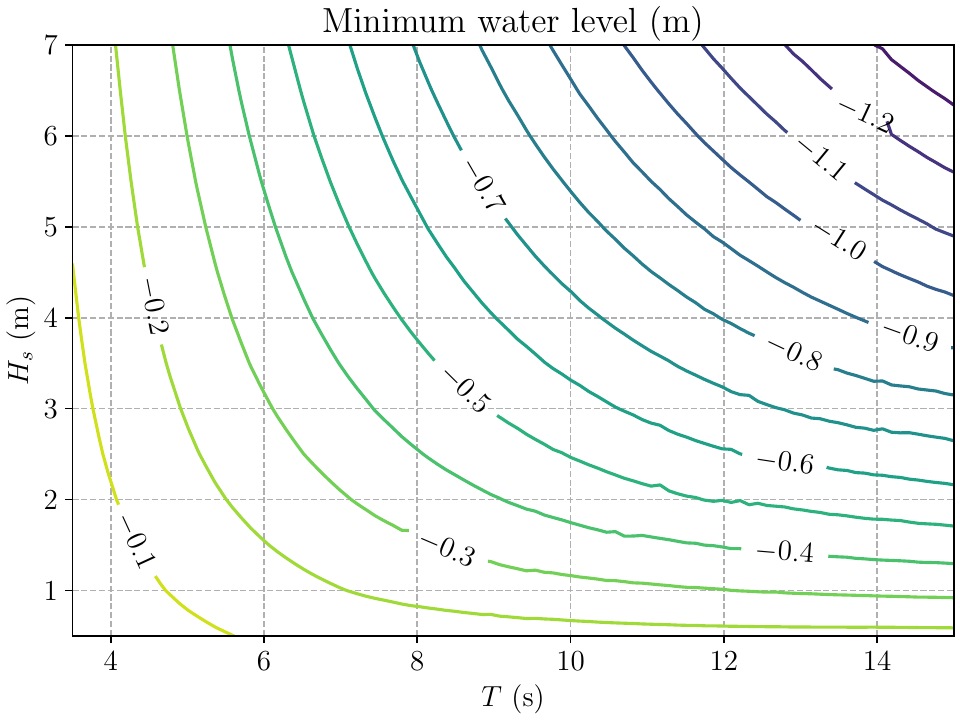}
}
\hfill
\subfloat{
\includegraphics[width=0.45\textwidth]{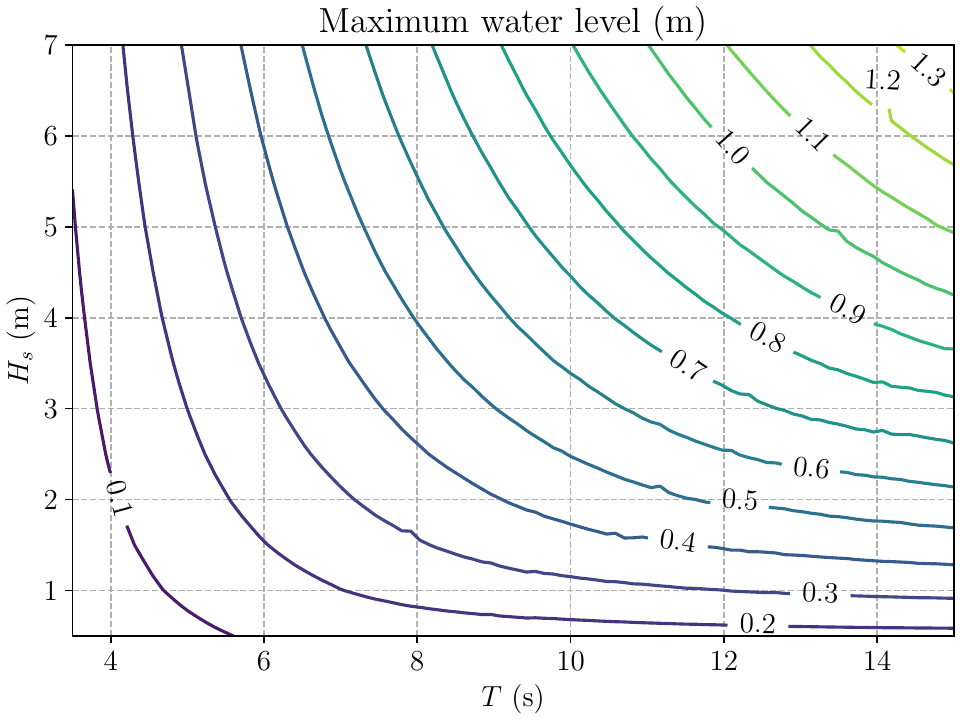}
}
\\
\subfloat{
\includegraphics[width=0.45\textwidth]{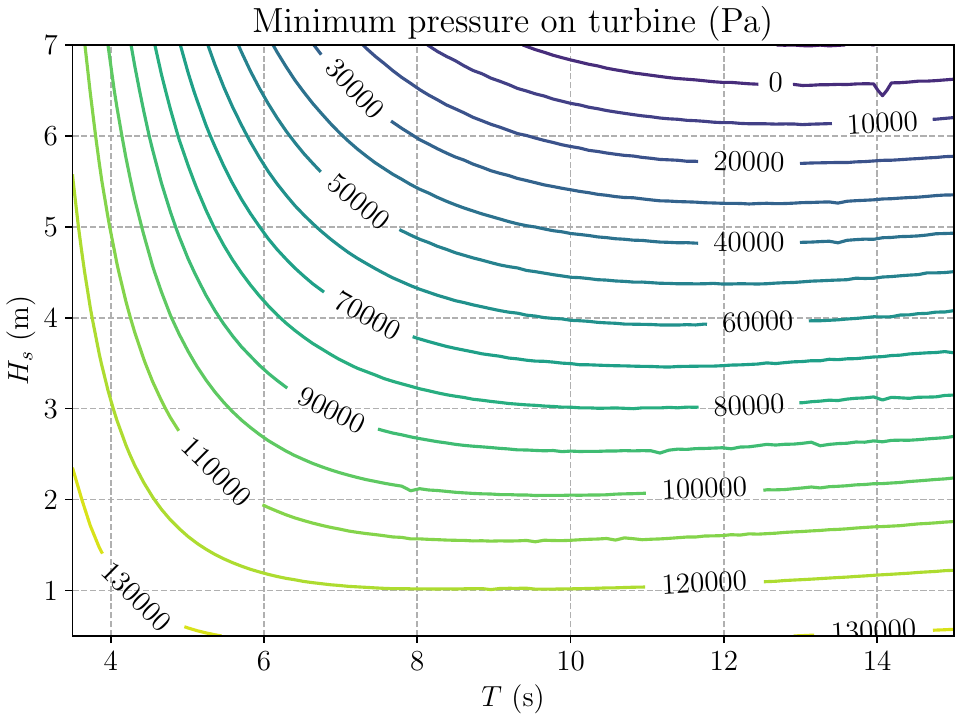}
}
\hfill
\subfloat{
\includegraphics[width=0.45\textwidth]{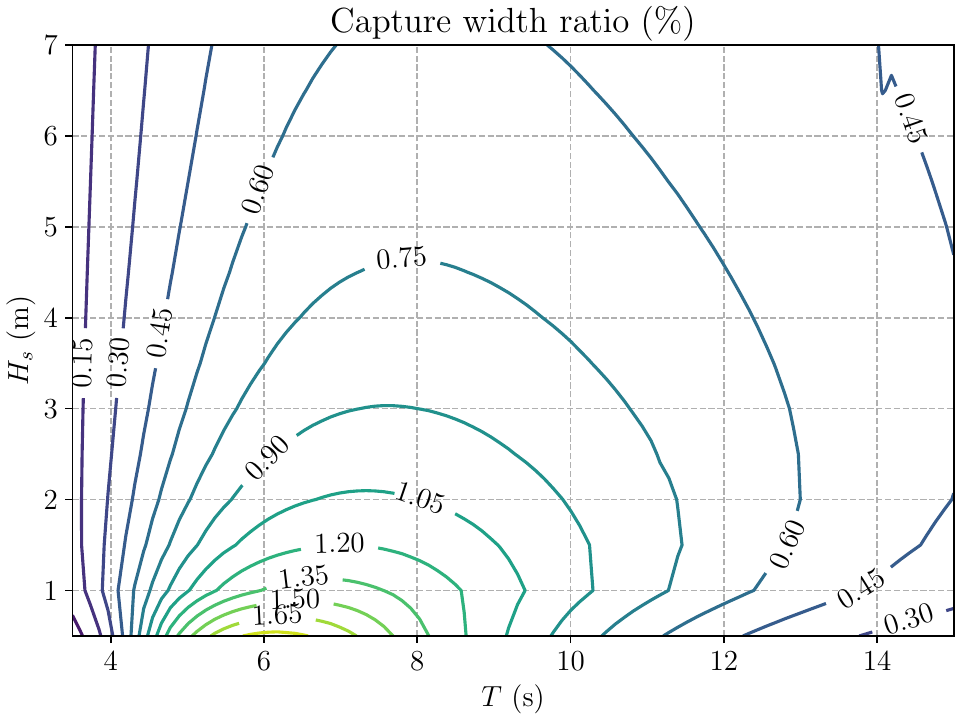}
}
\caption{Isolated device results, original device, 5-blade turbine}
\label{fig:pmatrix5pale-corr}
\end{figure}

%

\begin{figure}[p!]
\subfloat{
\includegraphics[width=0.45\textwidth]{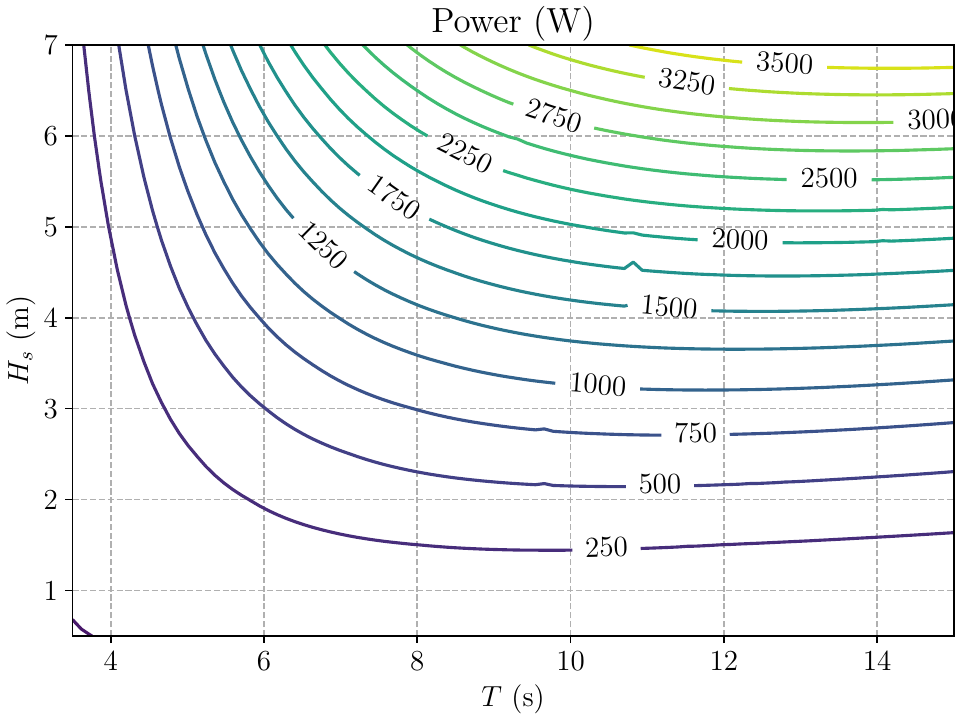}
}
\hfill
\subfloat{
\includegraphics[width=0.45\textwidth]{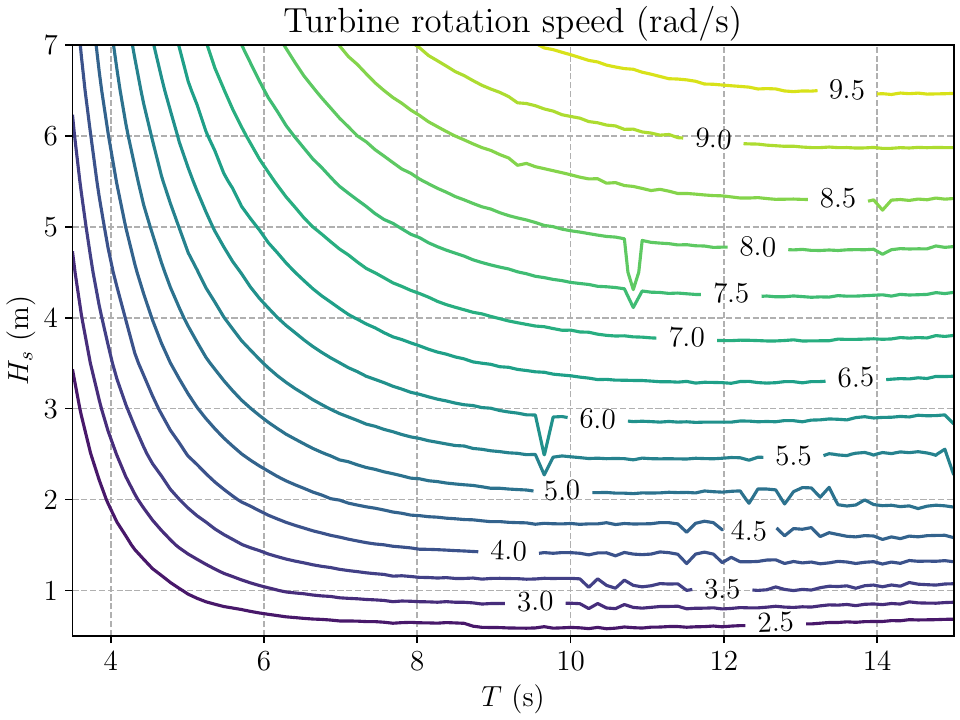}
}
\\
\subfloat{
\includegraphics[width=0.45\textwidth]{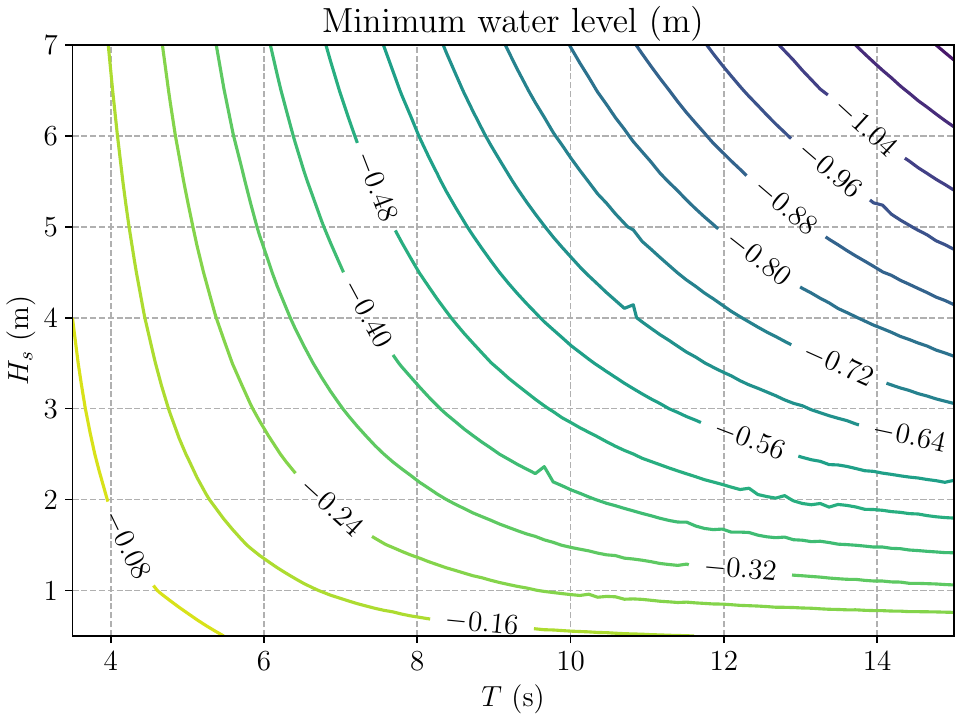}
}
\hfill
\subfloat{
\includegraphics[width=0.45\textwidth]{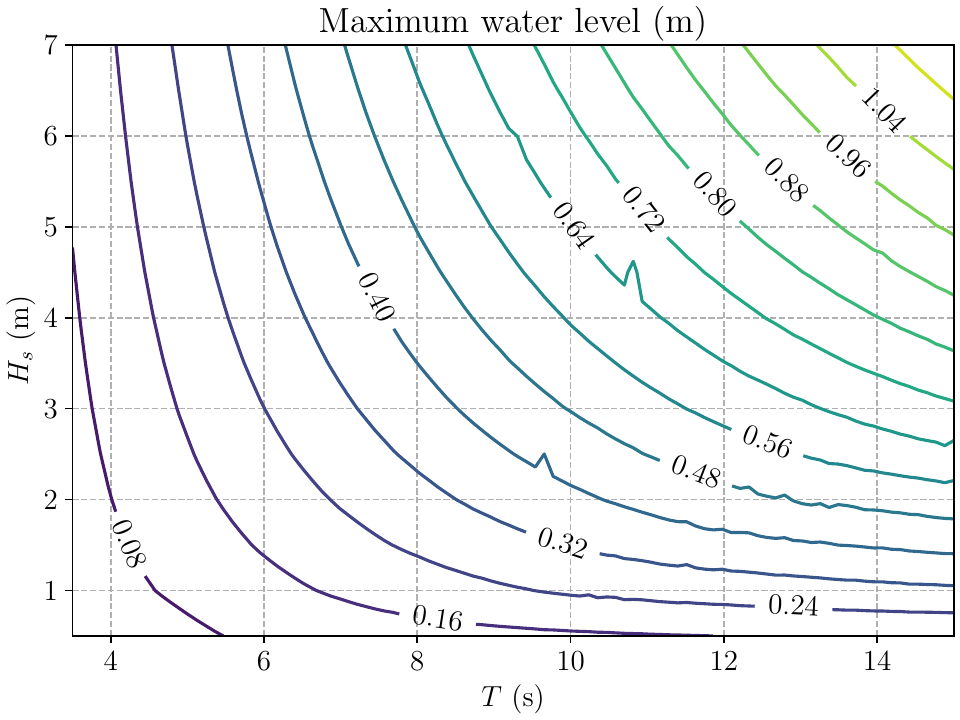}
}
\\
\subfloat{
\includegraphics[width=0.45\textwidth]{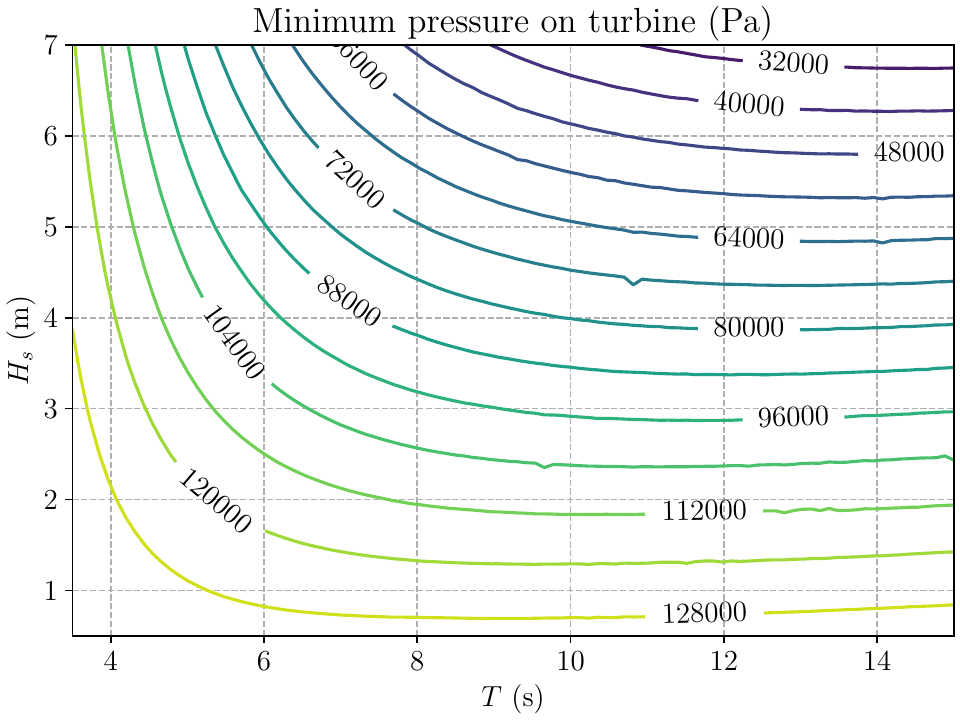}
}
\hfill
\subfloat{
\includegraphics[width=0.45\textwidth]{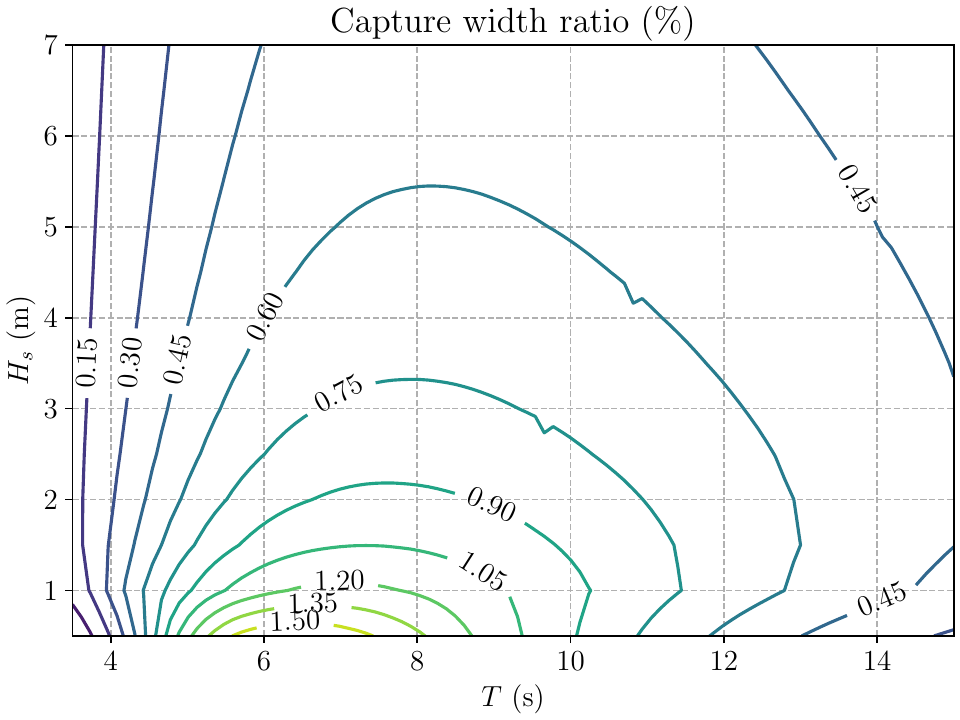}
}
\caption{Isolated device results, original device, 7-blade turbine}
\label{fig:pmatrix7pale-corr}
\end{figure}

\begin{figure}[p!]
\subfloat{
\includegraphics[width=0.45\textwidth]{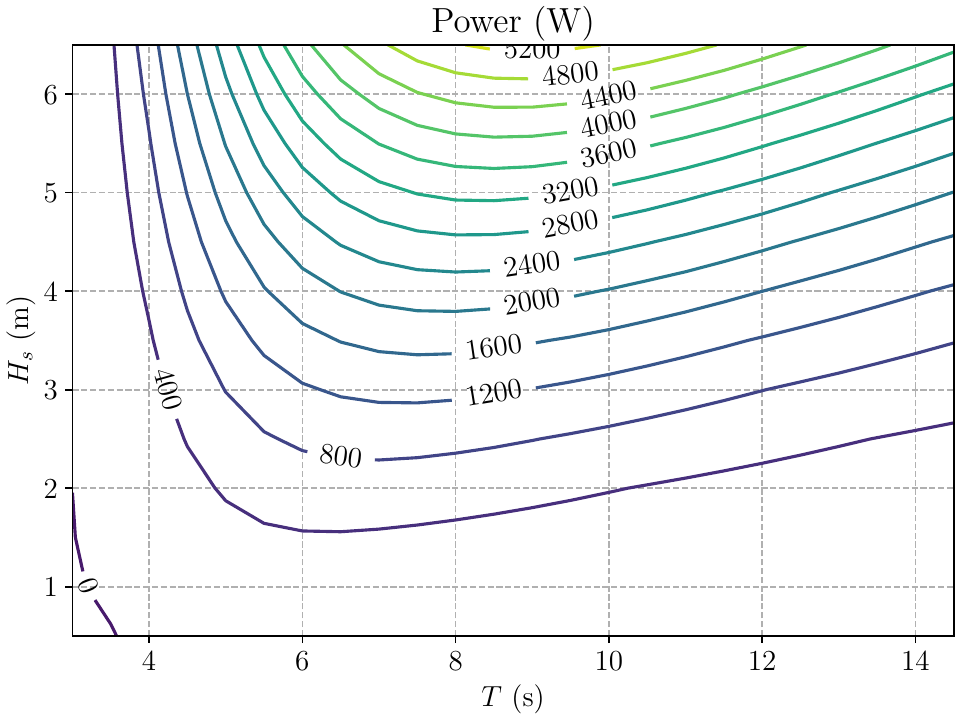}
}
\hfill
\subfloat{
\includegraphics[width=0.45\textwidth]{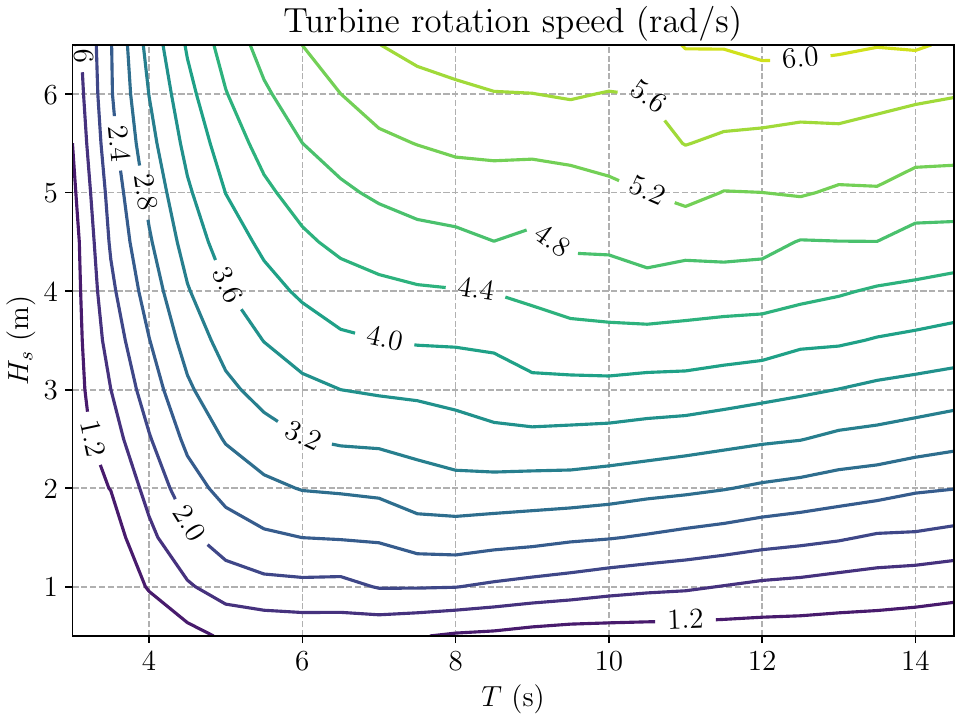}
}
\\
\subfloat{
\includegraphics[width=0.45\textwidth]{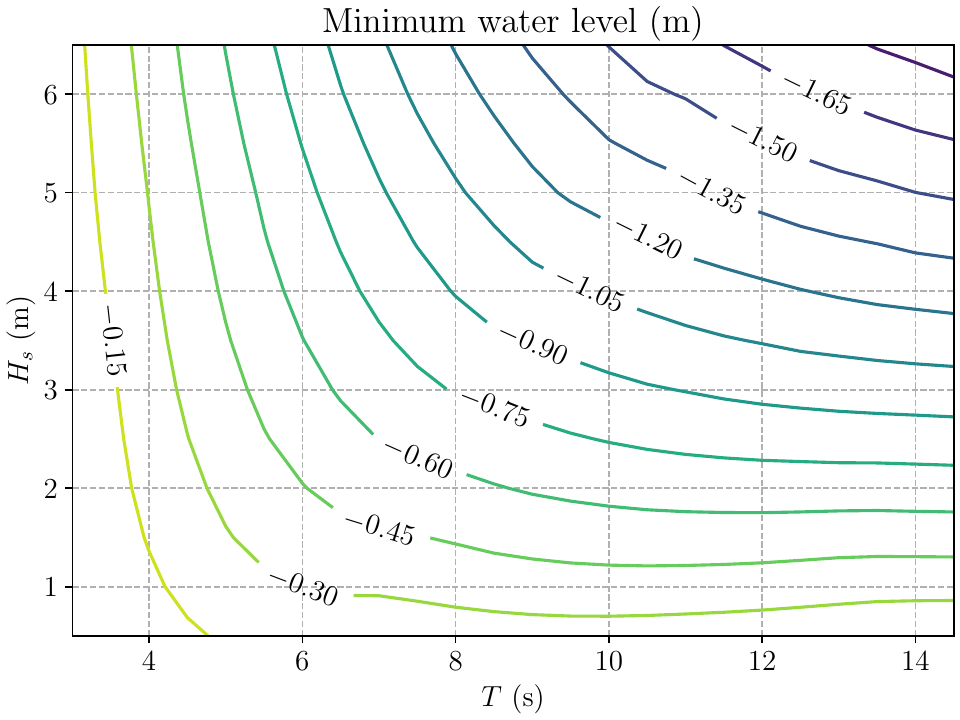}
}
\hfill
\subfloat{
\includegraphics[width=0.45\textwidth]{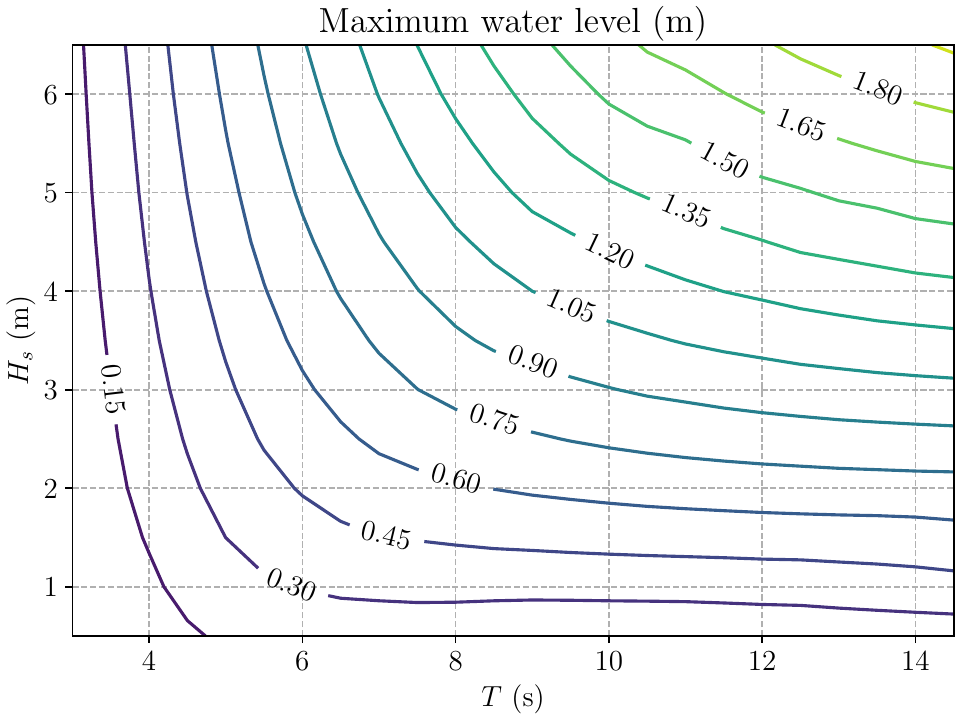}
}
\\
\subfloat{
\includegraphics[width=0.45\textwidth]{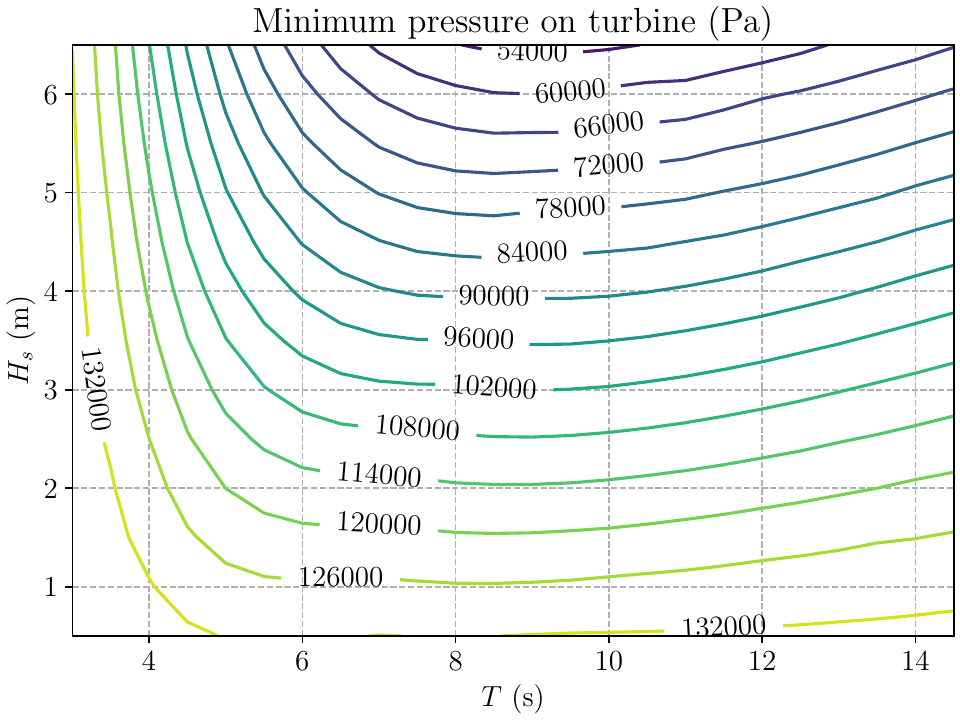}
}
\hfill
\subfloat{
\includegraphics[width=0.45\textwidth]{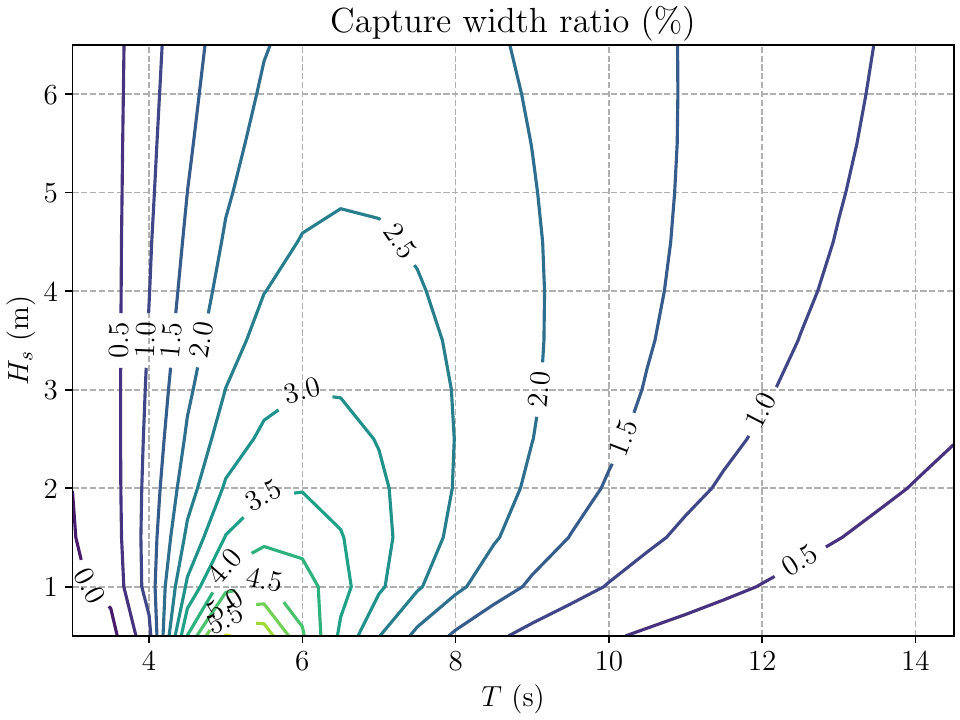}
}
\caption{Isolated device results, constant duct section, 7-blade turbine}
\label{fig:pmatrixsezcost}
\end{figure}

\begin{figure}[p!]
\subfloat{
\includegraphics[width=0.45\textwidth]{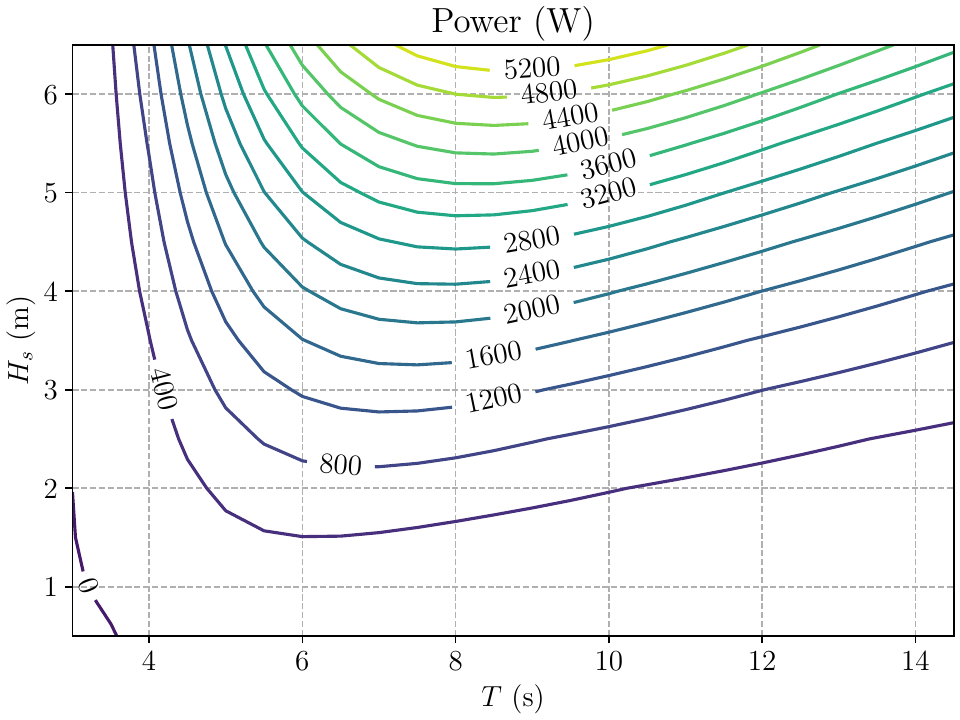}
}
\hfill
\subfloat{
\includegraphics[width=0.45\textwidth]{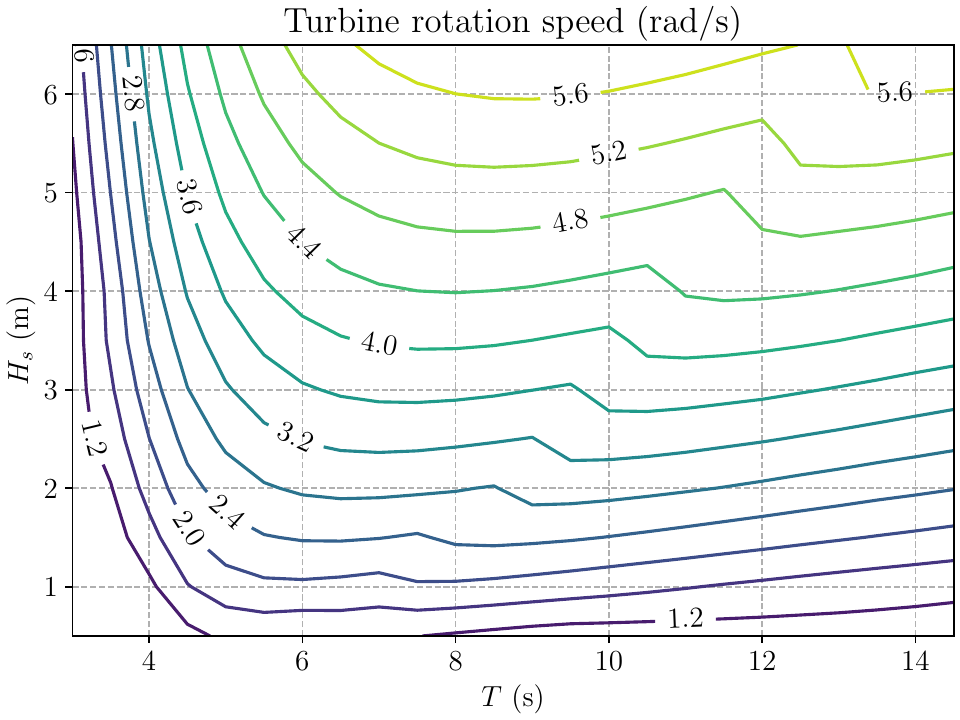}
}
\\
\subfloat{
\includegraphics[width=0.45\textwidth]{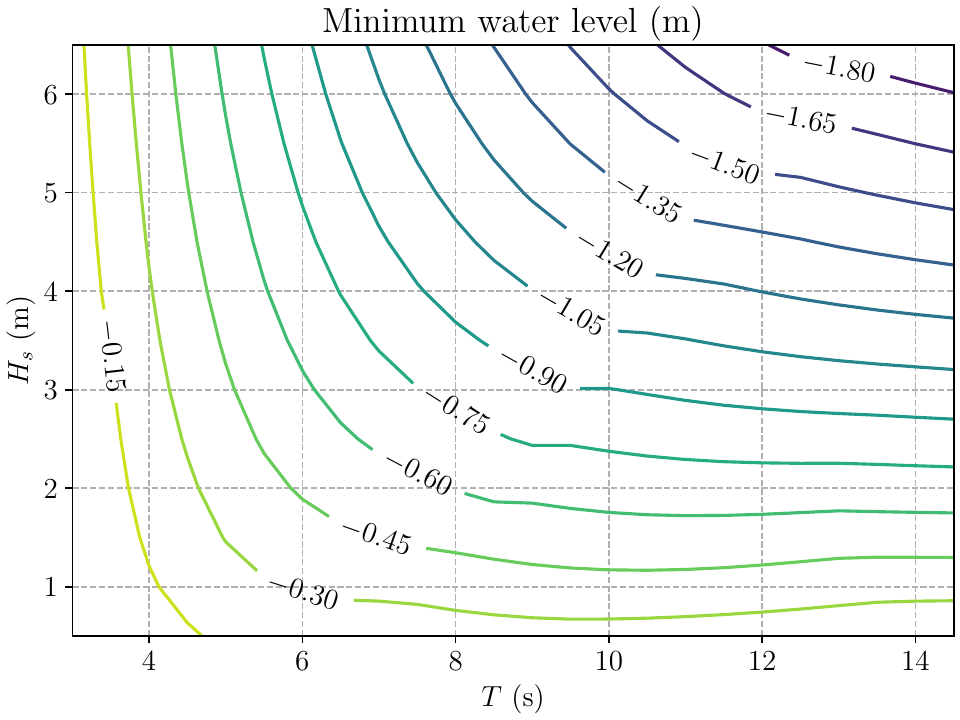}
}
\hfill
\subfloat{
\includegraphics[width=0.45\textwidth]{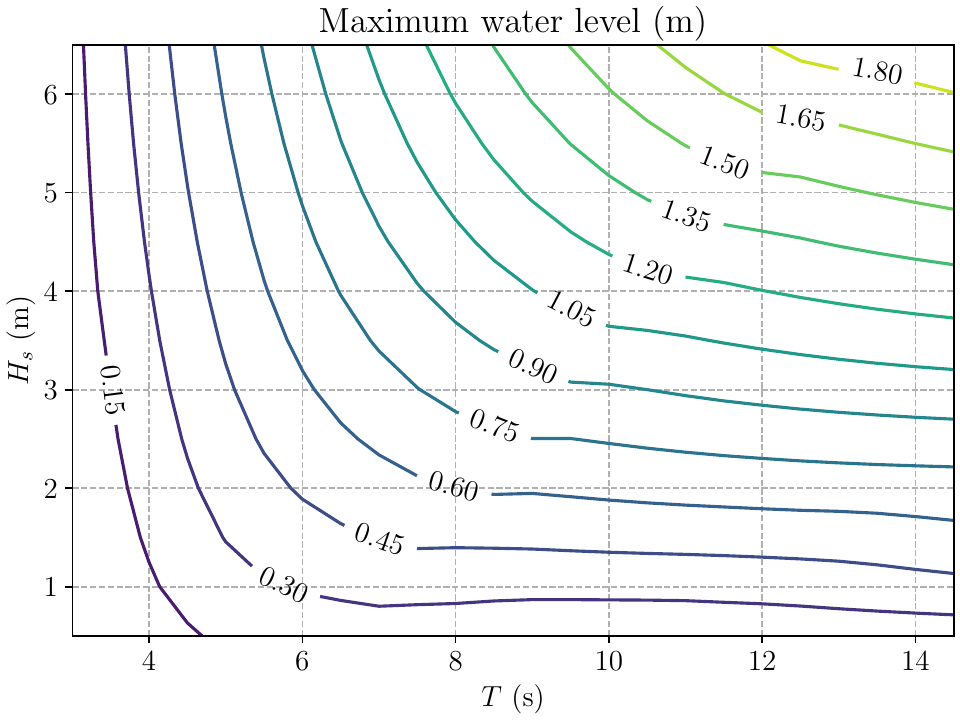}
}
\\
\subfloat{
\includegraphics[width=0.45\textwidth]{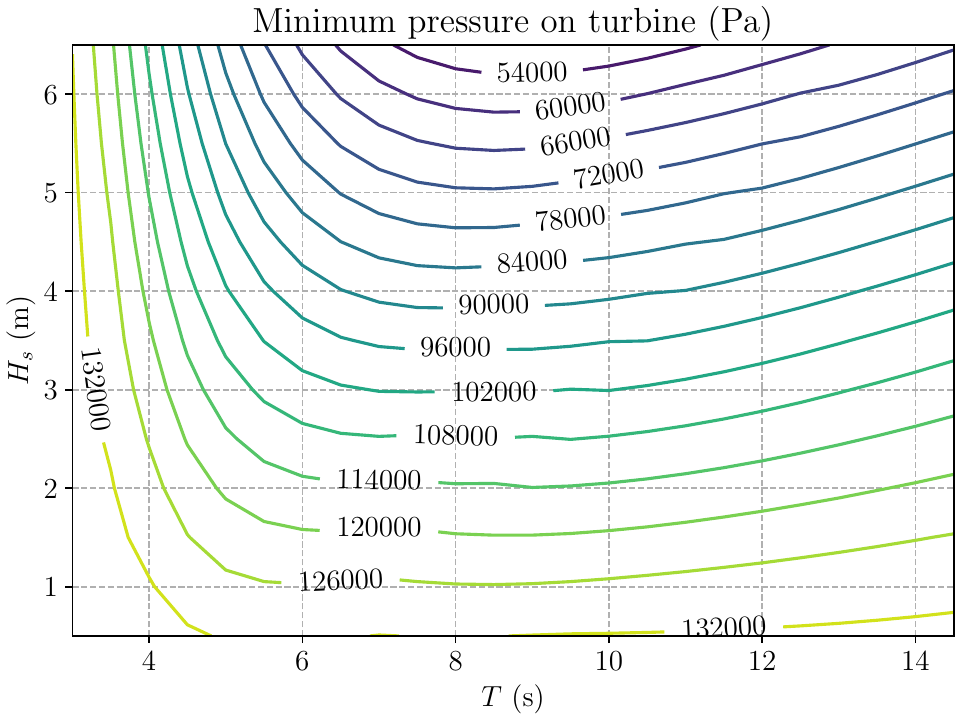}
}
\hfill
\subfloat{
\includegraphics[width=0.45\textwidth]{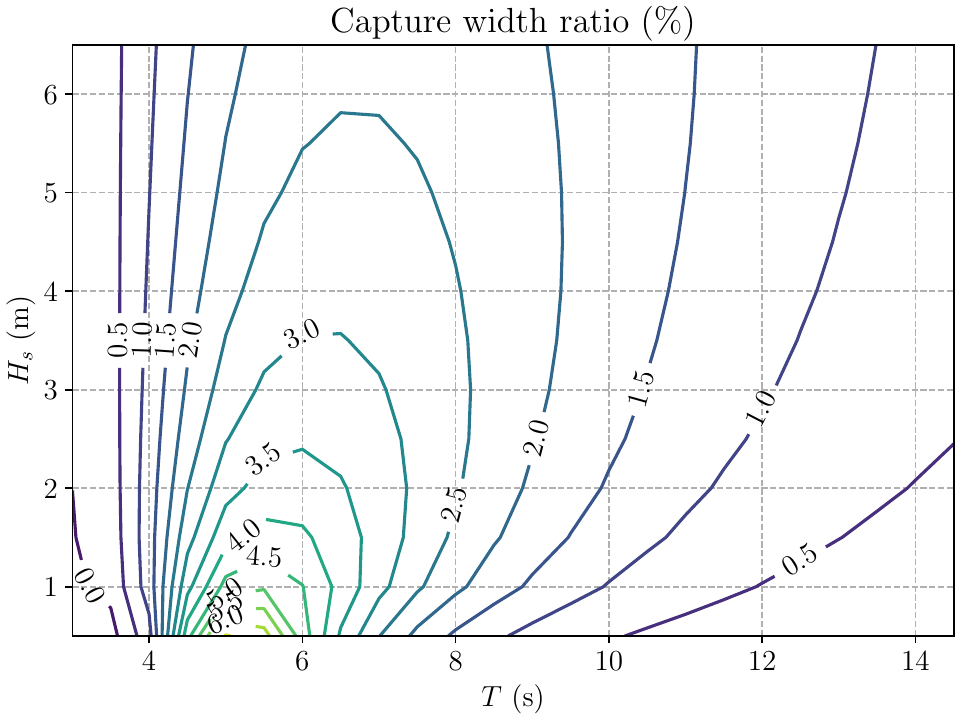}
}
\caption{Isolated device results, constant duct section, 7-blade turbine; linear model}
\label{fig:pmatrixsezcostlinear}
\end{figure}

\clearpage
\section{Power matrix parameters}
\label{sec:powmatparams}
For the purposes of simulation, a wave climate can be approximated as a finite superposition of monochromatic waves \cite{Chakrabarti1987}. The computation of a power matrix in a reasonable time requires considering a rather large number of wave conditions. To obtain first approximation results cheaply, it is possible to represent each wave condition as a single monochromatic wave with the same energy flux.

For a monochromatic wave, the power flux per unit width is $J  =v_g E$, where $v_g$ is the group velocity and $E$ the energy density per unit surface \cite{Falnes2020}. The group velocity is 
\begin{equation}
v_g = \frac{\partial \omega}{\partial k},
\end{equation}
from which, considering the deep water dispersion relation $\omega^2=kg$, we have $v_g = g/(2\omega)$. The energy density is in general $E = \rho g \overline{\eta^2}$. A monochromatic wave can be expressed as $\eta = A\cos(kx - \omega t)$, which yields $\overline{\eta^2} = A^2/2.$ Since the wave height $H$ is conventionally defined as twice the wave amplitude $A$, the power flux of a monochromatic waves of height $H$ and period $T$ is
\begin{equation}
J = \frac{\rho g^2}{32 \pi} T H^2.
\label{eq:powfluxmonochr}
\end{equation}
For irregular waves, energy is expressed in terms of spectral density $S(f)$:
\begin{equation}
E = \rho g \int_0^\infty S(f) {df},
\end{equation}
where in particular 
\begin{equation}
\int_0^\infty S(f) {df}  =\overline{\eta^2},
\end{equation}
and it is customary to define the significant wave height as $H_s = 4\eta_\text{rms}$, so that
\begin{equation}
H_s^2 = 16 \eta_\text{rms}^2 = 16 \overline{\eta^2}  = 16 \int_0^\infty S(f) {df}.
\end{equation}
The power flux is
\begin{equation}
J = \rho g \int_0^\infty v_g(f) S(f) {df} = \rho g \int_0^\infty \frac{g}{4\pi f} S(f) {df}  =\frac{\rho g^2}{4\pi} \int_0^\infty S(f) f^{-1} {df}.
\end{equation}
Finally, the energy period is defined as
\begin{equation}
T_e = \frac{\int_0^\infty S(f) f^{-1}{df}}{\int_0^\infty S(f) {df}},
\end{equation}
so that
\begin{equation}
J = \frac{\rho g^2}{64 \pi} T_e H_s^2.
\label{eq:powfluxirreg}
\end{equation}
Comparing \eqref{eq:powfluxmonochr} and \eqref{eq:powfluxirreg} we conclude that, in order to obtain the same power flux, it is possible to approximate a spectrum described by parameters $(H_s, T_e)$ by a monochromatic wave of height $H = H_s/\sqrt{2}$ and period $T = T_e$. We mention that if the peak period is considered instead of the energy period, that is, if the spectrum is defined by parameters $(H_s, T_p)$, then the power is not uniquely determined \cite{Guillou2020}.

\bibliography{refers}{}
\bibliographystyle{abbrv}

\end{document}